\documentclass[11pt]{article}

\usepackage{arxiv}

\usepackage{amsmath,amsfonts,latexsym,fullpage,graphicx,vmargin, mathrsfs}
\usepackage{pstricks}
\usepackage{comment}
\usepackage[inline]{enumitem}
\usepackage{cleveref}

\usepackage[ruled]{algorithm2e}
\usepackage{algpseudocode}

\RequirePackage[OT1]{fontenc}

\usepackage{amsmath}
\usepackage{amsfonts}
\usepackage{caption}
\usepackage{subcaption}
\usepackage{multirow}
\usepackage{tikz}
\usetikzlibrary{graphs}
\usetikzlibrary{cd}
\usepackage{bm}
\usepackage{url}
\usepackage{placeins}

\usepackage{float}

\RequirePackage{natbib}

\newtheorem{teo}{Theorem}
\newtheorem{prop}{Proposition}
\newtheorem{defi}{Definition}
\newtheorem{rmk}{Remark}
\newtheorem{lem}{Lemma}
\newtheorem{cor}{Corollary}

\newtheorem{assump}{Assumption}

\DeclareMathOperator*{\rank}{rank}

\DeclareMathOperator*{\Mapp}{Mapp}

\DeclareMathOperator*{\Top}{Top}

\DeclareMathOperator*{\Sets}{Sets}
\DeclareMathOperator*{\FSet}{FSet}
\DeclareMathOperator*{\Vect}{Vec}

\DeclareMathOperator*{\Tr}{Tr}

\usepackage{scalerel}
\DeclareMathOperator*{\Bigcdot}{\scalerel*{\cdot}{\bigodot}}

\newcommand{\X}{X_{\Bigcdot}}
\newcommand{\Y}{Y_{\Bigcdot}}

\newcommand{\T}{\pi_0(\X)}
\newcommand{\G}{\pi_0(\Y)}

\newcommand{\R}{\mathbb{R}}
\newcommand{\N}{\mathbb{N}}
\newcommand{\F}{\mathcal{F}}

\newcommand{\MT}{\mathcal{MT}}

\renewcommand{\S}{\mathcal{S}}

\newcommand{\virgolette}[1]{``#1''}

\firstpageno{1}

\begin{document}

\title{A Finitely Stable Edit Distance for Merge Trees}

\author{Matteo Pegoraro\thanks{Department of Mathematical Sciences, Aalborg University}}

%
\maketitle

\begin{abstract}
In this paper we define a novel edit distance for merge trees, which we argue to be suitable for a good range of applications. 
Relying also on some technical results contained in other works, we investigate its stability properties, which end up being analogous to the ones of the $1$-Wasserstein distance between persistence diagrams. In the appendix, we extensively compare our metric in relationship with other metrics appearing in the literature, with both theoretic and practical considerations and a simulation. 
\end{abstract}


\section{Introduction}\label{sec:intro}

Topological Data Analysis (TDA) is a particular set of techniques 
within the field of geometric data analysis which aim at including topological information 
into data analysis pipelines.
Topological information is usually understood in terms of generators of homology groups \citep{hatcher} with coefficients in some field. With \emph{persistent homology} these generators are extracted along a filtration of topological spaces to capture the shape of the initial datum, typically a function or a point cloud, at \virgolette{different resolutions} \citep{PH_survey}.
To proceed with the analysis, this ordered family of vector spaces is then represented with a topological summary. 
There are many different kinds of topological summaries such as persistence diagrams \citep{PD_1} (PDs), persistence images \citep{pers_img}, persistence silhouettes \citep{silhouettes}, and persistence landscapes \citep{landscapes}. Each of these summaries live in a space with different properties and purposes: for instance persistence diagrams are highly interpretable and live in a metric space, persistence landscapes are embedded in a linear space of functions, but the embedding is not closed under linear combinations, persistence images are instead obtained as vectors in $\mathbb{R}^n$, making them suitable for many Machine Learning techniques.

Along with the aforementioned summaries, there are also tree-shaped objects called merge trees. Merge trees arise naturally within the framework of TDA when dealing with zero dimensional homology groups, as they capture the merging structure of path connected components along a filtration of topological spaces. 
Originally, such objects stem out of Morse Theory \citep{milnor2016morse} as a topological summary related to 
 Reeb graphs \citep{reeb_1,reeb_2} and are frequently used for data visualization purposes \citep{merge_visual_1, merge_visual_2}.  
 Analogously, other different but related kinds of trees like hierarchical clustering dendrograms \citep{dendro_1} or phylogenetic trees \citep{phylo}  have also been used extensively in statistics and biology to infer information about a fixed set of labels. When considered as unlabelled objects, however, both clustering dendrograms and phylogenetic trees can be obtained as particular instances of merge trees.
The uprising of TDA has propelled works aiming at using merge trees and Reeb graphs as topological summaries in data analysis contexts and thus developing metrics and frameworks to analyze populations of such objects \citep{merge, merge_interl, merge_intrins, merge_frechet, merge_farlocca, merge_farlocca_2, di2016edit, bauer2020reeb}.

\subsection*{Previous Works on Merge Trees}

The works that have been dealing with the specific topic of merge trees can be divided into two groups: the first group is more focused on the definition of a suitable metric structure on merge trees and the second is more focused on the properties of merge trees and their relationships with PDs \citep{kanari2020trees,curry2021trees,beers2023fiber}.  
The first group, in turns, splits into a) works dealing with the interleaving distance between merge trees and other metrics with very strong theoretical properties
\citep{merge, merge_interl, interl_approx, merge_intrins, merge_frechet, curry2021decorated, cardona2021universal}, b) metrics which are more focused on the computational efficiency 
\citep{merge_farlocca, merge_farlocca_2, merge_wass, sridharamurthy2021comparative, wetzels2022deformation} at the cost of sacrificing stability properties and trying to mitigate the resulting problems with pre-processing and other computational solutions.
In particular, all these computationally-friendly metrics are framed as different variations stemming from the classical tree edit distance between unlabeled trees, which is obtained minimizing the edit distance for labeled trees \citep{tai1979tree, survey_ted} over all possible sets of labels. We further comment on this fact in the following section.

\subsection*{Main Contributions}

The aim of this manuscript is to propose a novel edit distance for merge trees and to investigate it's stability properties with a data analysis perspective. Plus, we extensively compare it with the other edit distances for merge trees appeared in literature. A very concise summary of this comparison can be found in \Cref{tab:edit}, which we also comment in the upcoming lines. 

\begin{itemize}
    \item The metric we propose has stability properties analogous to the ones of $1$-Wasserstein distance between PDs.
    In particular, the $1$-Wasserstein distance ($W_1$) and the bottleneck distance between PDs ($d_B$) enjoy the following relationship for every couple of diagrams $D_1,D_2$:
    \begin{equation}\label{eq:wass}
    d_B(D_1,D_2)\leq W_1(D_1,D_2)\leq (\rank(D_1)+\rank(D_2))d_B(D_1,D_2).
    \end{equation}    
    With  $\rank(D_i)$ being the cardinality of the diagrams $D_i$.
    Similarly, we prove that the edit distance that we define ($d_E$) and the interleaving distance between merge trees ($d_I$) satisfy:
    \begin{equation}\label{eq:interl}
    d_I(T_1,T_2)\leq d_E(T_1,T_2)\leq 2(\dim(T_1)+\dim(T_2))d_I(T_1,T_2).
    \end{equation}
        With  $\dim(T_i)$ being the number of edges in the merge tree $T_i$. We point out that the rightmost inequality is obtained in \cite{pegoraro2024functional}, which investigates some of the statistical properties of the metric we define here. While we prove the leftmost inequality in \Cref{prop:d_I<d_E}.
        
    In \Cref{sec:stability} we claim that \Cref{eq:wass} and \Cref{eq:interl} make $W_1$ and $d_E$ better suited for general data analysis purposes, compared to their \emph{universal} counterparts, that is, $d_B$ and $d_I$. \Cref{tab:edit} highlights that no other edit distance for merge trees is stable. We expand on that in \Cref{sec:edit_comparison} giving precise examples on how the metrics in \cite{merge_farlocca, merge_farlocca_2, merge_wass, sridharamurthy2021comparative, wetzels2022deformation} fail to be stable.
    
    \item Call $D_E$ the classical edit distance between unlabeled trees.
    \Cref{tab:costs} compares the computational costs of $D_E$ \citep{TED} and $d_E$ in terms of binary linear programming (BLP).      In particular, it shows that $d_E$ differs from $D_E$ only by two \emph{log} factors in the number of variables. thus, we argue that $d_E$ is relatively efficient. To the point that it has already been used in applications by other works \citep{pegoraro2024functional, cavinato2022imaging}.

\begin{table}[H]
\begin{tabular}{llll}
      & BLP Problems  & Variables                       & Constraints                   \\
$D_E(T,G)$ & $O(M\cdot N)$ & $O(M\cdot N)$                   & $O(\log_2(M)\cdot \log_2(N))$ \\
$d_E(T,G)$ & $O(M\cdot N)$ & $O(M\cdot\log_2(M)\cdot N\cdot\log_2(N))$ & $O(\log_2(M)\cdot \log_2(N))$
\end{tabular}
\caption{Computational costs of the classical edit distance for unlabeled trees \citep{TED} and our edit distance for merge trees. With $N = \dim T$ and $M = \dim G$.}
\label{tab:costs}
\end{table}

    We highlight that \cite{merge_farlocca, merge_farlocca_2, merge_wass, sridharamurthy2021comparative, wetzels2022deformation} all have a polynomial time algorithm/approximation. These algorithms are obtained in similar ways (see also \Cref{tab:edit}). \cite{selkow1977tree, jiang1995alignment, zhang1996constrained} contain some approximations of $D_E$, which can also be regarded as a metrics on their own. Such  \virgolette{approximated} edit distances have polynomial time complexity.  \cite{merge_farlocca, merge_farlocca_2, merge_wass, sridharamurthy2021comparative, wetzels2022deformation} frame their metrics as edit distances and then resort to these variations to obtain their feasible algorithms. \Cref{tab:edit} and \Cref{rmk:constr} highlight that, for instance, a constrained variation \citep{zhang1996constrained} is easily allowed also for our metric, but we still haven't worked on an algorithm to compute it. 
    Since incorporating our modifications into $D_E$ to obtain $d_E$ does not significantly increase its computational cost, we strongly believe that also the computational complexity of the constrained edit distance is only marginally affected by such modifications. And so, in particular, a constrained version of $d_E$ should admit a polynomial time algorithm.
\end{itemize}

\begin{table}[]
\begin{tabular}{llllll}
      & Metric & Stability                    & Cost       & Approx. Def. & Approx. Alg. \\
$d_E$ & Yes    & $d_I\leq d_E \leq 2(N+M)d_I$ & $\geq D_E$ & Yes          & $-$          \\
Wass 
  & Yes    & No                           & $\sim D_E$ & Yes          & Yes          \\
BDI 
  & No     & No                           & $\geq D_E$ & Yes          & Yes          \\
DB    & Yes    & No                           & $\geq D_E$ & Yes          & Yes          \\
Local 
 & Yes    & No                           & $\sim D_E$ & Yes          & Yes         
\end{tabular}
\caption{A summary of Edit Distances for Merge Trees. As in the main text, $D_E$ refers to the classical edit distance between unlabeled trees. \virgolette{Wass} stands for the metrics in \citep{merge_farlocca, merge_wass}; 
\virgolette{BDI} is for the one in \cite{merge_farlocca_2}; \virgolette{DB} stands for \cite{wetzels2022deformation}, and \virgolette{Local} for \citep{sridharamurthy2021comparative}.
Column \virgolette{Metric} states which distances are also metrics. Column \virgolette{Stability} underlines that $d_E$ is the only edit distance satisfying stability properties. Column \virgolette{Cost} compares the cost of the other edit distances with $D_E$:  $\sim D_E$ stands for the cost being equivalent, up to nuances, to the one of $D_E$; while $\geq D_E$ stands for the cost being higher.  \virgolette{Approx. Def.} and \virgolette{Approx. Alg.} tell for which edit distances an approximated version was given, along with a polynomial time algorithm. In particular: \cite{merge_farlocca, merge_wass} are obtained employing the metric $D_E$ (with some modifications on the costs) on trees built from the original merge trees, which encode the hierarchical properties of persistence pairs, and use \cite{zhang1996constrained} to find a polynomial time algorithm. \cite{merge_farlocca_2} tries to overcome the cons that come with fixing such hierarchical structure, but ends up losing the triangular inequality. The tractable algorithm is again inspired by \cite{zhang1996constrained}. \cite{wetzels2022deformation} shares some similarities with our approach, see \Cref{sec:edit_comparison}, but fails to achieve stability, and obtains a polynomial algorithm via \cite{selkow1977tree}. \cite{sridharamurthy2021comparative} uses \cite{zhang1996constrained} to compare local structures in merge trees via a modified constrained edit distance.}
\label{tab:edit}
\end{table}

To summarize, we aim at providing a metric which is interpretable and with good discriminative power due to its stability properties, computable at least in small data regimes and which has a potential polynomial time approximation to be studied in future works. 
 
%

\subsection*{Outline}

This paper is organized as follows. 
\Cref{sec:preliminary_defi}, and \Cref{sec:merge_and_edit} contain the preliminary definitions needed in the paper, which are collected from previous works in the field.
In the latter one, in particular, we review the definition of an edit distance between weighted trees, which we use in \Cref{sec:tree_edit} to define a metric structure for merge trees and to investigate its stability properties.
In \Cref{sec:conclusions_geom} we draw some conclusions.

The appendix reports definitions taken from literature (\Cref{sec:abstract_merge}, \Cref{sec:interl}, \Cref{sec:wass}) which can help the reader in navigating the manuscript. Plus, it features our comments on the potential polynomial time approximation we foresee for our metric (\Cref{rmk:constr}),
and the detailed comparison between the metrics for merge trees we mention in the introduction (\Cref{sec:edit_comparison}).
 \Cref{sec:proofs_geom} contains the proofs of the results in the manuscript.

\section{Preliminary Definitions - Merge Trees}
\label{sec:preliminary_defi}

We introduce merge trees coherently with most of the scientific literature dealing with such topics \citep{merge_intrins, merge_farlocca, merge_farlocca_2, merge_wass}: a combinatorial object, which we call tree structure, with a monotone increasing function defined on its vertices.

\begin{defi}
A tree structure $T$ is given by a set of vertices $V_T$ and a set of edges $E_T\subset V_T\times V_T$ which form a connected rooted acyclic graph.  We indicate the root of the tree with  $r_T$. We say that \(T\) is finite if \(V_T\) is finite. The order of a vertex $v\in V_T$ is the number of edges which have that vertex as one of the extremes, and is called $ord_T(v)$. 
Any vertex with an edge connecting it to the root is its child and the root is its father: this is the first step of a recursion which defines the father and children relationship for all vertices in \(V_T.\)
The vertices with no children are called leaves  or taxa and are collected in the set $L_T$. The relation $child < father $ generates a partial order on $V_T$. The edges in $E_T$ are identified in the form of ordered couples $(a,b)$ with $a<b$.
A subtree of a vertex $v$, called $sub_T(v)$, is the tree structure whose set of vertices is $\{x \in V_T \mid x\leq v\}$. 
\end{defi}

Given a tree structure $T$, identifying an edge $(v,v')$ with its lower vertex $v$, gives a bijection between $V_T-\{r_T\}$ and $E_T$, that is $E_T\cong V_T-\{r_T\}$ as sets. 
Given this bijection, we often use $E_T$ to indicate the vertices $v\in V_T-\{r_T\}$, and viceversa, to simplify the notation.

Now we want to identify merge trees independently of their vertex set.

\begin{defi}
Two tree structures $T$ and $T'$ are isomorphic if exists a bijection $\eta:V_T\rightarrow V_{T'}$ that induces a bijection between the edges sets $E_T$ and $E_{T'}$: $(a,b)\mapsto (\eta(a),\eta(b))$. Such $\eta$ is an isomorphism of tree structures.
\end{defi}

Then, we can give the definition of a merge tree.

\begin{defi}\label{defi:merge_TDA}
A merge tree is a finite tree structure $T$ with a monotone increasing height function $h_T:V_T\rightarrow \mathbb{R}\cup\{+\infty\}$ and such that 1) $ord_T(r_T)=1$ 2) $h_T(r_T)=+\infty$ 3) $h_T(v)\in\mathbb{R}$ for every $v<r_T$.
Two merge trees $(T,h_T)$ and $(T',h_{T'})$ are isomorphic if $T$ and $T'$ are isomorphic as tree structures and the isomorphism $\eta:V_T\rightarrow V_{T'}$ is such that $h_T = h_{T'} \circ \eta$. Such $\eta$ is an isomorphism of merge trees. We use the notation $(T,h_T)\cong (T',h_{T'})$. The set of all merge trees up to isomorphism is called 
$\mathcal{MT}$.
\end{defi}

With some slight abuse of notation we set $\max h_T =\max_{v\in V_{T}\mid v<r_T}h_{T}(v)$ and $\arg \max h_T =\max\{v\in V_{T}\mid v<r_T\}$. Note that, given $(T,h_T)$ merge tree, there is only one edge of the form $(v,r_T)$ and we have $v=\arg\max h_T$.

\begin{defi}\label{defi:ghosting}
Given a tree structure $T$, we can eliminate an  order two vertex, connecting the two adjacent edges which arrive and depart from it. 
Suppose we have two edges $e=(v_1,v_2)$ and $e'=(v_2,v_3)$, with $v_1<v_2<v_3$. And suppose $v_2$ is of order two. Then, we can remove $v_2$ and merge $e$ and $e'$ into a new edge $e''=(v_1,v_3)$.
This operation is called the \emph{ghosting} of the vertex. Its inverse transformation is called the \emph{splitting} of an edge.

Consider a merge tree $(T,h_T)$ and obtain $T'$ by ghosting a vertex of $T$. Then $V_{T'}\subset V_T$ and thus we can define $h_{T'}:=h_{T\mid V_{T'}}$.
\end{defi}

Now we can state the following definition.

\begin{defi}
Merge trees are equal up to order $2$ vertices if they become isomorphic after applying a finite number of ghostings or splittings. We write $(T,h_T)\cong_2(T',h_{T'})$. The set of all merge trees up to order $2$ vertices is called 
$\mathcal{MT}/\sim_2$.
\end{defi} 

Merge trees are usually employed go give a combinatorial representation of \emph{persistent sets} \citep{patel2018generalized,curry2021decorated} which are generated via usual TDA pipelines, applying the functor $\pi_0$ on filtrations of topological spaces. Going through such details is not essential for defining and understanding the metric $d_E$, however it is necessary to introduce the interleaving distance between merge trees and discuss stability properties. Thus, we report them in \Cref{sec:interl}.

\section{Preliminary Definitions - Weighted Trees Edit Distance}\label{sec:merge_and_edit}

In this section we introduce the last pieces of notation we need to define the merge tree edit distance.
In \Cref{sec:weights_edits} we report how 
\cite{pegoraro2023edit} builds a metric for \emph{weighted trees} and in \Cref{sec:mappings} we recall from \cite{pegoraro2023edit} the definition of \emph{mapping}, a combinatorial object which will be used throughout the manuscript.

\subsection{Weighted Trees and Edits}
\label{sec:weights_edits}

\begin{defi}
A tree structure $T$ with a weight function $w_T:E_T\rightarrow \mathbb{R}_{>0}$ is called weighted tree. Isomorphisms of weighted trees are defined as in \Cref{defi:merge_TDA}. The set of all weighted trees up to isomorphism is called $(\mathcal{T},\mathbb{R}_{\geq 0})$. 
\end{defi}

Given a weighted tree $(T,w_T)$ we can modify its edges $E_T \cong V_T-\{r_T\}$ with a sequence of the following edit operations:

\begin{itemize}
\item we call \emph{shrinking} of a vertex/edge 
a change of the weight function. The new weight function must be equal to the previous one on all edges, a part from the \virgolette{shrunk} one, whose weight can become bigger or smaller. In other words,  for an edge $e$, this means changing the value $w_T(e)$ with another value in $\mathbb{R}_{>0}$. 
\item A \emph{deletion} is an edit with which a vertex/edge is deleted from the dendrogram. Consider an edge $(v_1,v_2)$. The result of deleting $v_1$ is a new tree structure, with the same vertices a part from $v_1$ (the smaller one), and with the father of the deleted vertex which gains all of its children.
The inverse of the deletion is the \emph{insertion} of an edge along with its lower vertex. 
We can insert an edge at a vertex \(v\) specifying the name of the new child of \(v\), the children of the newly added vertex (that can be either none, or any portion of the children of \(v\)), and the value of the weight function on the new edge. 
\item Lastly, we can remove or add order two vertices via the \emph{ghosting} and \emph{splitting} edits, which have already been defined in \Cref{defi:ghosting}. The weight function is obtained by summation of the merged edges in case of ghostings, and must be arbitrarily defined in case of splittings. With the condition that splitting an edge and then ghosting the order two vertex obtained must restore the original edge weight.
\end{itemize}

\begin{figure}
	\begin{subfigure}[c]{0.49\textwidth}
    	\centering
    	\includegraphics[width = \textwidth]{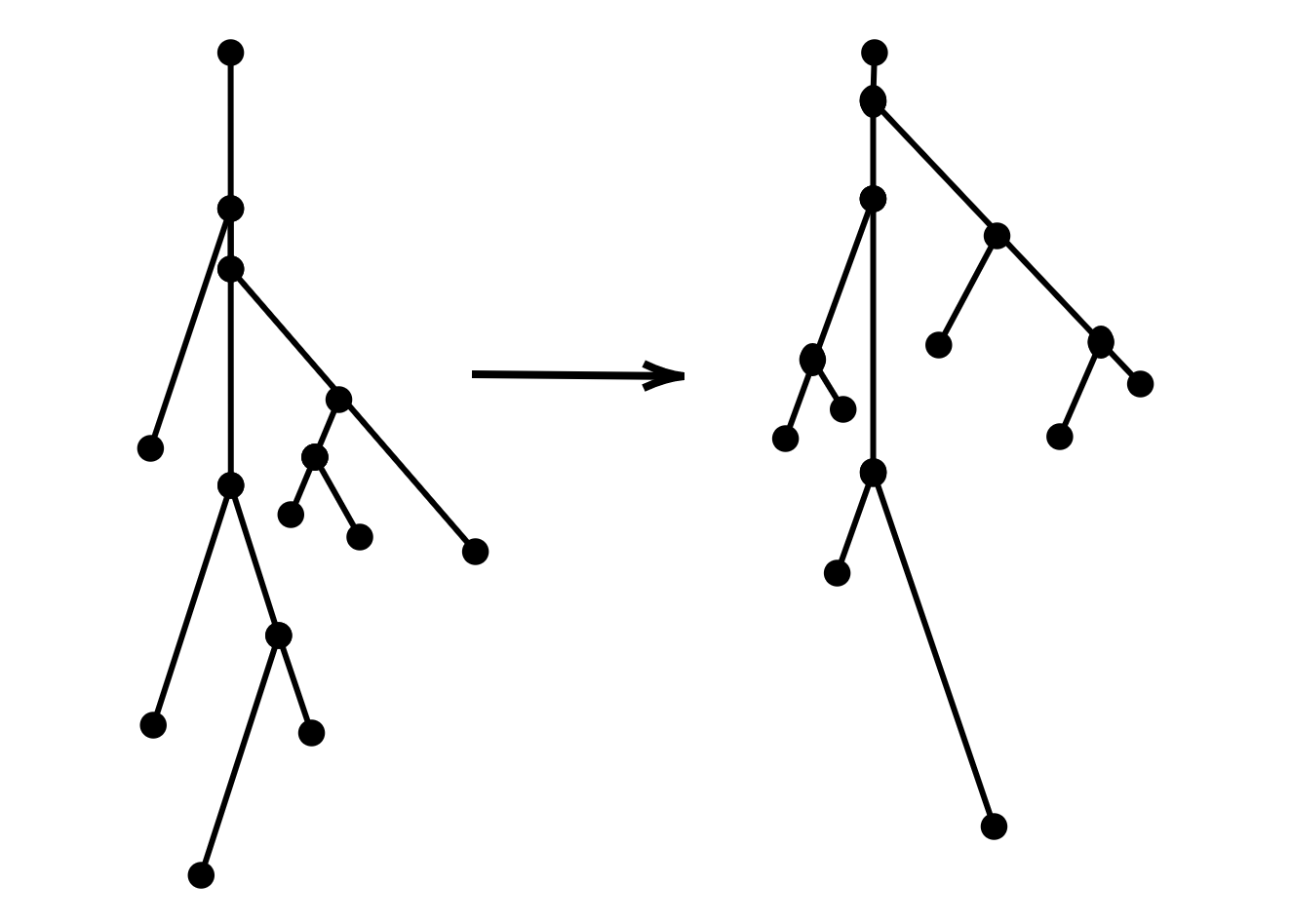}
    	\caption{Starting and target weighted trees.}
    	\label{fig:recap}
    \end{subfigure}
   	\begin{subfigure}[c]{0.49\textwidth}
    	\centering
    	\includegraphics[width = \textwidth]{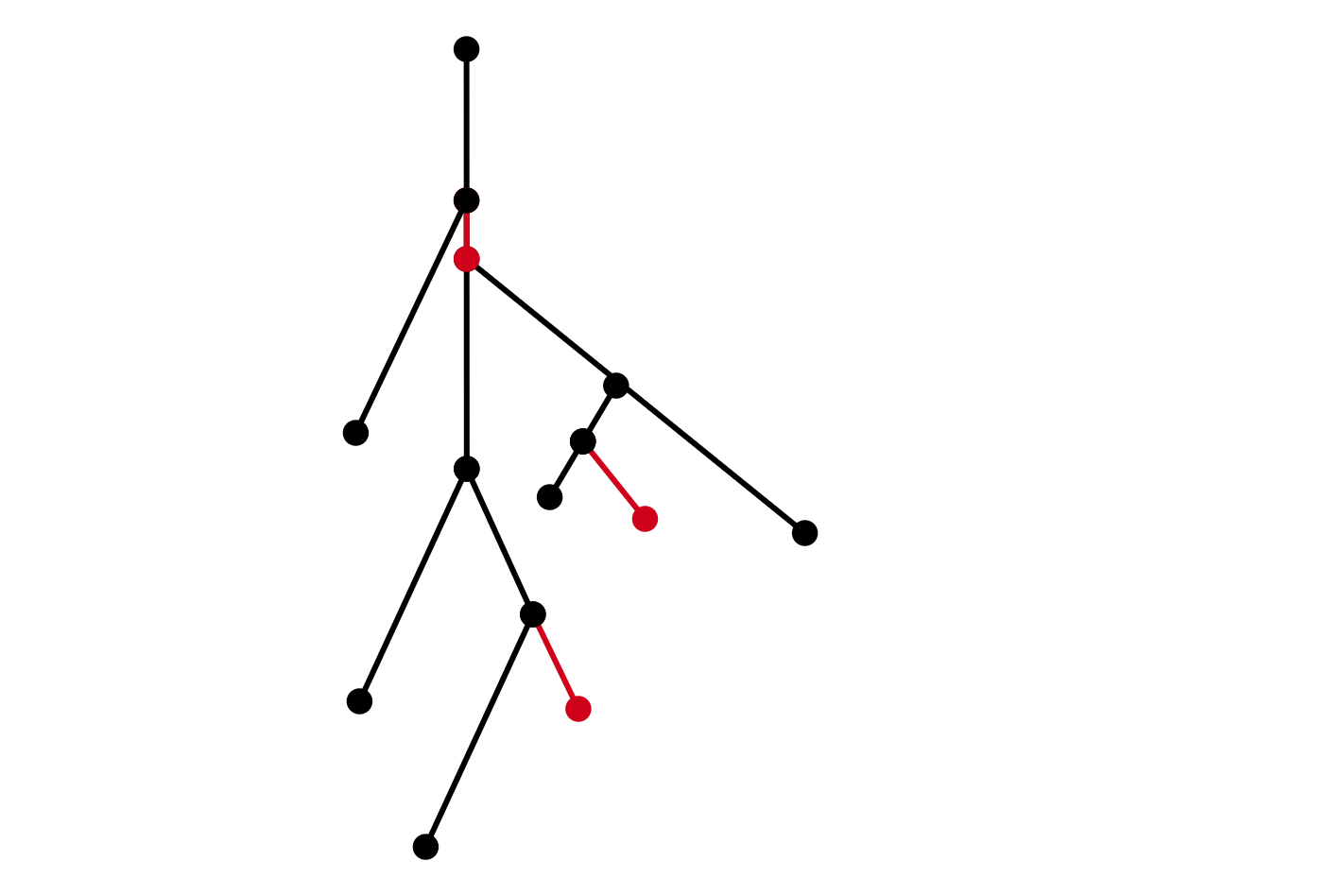}
    	\caption{Deletions in red.}
    	\label{fig:deletions}
    \end{subfigure}

	\begin{subfigure}[c]{0.49\textwidth}
		\centering
		\includegraphics[width = \textwidth]{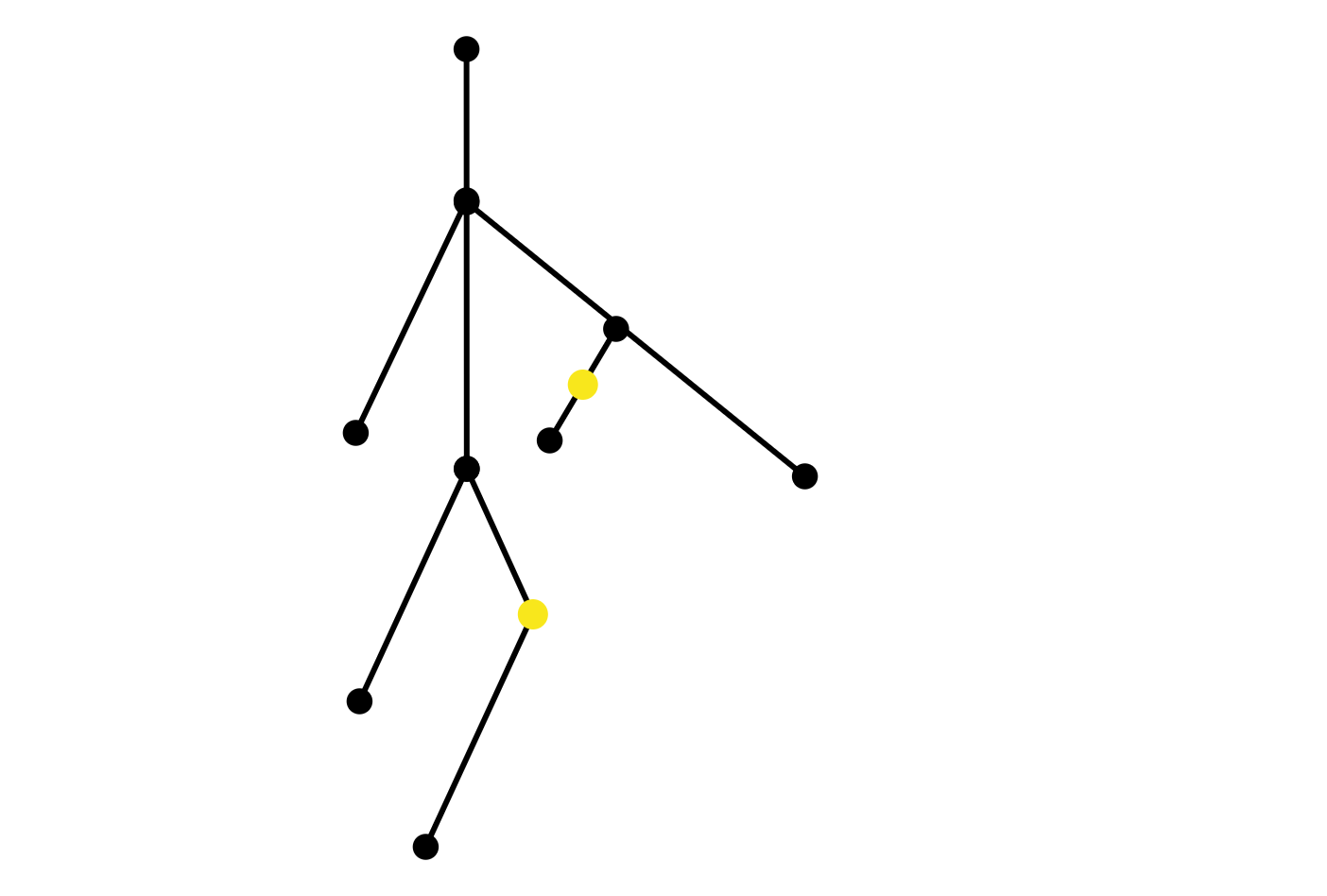}
		\caption{Ghostings in yellow.}
		\label{fig:ghostings}
	\end{subfigure}
	\begin{subfigure}[c]{0.49\textwidth}
    	\centering
    	\includegraphics[width = \textwidth]{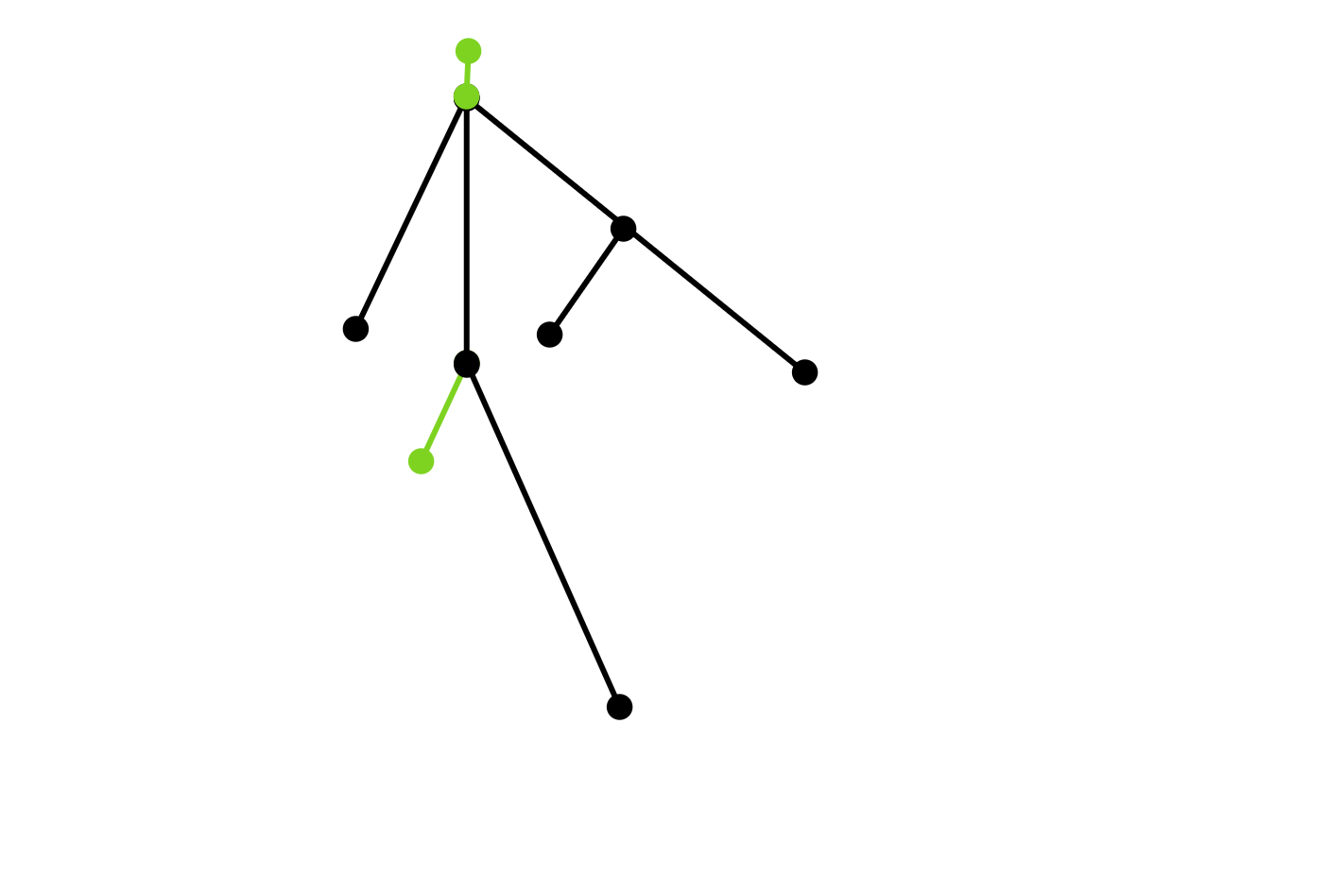}
    	\caption{Shrinkings in green.}
    	\label{fig:shrinkings}
    \end{subfigure}

    \begin{subfigure}[c]{0.49\textwidth}
    	\centering
    	\includegraphics[width = \textwidth]{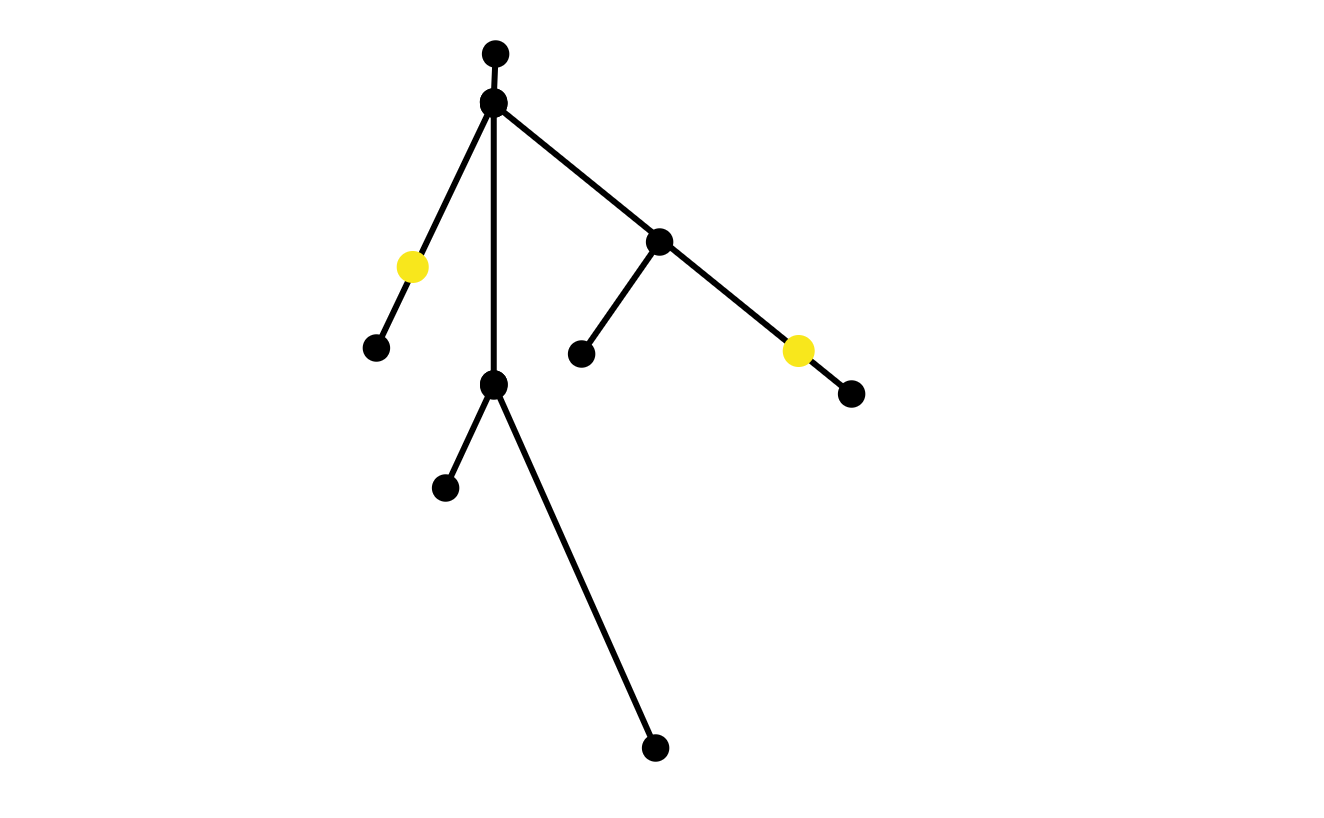}
    	\caption{Splittings in yellow.}
    	\label{fig:splittings}
    \end{subfigure}
    \begin{subfigure}[c]{0.49\textwidth}
    	\centering
    	\includegraphics[width = \textwidth]{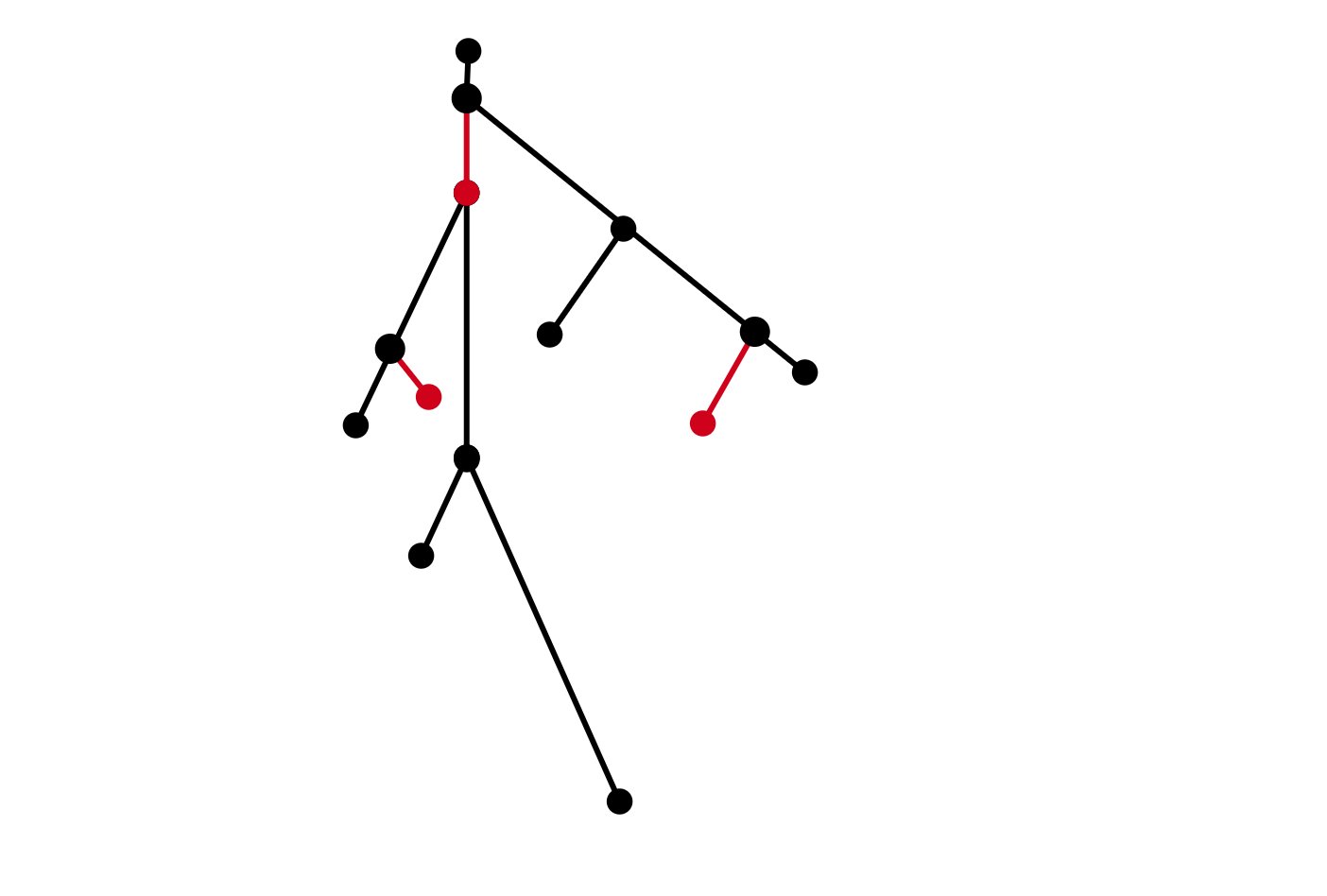}
    	\caption{Insertions in red.}
    	\label{fig:insertions}
    \end{subfigure}
    \hfill
\caption{(b)$\rightarrow$(e) form an edit path made from the left weighted tree in \Cref{fig:recap} to the right one. Each time the edges involved in the editing are highlighted with different colors. In the following plot such vertices return black. This edit path can clearly be represented with a mapping - \Cref{sec:mappings} - made by 
couples $(v,"D")$ for all the red vertices in \Cref{fig:deletions}, $(v,"G")$ for all the yellow vertices in \Cref{fig:ghostings}, $(v,w)$ for all the vertices associated via the green color in \Cref{fig:shrinkings}, $("G",w)$ for all the yellow vertices in \Cref{fig:ghostings} and $("D",w)$ for all the red vertices in \Cref{fig:insertions}. }
\label{fig:edits_TDA}
\end{figure}

The set of all weighted trees  considered up to ghostings and splittings is called $(\mathcal{T}_2,\mathbb{R}_{\geq 0})$.

A weighted tree $T$ can be edited to obtain another weighted tree, on which one can apply a new edit to obtain a third tree and so on. 
Any finite composition of edits is called \emph{edit path}. See \Cref{fig:edits_TDA} for an example of an edit path. The set of finite edit paths between $T$ and $T'$ is called $\Gamma(T,T')$.

The cost of the edit operations is defined as follows:
\begin{itemize}
\item the cost of shrinking an edge is equal to the absolute value of the difference of the two weights;
\item for any deletion/insertion, the cost is equal to the weight of the edge deleted/inserted;
\item the cost of ghostings and splittings is zero.
\end{itemize}

The cost of an edit path is the sum of the costs of its edit operations.
Putting all the pieces together, the edit distance $d_E$ between weighted trees is defined  as $d_E(T,T')=\inf_{\gamma\in\Gamma(T,T')} cost(\gamma)$,
where $\Gamma(T,T')$ indicates the set of edit paths which start in $T$ and end in $T'$.
In \cite{pegoraro2023edit} it is proved that $d_E$ is a metric on the space
of weighted trees considered up to order two vertices.

\subsection{Mappings}\label{sec:mappings}
Given an edit path between two weighted trees, its cost 
is often invariant up to many permutations of the edits. 
To better work in such environment, we start considering paths up to some permutation of the edits. Objects called mappings, as defined in \cite{pegoraro2023edit}, Section \virgolette{Trees}, help us in doing this, as well as making the metric $d_E$ more tractable. For this reason now we report their definition.    
As in \cite{pegoraro2023edit}, $"D"$ and $"G"$ are be used to indicate \virgolette{deletion} and \virgolette{ghosting}.

A \emph{mapping} between $T$ and $T'$ is a set 
$M\subset (E_T \cup \{"D","G"\})\times (E_{T'} \cup \{"D","G"\})$ 
satisfying:

\begin{itemize}
\item[(M1)] consider the projection of the Cartesian product $(E_T \cup \{"D","G"\})\times (E_{T'} \cup \{"D","G"\})\rightarrow (E_{T} \cup \{"D","G"\})$; we can restrict this map to $M$ obtaining $\pi_T:M\rightarrow (E_T \cup \{"D","G"\})$. The maps $\pi_T$ and $\pi_{T'}$ are surjective on $E_T\subset (E_T \cup \{"D","G"\})$ and 
$E_{T'}\subset (E_{T'} \cup \{"D","G"\})$ respectively;
\item[(M2)]$\pi_T$ and $\pi_{T'}$ are injective;
\item[(M3)] $M\cap (V_T\times V_{T'})$ is such that, given $(a,b)$ and $(c,d)\in M\cap (V_T\times V_{T'})$, $a>c$, if and only if $b>d$;
\item[(M4)] if $(a,"G")\in M$ (or analogously $("G",a)$), 
then after applying all deletions of the form $(v,"D")\in M$, the vertex $a$ becomes an order $2$ vertex. In other words: let $child(a)=\{b_1,..,b_n\}$. Then there is one and only one $i$ such that for all $j \neq i$, for all $v \in V_{sub(b_j)}$, we have $(v,"D")\in M$; and there is one and only one $c$ such that $c=\max\{x<b_i\mid (x,y)\in M$ for any $y \in V_{T'}\}$. 
\end{itemize}  

We call $\Mapp(T,T')$ the set of all mappings between $T$ and $T'$. We may refer to edges which appear in the couples in $M\cap (V_T\times V_{T'})$ as the \emph{coupled} or \emph{matched} edges/vertices.

As in \cite{pegoraro2023edit}, we use the properties of $M\in \Mapp(T,T')$ to parametrize a set of edit paths starting from $T$ and ending in $T'$, which are collected under the name $\gamma_M$.

\begin{itemize}
\item $\gamma_{d}^T$ a path made by the deletions to be done in $T$, that is, the couples $(v,"D")$, executed in any order. So we obtain $T^M_d=\gamma_{d}^T(T)$, which, instead, is well defined and not depending on the order of the deletions.
\item One then proceeds with ghosting all the vertices $(v,"G")$ in $M$, in any order, getting a path $\gamma^T_g$ and the dendrogram $T_M:= \gamma^T_g\circ \gamma^T_d (T)$.
\item Since all the remaining points in $M$ are coupled, the two dendrograms $T'_M$ (defined in the same way as $T_M$, but starting from $T'$) and $T_M$ must be isomorphic as tree structures. This is guaranteed by the properties of $M$. So one can shrink $T_M$ onto $T'_M$, and the composition of the shrinkings, executed in any order is an edit path $\gamma_s^T$.
\end{itemize}

By construction $\gamma_s^T\circ \gamma_g^T\circ \gamma_d^T(T)=T'_M$,
and $(\gamma_d^{T'})^{-1}\circ (\gamma_g^{T'})^{-1}\circ \gamma_s^T\circ \gamma_g^T\circ \gamma_d^T (T)=T'.
$.
Where the inverse of an edit path is thought as the composition of the inverses of the single edit operations, taken in the inverse order.

Lastly, we call $\gamma_M$ the set of all possible edit paths of the form $(\gamma_d^{T'})^{-1}\circ (\gamma_g^{T'})^{-1}\circ \gamma_s^T\circ \gamma_g^T\circ \gamma_d^T$, 
obtained by changing the order in which the edit operations are executed inside $\gamma_d$, $\gamma_g$ and $\gamma_s$.
Even if $\gamma_M$ is a set of paths, its cost is well defined:

\[
cost(M):= cost(\gamma_M)=cost(\gamma_d^T)+cost(\gamma_s^T)+cost(\gamma_d^{T'}).
\]

See \Cref{fig:edits_TDA} for an example of a mapping between weighted trees.
We conclude this section by recalling that \cite{pegoraro2023edit}
proves that given two weighted trees $T$ and $T'$, for every finite edit path $\gamma\in \Gamma(T,T')$, there exists a mapping $M \in \Mapp(T,T')$ such that $cost(M)\leq cost(\gamma)$.

Lastly, consider $M_2(T,T')\subset \Mapp(T,T')$ defined as follows.

\begin{defi}[\cite{pegoraro2023edit}]\label{defi:M_2_mapp}
A mapping $M\in \Mapp(T,T')$ has maximal ghostings if $(v,"G")\in M$ if and only if $v$ is of order $2$ after the deletions in $T$ and, similarly $("G",w)\in M$ if and only if $w$ is of order $2$ after the deletions in $T'$.

A mapping $M\in \Mapp(T,T')$ has minimal deletions if $(v,"D")\in M$ only if neither $v$ nor $father(v)$ are of order $2$ after applying all the other deletions in $T$ and, similarly, $("D",w)\in M$ only if $w$ is not of order $2$ after applying all the other deletions in $T'$.

We collect all mappings with maximal ghostings and minimal deletions in the set $M_2(T,T')$.
\end{defi}

\begin{lem}[\cite{pegoraro2023edit}]\label{lemma:M_2}
\[
\min\{cost(M)\mid  M\in \Mapp(T,T') \}=\min\{cost(M)\mid M\in M_2(T,T')\}
\]
\end{lem}

\section{Merge Trees Edit Distance}
\label{sec:tree_edit}

In this section we finally exploit the notation established in the previous sections to obtain a (pseudo) metric for merge trees.

\subsection{Truncation Operators}
\label{sec:truncation}

First we need to bridge between merge trees and weighted trees, in order to induce a metric on merge trees by means of the edit distance defined in \Cref{sec:merge_and_edit}. The general idea is that we want to truncate the edge going to infinity of a merge tree to obtain a weighted tree. 

\begin{figure}
	\begin{subfigure}[c]{0.49\textwidth}
    	\centering
    	\includegraphics[width = \textwidth]{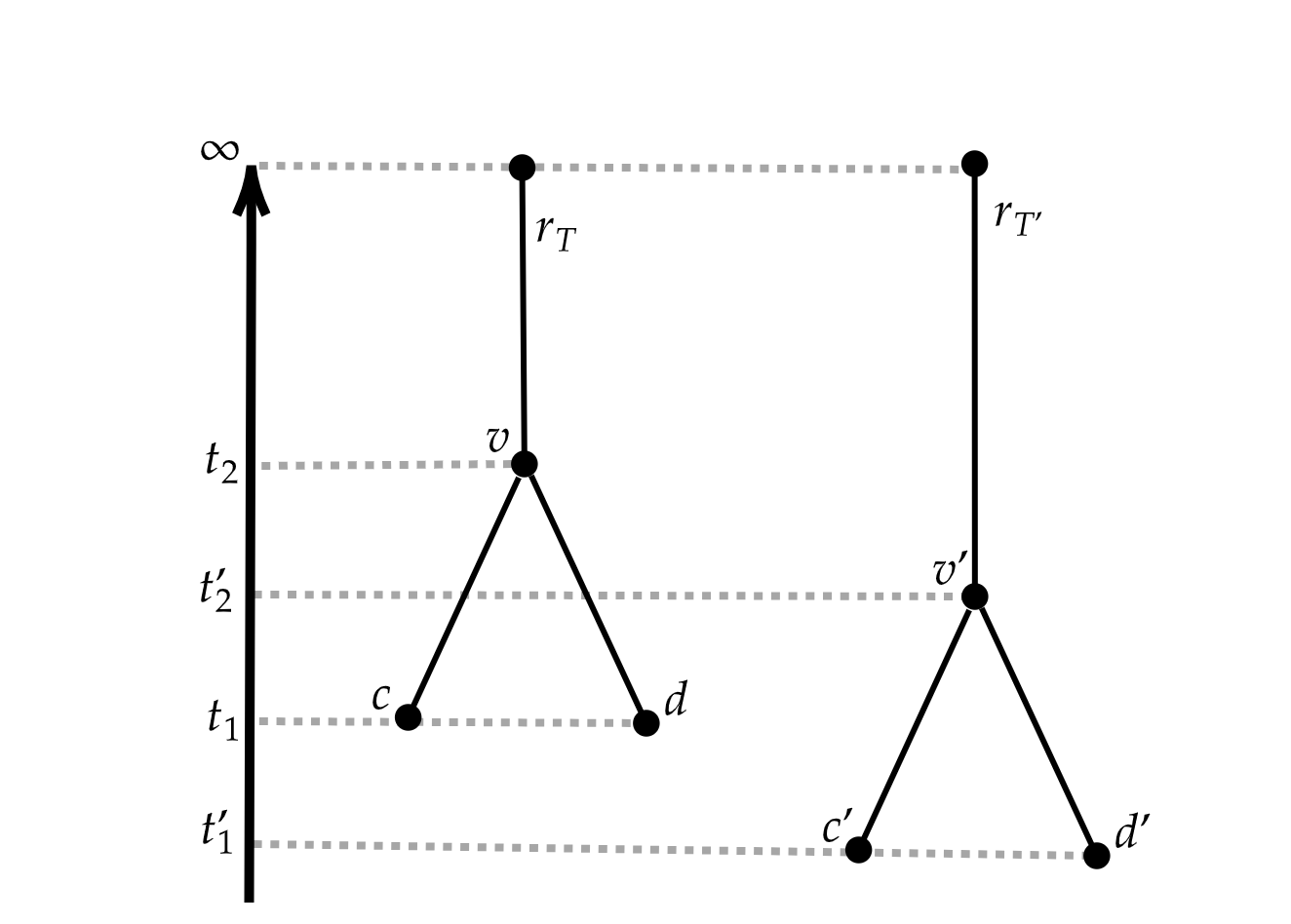}
    	\caption{Two merge trees $T$ and $T'$.}
    	\label{fig:truncate_1}
    \end{subfigure}
   	\begin{subfigure}[c]{0.49\textwidth}
    	\centering
    	\includegraphics[width = \textwidth]{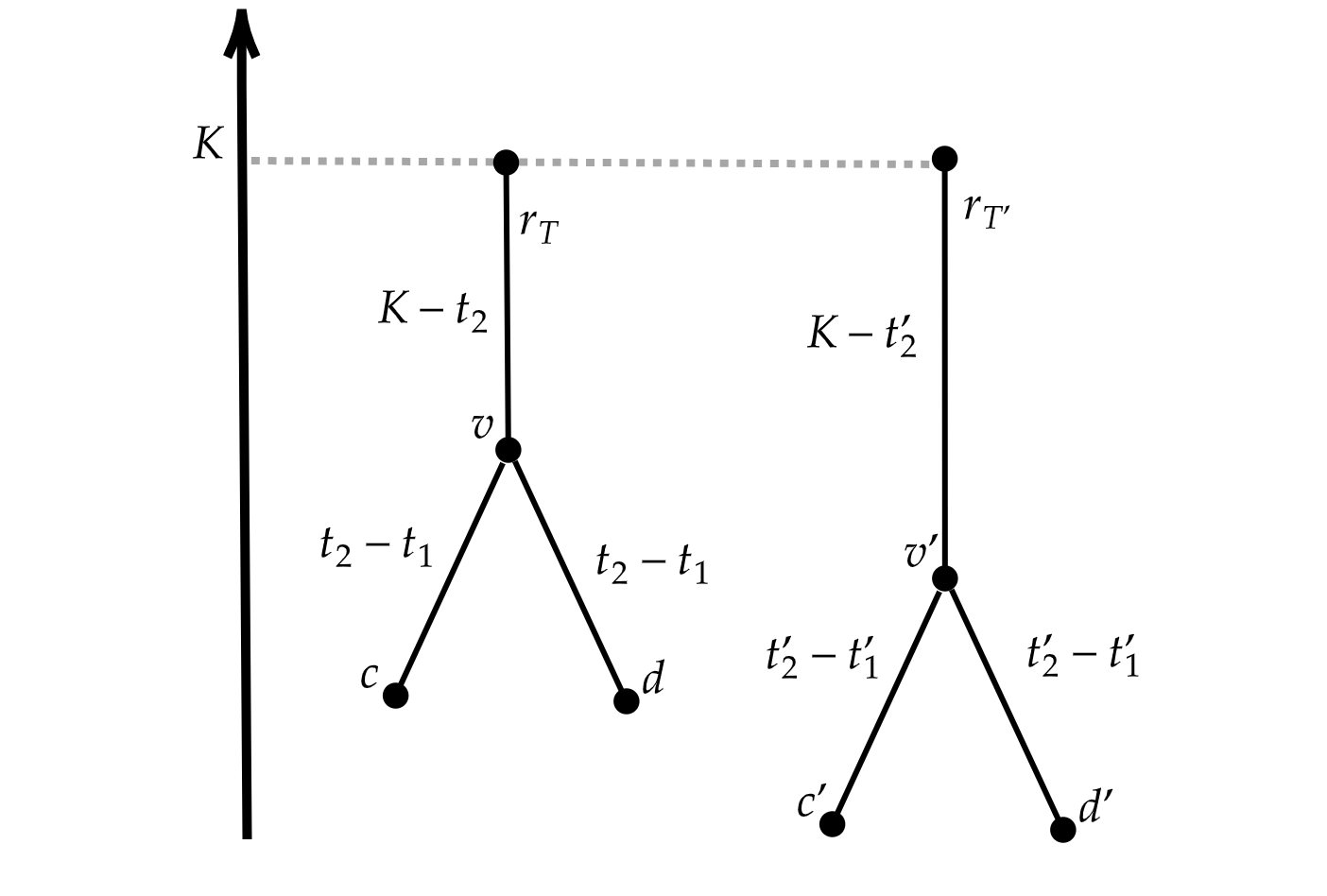}
    	\caption{A graphical representation of the weighted trees $\Tr_K(T)$ and $\Tr_K(T')$. Roughly speaking, the merge trees $T$ and $T'$ are recovered sending the two roots back to infinity.}
    	\label{fig:truncate_2}
    \end{subfigure}

\caption{A graphical representation of the truncation operator. }
\label{fig:truncation}
\end{figure}

Starting from a merge tree $(T,h_T)$ it is quite natural to turn the height function $h_T$ into a weight function $w_T$ via the rule $w_T((v,v')):=h_T(v')-h_T(v)$. The monotonicity of $h_T$ guarantees that $w_T$ take values in $\mathbb{R}_{> 0}$. However, we clearly have an issue with the edge $(v,r_T)$ as $h_T(r_T)=+\infty$. To solve this issue we need some novel tools - see \Cref{fig:truncation}. First consider the set of merge trees $\mathcal{MT}$ and build the subset $\mathcal{MT}_K :=\{(T,h_T)\in\mathcal{MT}\mid \max h_T\leq K\}$, for some $K\in\mathbb{R}$. Then define the truncation operator at height $K$ as follows:

\begin{equation}
\begin{aligned}
\textstyle\Tr_K:&\mathcal{MT}_K\xrightarrow{\hspace*{2cm}} (\mathcal{T},\mathbb{R}_{\geq 0})
\\
&(T,h_T)\mapsto (T,h_K)\mapsto (T,w_T)
\end{aligned}
\end{equation}

with $h_K(v)=h_T(v)$ if $v<r_T$ and $h_K(r_T)=K$. Then we set $w_T((v,v'))=h_K(v')-h_K(v)$. To avoid 
$w_T((v,r_T))=0$, if $\max h_T= K$, we take $(T,h_T)\mapsto (T',h_{T'})\mapsto (T',w_{T'})$ with $T'$ obtained from $T$ via the removal of $r_T$ from $V_T$ and $(v,r_T)$ from $E_T$. The map $h_{T'}$ is $h_{T\mid V_{T'}}$. Then $w_{T'}((v,v')):=h_{T'}(v')-h_{T'}(v)$.

In other words, with $\Tr_K$ we are fixing some height $K$, truncating the edge $(\arg\max h_T,r_T)$ at height $K$, and then obtaining a positively weighted tree, as in \Cref{fig:truncation}. To go back with $(\Tr_K)^{-1}$ we \virgolette{hang} a weighted tree at height $K$ and extend the edge $(v,r_T)$ to $+\infty$. We formally state these ideas in the following proposition which we state without proof.

\begin{prop}\label{rmk:inverting_repr}
The operator $\Tr_K:\mathcal{MT}_K\rightarrow (\mathcal{T},\mathbb{R}_{\geq 0})$ can be inverted via $(T,w_T)\mapsto \Tr_K^{-1}((T,w_T))=(T',h_{T'})$ with the following notation: the tree structure $T'$ is obtained from $T$ via adding $r_{T'}$ to $V_T$ and $(r_T,r_{T'})$ to $E_T$ and ghosting $r_T$ if it becomes an order $2$ vertex. Then we have $h_{T'}(r_T)=K$ (if it is not ghosted, i.e. if it is of order greater than $2$), $h_{T'}(r_{T'})=+\infty$ and, recursively, for $(v,v')\in E_{T}$, $h_{T'}(v)=h_{T'}(v')-w_T((v,v')) $. Clearly $\Tr_K^{-1}(\Tr_K((T,h_T)))\cong (T,h_T)$. Thus $\mathcal{MT}_K\cong (\mathcal{T},\mathbb{R}_{\geq 0})$ as sets, for every $K\in\mathbb{R}$.  Moreover $\Tr_K((T,h_T))\sim_2 \Tr_K((T',h_{T'}))$ if and only if 
$(T,h_{T})\sim_2 (T',h_{T'})$.
\end{prop}

\subsection{Edit Distance For Merge Trees}
\label{sec:MTED}

Consider $(T,h_T),(T',h_{T'})\in\mathcal{MT}$ and select $K$ such that $(T,h_T),(T',h_{T'})\in\mathcal{MT}_K$. Let 
$(G,w_G)=\Tr_K((T,h_T))$ and   
$(G',w_{G'})=\Tr_K((T',h_{T'}))$.
Lastly, set:
\[
d_K((T,h_T),(T',h_{T'})):= d_E((G,w_G),(G',w_{G'})).
\]

Despite being a promising and natural pseudo metric, $d_K$ is not defined on the whole $\mathcal{MT}$ and, on top of that, it is not clear if its value depends on the $K$ we choose. The following result solves these issues and proves that we can use $d_K$ to define a metric on
$\mathcal{MT}$.

\begin{prop}[Truncation]\label{prop:truncation}
Take $(T,w_T)$ and  $(T',w_{T'})$ weighted trees.
Suppose $r_T$ and $r_{T'}$ are of order $1$ and take two splitting operations: $\{(v,r_T)\}\rightarrow \{(v,v'),(v',r_T)\} $ and $\{(w,r_{T'})\}\rightarrow \{(w,w'),(w',r_{T'})\} $, with which we obtain the weighted trees $(G,w_{G})$ and $(G',w_{G'})$. Suppose that $w_{G}((v',r_{T}))=w_{G'}((w',r_{G}))$. Then $d_E(T,T')=d_E(sub_{G}(v'),sub_{G'}(w'))$.
\end{prop}

We now exploit \Cref{prop:truncation} to define the edit distance between merge trees.

\begin{teo}[Merge Tree Edit Distance]\label{teo:d_K}
Given two merge trees $(T,h_T),(T',h_{T'})$
such that $(T,h_T),(T',h_{T'})\in\mathcal{MT}_K$ and $(T,h_T),(T',h_{T'})\in\mathcal{MT}_{K'}$, then $d_{K}(T,T')=d_{K'}(T,T')$. 

Thus, for any couple of merge trees in $\mathcal{MT}$, we can define the \emph{merge tree edit distance} 
\[d_E((T,h_T),(T',h_{T'})):= d_K((T,h_T),(T',h_{T'}))
\]
 for any $K\geq \max\{\max h_T, \max h_{T'}\}$.
\end{teo}

For ease of notation we call $d_E$ both the edit distance between weighted trees and the one between merge trees. The arguments of the distances should clarify which of the two metrics we are referring to.
By \Cref{rmk:inverting_repr} we also have the following corollary of \Cref{teo:d_K}. 

\begin{cor}\label{cor:d_E}
The distance $d_E$ is a pseudo metric on $\mathcal{MT}$ and a metric on $\mathcal{MT}/\sim_2$: given two merge trees $(T,h_{T})$ and $(T',h_{T'})$,
$d_E((T,h_T),(T',h_{T'}))=0$ if and only if $(T,h_T)\sim_2(T',h_{T'})$.
\end{cor}

\begin{rmk}\label{rmk:d_K}
\Cref{teo:d_K} and \Cref{cor:d_E} have a series of important implications. First, they say that 
$(\mathcal{MT}_K/\sim_2,d_E)$ is isometric and isomorphic to $(\mathcal{T}_2,d_E)$ and thus, if we have a subset of merge trees contained in $\mathcal{MT}_K/\sim_2$, for some $K$, 
we can map them in $(\mathcal{T}_2,d_E)$ via $\Tr_K$ and carry out our analysis there.
Second, suppose we are given a merge tree $T''$ with $\max h_{T''}>K$.
For any two merge trees $T,T'$ with $K\geq \max h_T,\max h_{T'}$, we can consider $K'\geq \max h_{T''}$ and compute $d_E(T,T'')=d_{K'}(T,T'')$ and $d_E(T',T'')=d_{K'}(T',T'')$. But we do not have to compute again $d_{K'}(T,T')$ for we have $d_E(T,T')=d_{K'}(T,T')=d_{K}(T,T')$.  
\end{rmk}

See \Cref{sec:order_2} for some examples.
We close this section with the following remark. 

\begin{figure}[!h]
\centering
	\begin{subfigure}[c]{0.49\textwidth}
    	\includegraphics[width = \textwidth]{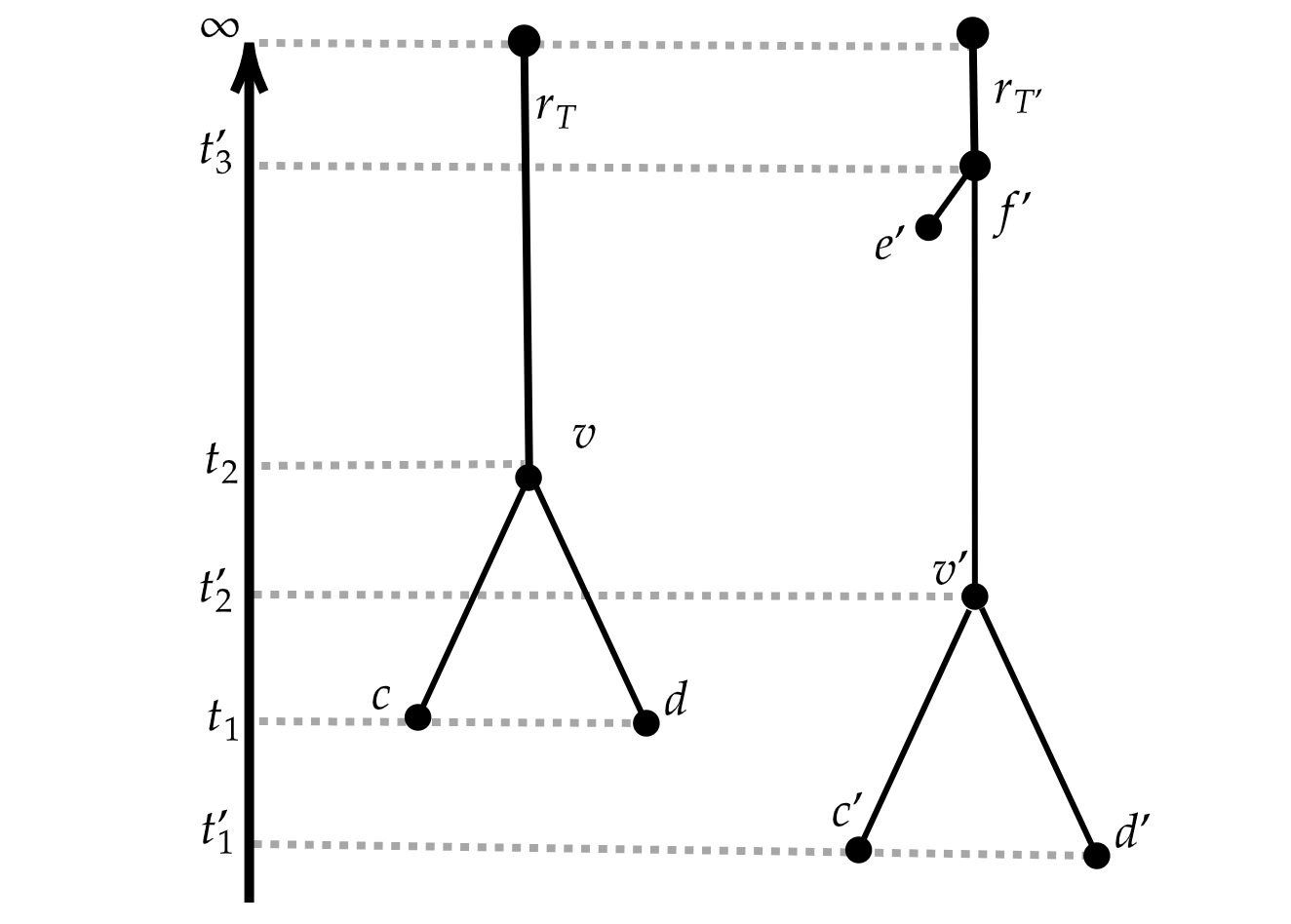}
    	\caption{Two merge trees $T$ (left) and $T'$ (right).}
    	\label{fig:unpleasant_0}
    \end{subfigure}
    
    	\begin{subfigure}[t]{0.49\textwidth}
    	\centering
    	\includegraphics[width = \textwidth]{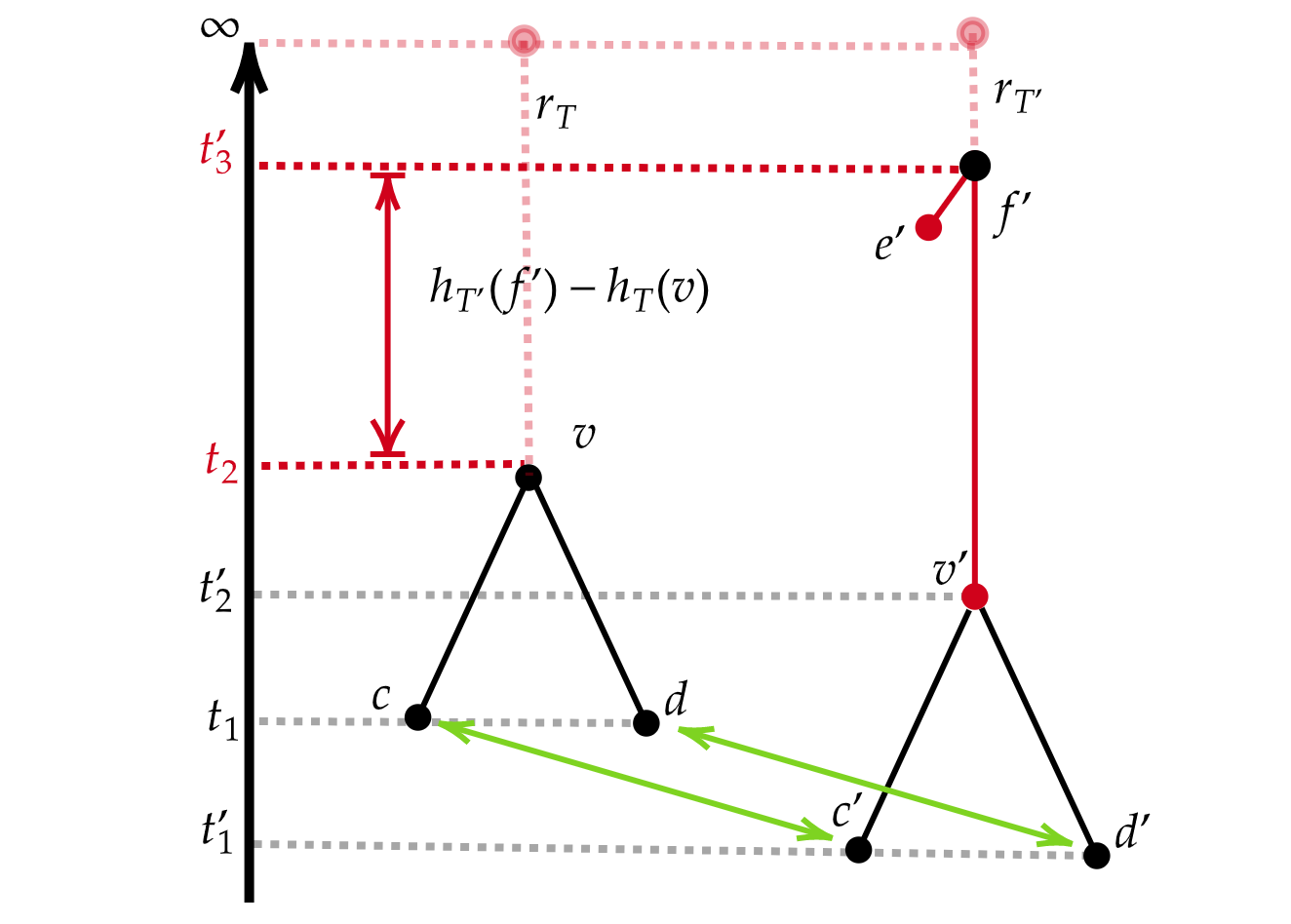}
    	\caption{The two merge trees $T$ (left) and $T'$ (right) in \Cref{fig:unpleasant_0} represented as triplets of the form $(G,w_{G},h_T(v))$ - see \Cref{rmk:naive}. The shaded red edges represent the merge trees $T$ and $T'$, which are removed to obtained the weighted trees. The height values $h_T(v)$ are colored in red. The red edges represent deletions and the green arrows couple the remaining edges, giving a mapping. The cost of such mapping is given by the height differences $\mid h_T(v)-h_{T'}(f')\mid $ plus the two deletions and the weight differences between the coupled edges (which visually appear to be very small). Clearly such distance is inflated by having to account both for $\mid h_T(v)-h_{T'}(f')\mid $ and for the deletion of $(v',f')$.}
    	\label{fig:unpleasant_1}
    \end{subfigure}
   	\begin{subfigure}[t]{0.49\textwidth}
    	\centering
    	\includegraphics[width = \textwidth]{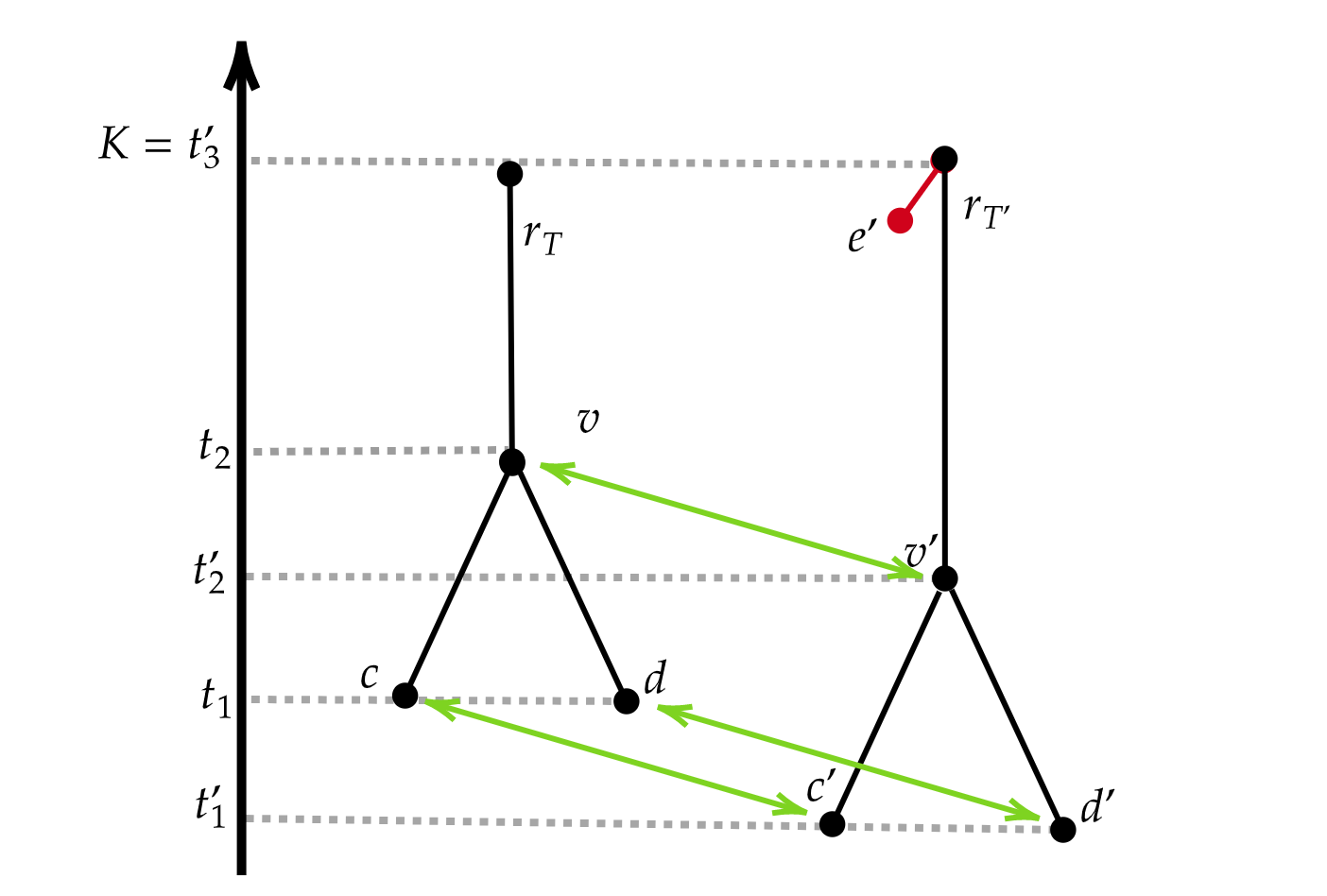}
    	\caption{The weighted trees $\Tr_K(T)$ (left) and $\Tr_K(T')$ (right) with $K=t'_3$, obtained from the merge trees $T$ and $T'$ as in \Cref{fig:unpleasant_0}.  The red edges represent deletions and the green arrows couple the remaining edges, giving a mapping. The cost of such mapping is given by the weight differences between the coupled edges. Note that the weight difference between $(v,r_T)$ and $(v',r_{T'})$ - i.e. $h_T(v)-h_{T'}(v')$ - is surely smaller than the cost of deleting $(v',f')$}.
    	\label{fig:unpleasant_2}
    \end{subfigure}
\caption{A situation in which the metric defined in \Cref{rmk:naive} shows an unstable behavior, while $d_K$ is still well behaved.}
\label{fig:unpleasant}
\end{figure}

\begin{rmk}\label{rmk:naive}
A more naive approach could have been to model each merge tree as a triplet $(G,w_{G},h_T(v))\in (\mathcal{T},\mathbb{R}_{\geq 0})\times \mathbb{R}$, with $v = \arg\max h_T$ and $G=sub_T(v)$ - i.e. to record the height of the last merging point in $T$ and then remove $(v,r_T)$ from $E_T$. Then one could define a pseudo metric on $\mathcal{MT}$ via the rule: 
\[
d((T,h_T),(T',h_{T'}))=\mid h_T(v)-h_{T'}(v')\mid + d_E((G,w_G),(G',w_{G'})).
\]
but this leads to unpleasant and unstable behaviors, as in \Cref{fig:unpleasant}. 
\end{rmk}

By definition, computing $d_E$ for merge trees amounts to computing $d_K$, and so we can exploit the algorithm presented in \cite{pegoraro2023edit} to compute $d_E$ for weighted trees.

\begin{prop}[\cite{pegoraro2023edit}]\label{prop:complexity}
Let $T$ and  $T'$ be two merge trees with full binary tree structures with $\#E_T = N$ and $\# E_{T'}=M$.
Then $d_E(T,T')$ can be computed by solving $O(N\cdot M)$ BLP problems with $O(N\cdot \log(N)\cdot M\cdot \log(M))$ variables and $O(\log_2(M)+ \log_2(N))$ constraints. 
\end{prop}

As anticipated in \Cref{sec:intro}, the classical edit distance can be computed by solving $O(N\cdot M)$ BLP problems with $O(N\cdot  M)$ variables and $O(\log_2(M)+ \log_2(N))$ constraints \citep{TED}. Thus, \Cref{prop:complexity} also suggests that there could be polynomial time approximation algorithms also for the metric $d_E$, as for the classical edit distance \citep{zhang1996constrained}.

\subsection{Stability}
\label{sec:stability}

Merge trees are trees which are used to represent data and, considering populations of such objects, to explore the variability of a data set. Thus, any metric which is employed on merge trees (and on any data representation in general) needs to measure the variability between merge trees in a \virgolette{sensible} way, coherently with the application considered. A formal way to assess such properties often comes in the form of continuity results of the operator mapping data into representations, which, in TDA's literature, are usually referred to as stability results. 

Stability results are usually understood in terms of interleavings between persistence modules \citep{chazal2008persistencemodules}: the distance one measures between (a summary of) the persistence modules should not deviate too much from the interleaving distance between the modules themselves. The rationale behind these ideas follows from the fact that given two functions $f,g:X\rightarrow \mathbb{R}$ with suitable properties and such that $\parallel f-g\parallel_\infty<\varepsilon$, then their sublevel set filtrations induce $\varepsilon$-interleaved persistence modules via homology functors. Moreover, interleaving distances are taken as reference because they are \emph{universal} among the metrics bounded from above by $\parallel f-g\parallel_\infty$ \citep{interl_multi, cardona2021universal}: for any other metric $d$ on persistence modules/merge trees  such that $d(S_f,S_g)\leq \parallel f-g\parallel_\infty$  
then $d(S_f,S_g)\leq d_I(S_f,S_g)$, with $S_f,S_g$ being the persistence modules/merge trees representing $f,g$ and $d_I$ being the interleaving distance between them. Thus, a good behavior in terms of interleaving distance implies a good handling of pointwise noise between functions and so also interpretability of the metric.

Being the edit distance a summation of the costs of local modifications of trees we expect that $d_E$ cannot be bounded from above by the interleaving distance between merge trees as such metric (see \Cref{sec:interl}) in some sense, measures only the biggest modification needed optimally match two merge trees. 
Thus, a suitable stability condition, which we name \emph{finite stability}, would be for $\varepsilon$-interleaved trees to be obtained one from the other just by means of edits with cost $\varepsilon$, and with edit paths whose number of operations is bounded by the sizes of the trees.
That is, the cost of the local modifications we need to match the two merge trees goes to $0$ as their interleaving distance gets smaller and smaller. As is the case for Wasserstein metrics between PDs.

In the following we need to use persistence modules, persistent sets and, in particular, those persistent sets that can be represented via merge trees, which we call \emph{regular abstract merge trees} (RAMTs). For the formal definitions, see \Cref{sec:abstract_merge}. Here we just say that an abstract merge tree is a functor $\T:(\R,\leq) \rightarrow \FSet$ that before some real number is always empty, while after some other real number is always a singleton, and its arrows are all isomorphisms but for a finite number of $t,t'\in \R$, which are called \emph{critical values}. Regularity refers to the behavior of $\T$ at critical values. To such $\T$ we can assign a merge tree $\mathcal{M}(\T)$ without order $2$ vertices. Informally, $\T$ is an algebraic representation of the geometric realization of the tree structure of a merge tree, with the height function defined on all the points of such geometric realization.
Similarly, we make use of the interleaving distance and $\varepsilon$-compatible maps between persistent sets, which are established in TDA literature and thus we just report their definition in \Cref{sec:interl}.

\begin{defi}
Given a persistence module $S:\mathbb{R}\rightarrow \Vect_{\mathbb{K}}$, we define its rank as 
$\rank(S):=\# PD(S)$ i.e. the number of points in its persistence diagram (if it exists). When $S$ is generated on $\mathbb{K}$ by a regular abstract merge tree $\T$ we have  $\rank(S):=\# PD(S) =\# L_{T}$, with $(T,h_T)=\mathcal{M}(\T)$. In this case, we use $\rank(T)$ to refer to $\rank(S)$. We also fix the notation $\dim(T):=\# E_T$. Note that $\dim(T)\leq 2\rank(T)$. Lastly, we define the space of merge trees with uniform upper bound on their dimension:
\[
\MT^N = \{T\in \MT \mid \dim(T)\leq N\}.
\]
\end{defi}

With these pieces of notation we can introduce the notion of \emph{finitely stable} metrics. 

\begin{defi}
A metric for merge trees is finitely stable if there is $C>0$ such that for every RAMTs $\T$ and $\G$, upon setting $T=\mathcal{M}(\T)$ and $T'=\mathcal{M}(\G)$, we have:
\[
d(T,T')\leq C(\rank(T)+\rank(T'))\varepsilon.
\]
\end{defi}

In view of these definitions we can rewrite the following theorem from \cite{pegoraro2024functional}.

\begin{teo}[\cite{pegoraro2024functional}]
\label{teo:stability}
If there are $\alpha,\beta$ $\varepsilon$-compatible maps between two RAMTs $\T$ and $\G$, then there exist a mapping $M$ between $T=\mathcal{M}(\T)$ and $G=\mathcal{M}(\G)$ such that $cost_M((a,b))\leq 2 \varepsilon$ for every $(a,b)\in M$.
Since $\dim(T)\leq 2\rank(T)$, $d_E$ is finitely stable.
\end{teo}

Now we prove a complementary result, which gives the complete picture of the stability properties of $d_E$, allowing the comparison with Wasserstein distances between PDs. Note that in the remaining part of the manuscript, to lighten the notation, we often deliberately confuse the a RAMT $\T$ and the associated merge tree $\mathcal{M}(\T)$.

\begin{teo}\label{prop:d_I<d_E}
We always have $d_I(T,G)\leq d_E(T,G)$. 
\end{teo}

The bound given by \Cref{prop:d_I<d_E} cannot be improved.

\bigskip

\subparagraph*{Example}
Consider $f:I\rightarrow \mathbb{R}$ and $g(x)=f(x)+h$ for some fixed $h\in\mathbb{R}$, and let $\X$ be the sublevel set filtration of $f$ and $\Y$ of $g$. Taking $(T,h_T)$ and $(G,h_G)$ as the merge trees representing $\T$ and $\G$, we have $ d_I(T,G)=d_E(T,G)=h$. 

\bigskip

Putting together \Cref{prop:d_I<d_E} and \Cref{teo:stability} we obtain the following inequalities.

\begin{cor}
\begin{equation}
d_I(T,G)\leq d_E(T,G)\leq 2(\dim(T)+\dim(G))d_I(T,G).
\end{equation}
\end{cor}

\Cref{eq:interl} is the inequality that better summarizes the stability properties of $d_E$.
Moreover, inequalities similar to the ones expressed by \Cref{eq:interl} relate also the bottleneck and the $1$-Wasserstein distances between PDs. See \Cref{sec:wass} for definitions. 
In fact, if $d_B$ is the bottleneck distance between PDs, and $W_1$ is the $1$-Wasserstein distance, for every pair of persistence diagrams $D_1,D_2$ we have:
\begin{equation}
d_B(D_1,D_2)\leq W_1(D_1,D_2)\leq (\rank(D_1)+\rank(D_2))d_B(D_1,D_2).
\end{equation}

The similarity between \Cref{eq:interl} and \Cref{eq:wass} qualifies $d_E$ as an analogous for merge trees of the $1$-Wasserstein distance for PDs, as they have analogous stability properties w.r.t. the universal distance of, respectively, merge trees and PDs. 
Those properties are in line with our expectations: when editing $\varepsilon$-interleaved merge trees we need to produce small local modification for each vertex of the merge tree, and then add up all the contributions. In the same way, with  the $1$-Wasserstein distance, one measures the difference between $\varepsilon$-interleaved diagrams aggregating the small discrepancies between all persistence pairs.

Universal distances like the interleaving and the bottleneck ones satisfy such strong stability properties because, in some sense, they measure the maximum cost of the modifications that we have to make on the considered objects. In other words, they are very stable because they are heavily insensitive: one could add an infinite amount of deformations to the considered objects, smaller than the biggest one (like, infinite points close enough to the diagonal of PDs) without these changes affecting the distance. 
That is, they are very poor in discriminating between different objects, and not just when local deformations are small, as shown in \Cref{sec:vs_interl}. More formally, the topology induced by $d_E$ is strictly richer than the one induced by $d_I$ as, for any $T'\in \mathcal{MT}$ and $C>0$, the following hold:
\[
\{T\in \mathcal{MT}\mid d_E(T,T')<C\}\subset\{T\in \mathcal{MT}\mid d_I(T,T')<C\}. 
\]
On the other hand it is not difficult to build a tree $T$ and a sequence $\{T_n\}_{n\in \mathbb{N}}$ such that $d_I(T_n,T)\rightarrow 0$ and $d_E(T_n,T)\rightarrow +\infty$. See also the upcoming example.

\bigskip

\subparagraph*{Example}
Let $h_n:[0,1]\rightarrow \R$ be $h_n(x)=1/n\cdot \sin(2\pi n^2x)$ and $f:[0,1]\rightarrow \R$ being $f(x)=x$. Let $T$ be the merge tree of $f$ and $T_n$ the merge tree of $f_n:=f+h_n$.
We have $d_I(T,T_n)\leq 1/n$ while $d_E(T,T_n)\sim O(n)$. Similarly, let $D$ be the PD of $f$ and $D_n$ the PD of $f_n$.
We have $d_B(D,D_n)\leq 1/n$ while $W_1(D,D_n)\sim O(n)$. 

\bigskip

This reinforces the fact that, in general, $d_E$ is potentially better than $d_I$ at discriminating trees, as is has more open sets to separate objects.
From the data analysis perspective these facts have interesting consequences, which we try to sum up merging some ideas which are commonly found in TDA theoretical works, with others which are diffused among practitioners.

Universal metrics are of great use whenever we need to work with topological summaries extracted from noisy data, that we cannot carefully denoise. 
Suppose for instance that we observe an iid sample $\{(x_i,y_i)\}_{i=1}^n$ from the model $\mathcal{Y} \mid \mathcal{X}=x \sim f(x) + \varepsilon$, with $\varepsilon$ being a noise random variable with zero mean and finite variance, $\mathcal{X}, \mathcal{Y}$ being two real valued random variables, and $f:[a,b]\rightarrow \R$ a smooth function. Which is a standard statistical model in functional data analysis and regression problems \citep{book_fda}. Suppose we can build an estimate $\hat{f}_n$ of $f(\cdot)=\mathbb{E}(\mathcal{Y} \mid \mathcal{X}=\cdot)$ from $\{(x_i,y_i)\}_{i=1}^n$. Universal stability implies that $\hat{f}_n$ can be a very naive estimator of $f$, as $d_I(T_f,T_{\hat{f}_n}) \leq \parallel f-\hat{f}_n\parallel_\infty.$

Thus, as long as the estimator provides a low point wise error, we know that $T_{\hat{f}_n}$, which we obtained from our samples, is a good estimate of $T_f$ in terms of $d_I$. 
Having a very rough estimate $\hat{f}_n$ (e.g. piecewise affine interpolation of $\{(x_i,y_i)\}_{i=1}^n$) would instead be very bad for $d_E$ and $W_1$ as $\hat{f}_n$ is likely to have many ancillary/noisy oscillations, which can blow up $\dim(T_{\hat{f}_n})$ and $\rank(PD(\hat{f}_n))$,  so that  $T_{\hat{f}_n}$ and $PD(\hat{f}_n)$ are potentially very poor approximations of $T_f$ and $PD(f)$ in terms of $d_E$ and $W_1$, respectively. So, whatever analysis is ran with such metrics will produce no reliable results.

However, as shown in \cite{pegoraro2024functional}, one can build $\hat{f}_n$ with good enough properties (e.g. not just uniform convergence, but also uniform convergence of the derivatives), so that $T_{\hat{f}_n}$ ends up being a good estimator of $T_f$ also in terms of $d_E$. Thus, upon designing careful denoising procedures, for instance relying on more refined statistical estimators for the kind of data considered, one can solve the problem of dealing with noise, and is not constrained to use universally stable metrics. Of course, the same applies if the analysts and the field experts can assume that the data generating process is affected by noise only in negligible terms. 

In any of the aforementioned situations, the analysts can build statistically grounded pipelines involving metrics like $d_E$ and $W_p$ which are much more sensitive and whose discriminative power is arguably more useful in many data analysis scenarios compared to their universal counterparts, to the point that, to the best of our knowledge, Wasserstein metrics are more frequently used in applications than the bottleneck one. On the other hand, if data is heavily corrupted by noise and cannot be effectively denoised, this will affect the results obtained with $W_p$ and $d_E$ much more than if the analysis had been carried out with $d_B$ and $d_I$. 

We complement these comments with two additional results. With the first easy one, which we state without proof, we motivate our terminology of \virgolette{finite} stability: it indicates that $d_I$-Lipschitz operators (like the one mapping tame functions with the sup norm into the merge trees of their sublevel sets) are also $d_E$-Lipschitz on all finite dimensional subspaces of merge trees.
With the second one, instead, we highlight that any $\varphi:(X,d_X)\rightarrow (\MT,d_I)$ which is continuous, is also arbitrarily close (in terms of $d_I$) to being $d$-continuous, for every finitely stable metric $d$. More precisely, $\varphi$ 
can locally be approximated via its smoothing: $\S_\varepsilon \circ \varphi$. Where $\S_\varepsilon$ is the smoothing operator used to define the interleaving distance in \Cref{sec:interl}.

\begin{prop}\label{prop:finitely}

    Let $d$ be a finitely stable metric between merge trees and let $\varphi:(X,d_X)\rightarrow (\MT,d_I)$ be a Lipschitz operator from some metric space $(X,d_X)$ to the space of regular merge trees considered with the interleaving distance. Then $\varphi$ is Lipschitz w.r.t. $d$ on all finite dimensional subspaces of $\MT$. That is, for every $N\in \N$ the following operator is Lipschitz:
    \[
    \varphi_{\mid \varphi^{-1}(\MT^N)}: (\varphi^{-1}(\MT^N),d_X)\rightarrow (\MT^N,d).
    \]
\end{prop}

\begin{prop}\label{prop:locally}

    Consider $T$ and $\{T_n\}_{n\in \N}$ merge trees such that $d_I(T,T_n)\rightarrow 0,$
    and let $d$ be a finitely stable distance between merge trees. For any $\epsilon>0$,  there exists $\varepsilon>0$ (depending also on $T$) such that, for $G_n:=\S_{\varepsilon}(T_n)$ we have:
    \begin{enumerate}
        \item $d_I(T_n,G_n)\leq \epsilon$ for all $n$;
        \item $\lim_{n\rightarrow \infty} d(T,G_n)\leq \epsilon.$
    \end{enumerate}
\end{prop}

\section{Conclusions}
\label{sec:conclusions_geom}
In this manuscript we introduce a novel edit distance for merge trees, and study its stability properties. To the best of our knowledge, our metric is the only edit distance for merge trees with proven stability properties.
We devote a good portion of the work to support the claim that, a priori, the stability properties we prove are more suited for general data analysis scenarios w.r.t. the universal ones of the interleaving distance. 
The appendix further expands on that, with a detailed comparison with the other edit distances for merge trees and a simulation involving the interleaving distance.  
Notably, other works employ this metric in applications \citep{pegoraro2024functional, cavinato2022imaging}. 

This works opens a series of further developments in complementary directions: 1) aiming at closing the computational gap with the other edit distances for merge trees by obtaining an algorithm for the constrained version of our edit distance. In the manuscript we claim that this algorithm should have polynomial time complexity, and we support this claim with the relationships between our metric and the classic edit distance between unlabeled trees; 2) developing data analysis tools to work in the space of merge trees. In particular, we are currently working on defining means of finite sets of trees and in defining embeddings of finite sets of trees into euclidean spaces, trying to preserve some metric properties of the original set of trees.

\newpage

\appendix{

\section{Abstract Merge Trees}\label{sec:abstract_merge}

Following \cite{patel2018generalized, curry2021decorated} we introduce persistent sets.
 \Cref{fig:preliminary} illustrates some of the objects we introduce in this section.

\begin{figure}
    \begin{subfigure}[c]{0.49\textwidth}
    	\centering
    	\includegraphics[width = \textwidth]{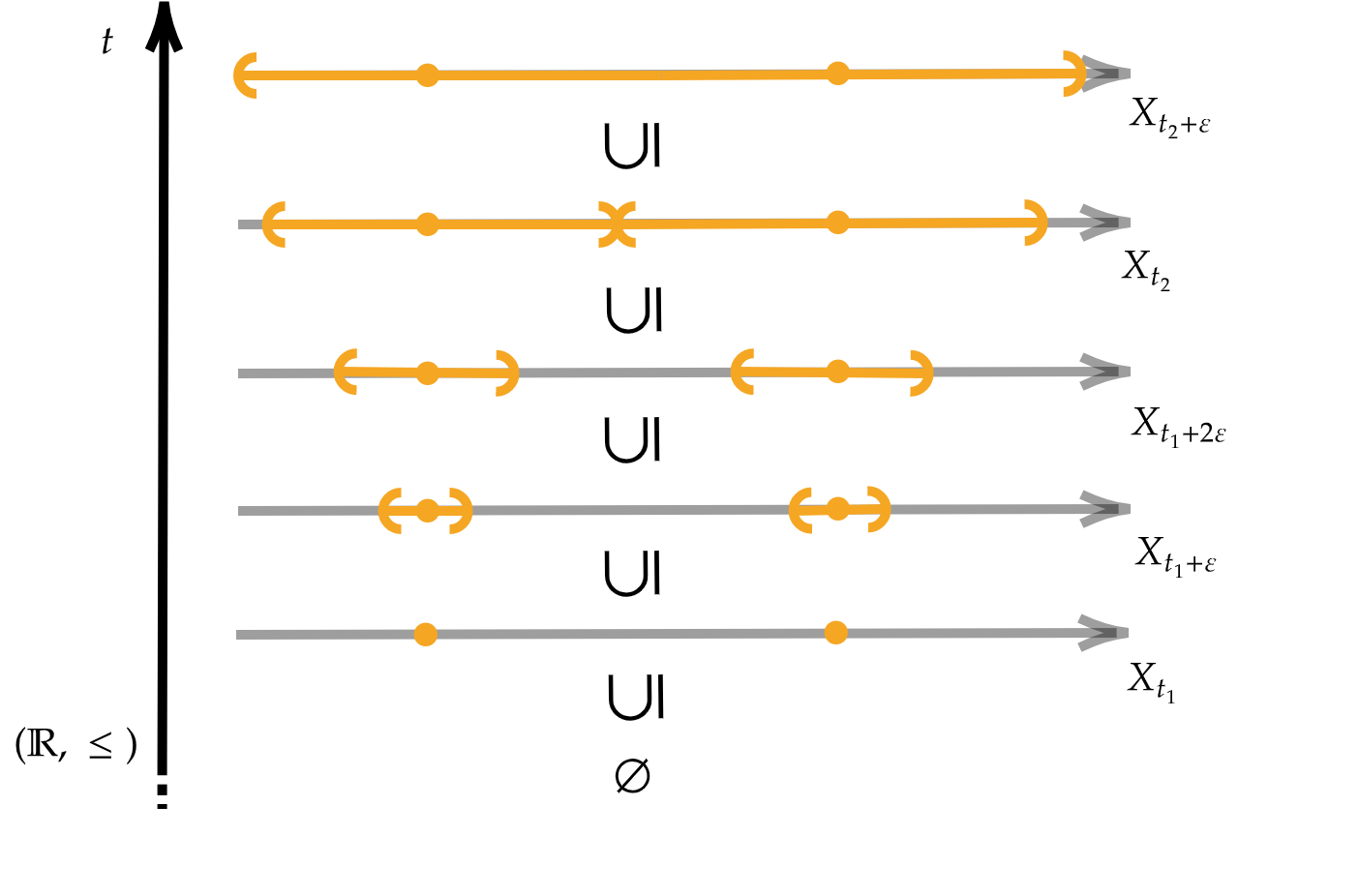}
    	\caption{A filtration $\X$.}
    	\label{fig:filtr}
	\end{subfigure}
    \begin{subfigure}[c]{0.49\textwidth}
    	\centering
	    \includegraphics[width = \textwidth]{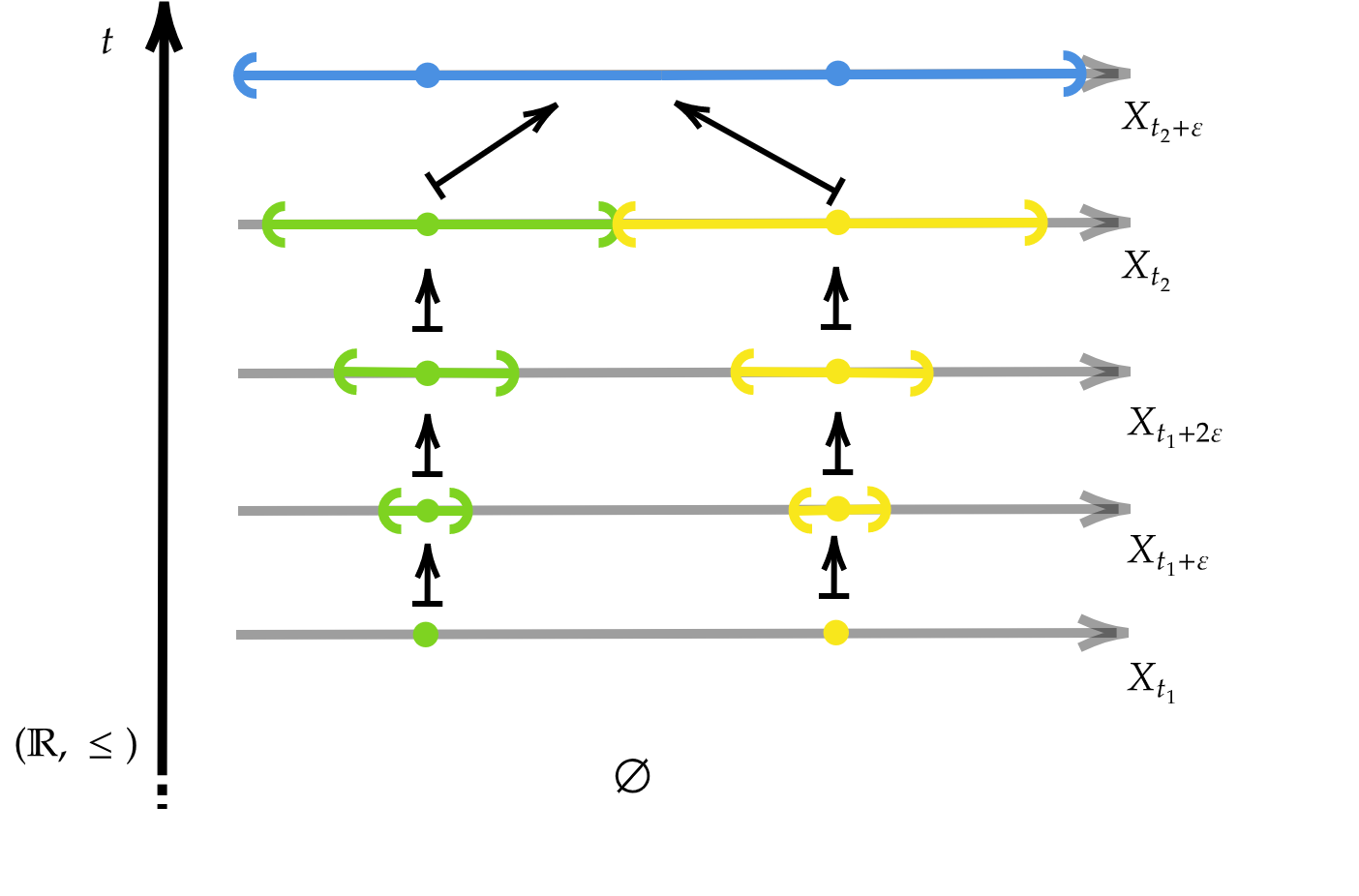}
    	\caption{An abstract merge tree $\T$.}
    	\label{fig:abs_MT}
	\end{subfigure}


    \begin{subfigure}[c]{0.49\textwidth}
    	\centering
    	\includegraphics[width = \textwidth]{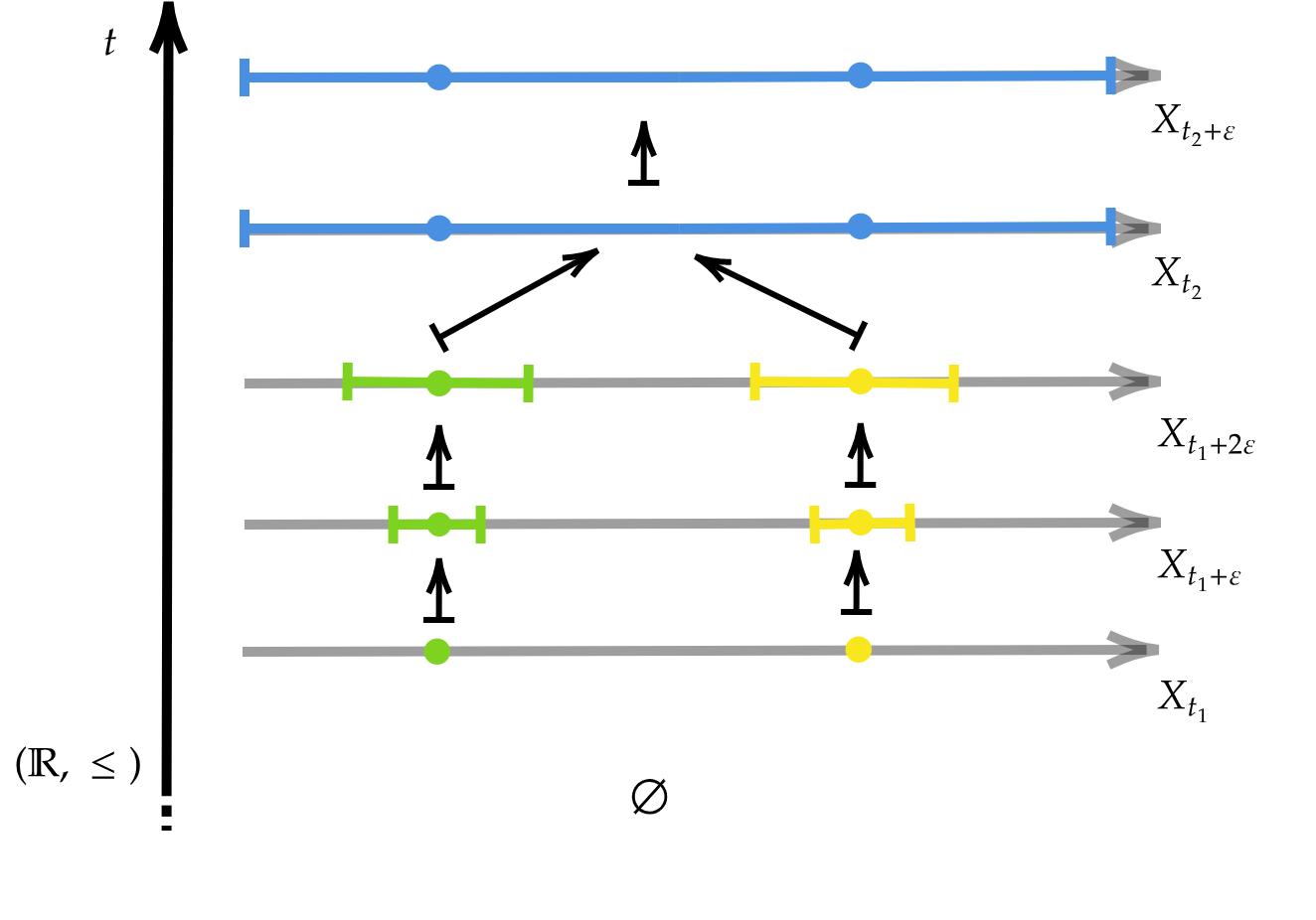}
    	\caption{A regular abstract merge tree $\T$.}
    	\label{fig:regular_MT}
	\end{subfigure}
    \begin{subfigure}[c]{0.49\textwidth}
    	\centering
	    \includegraphics[width = \textwidth]{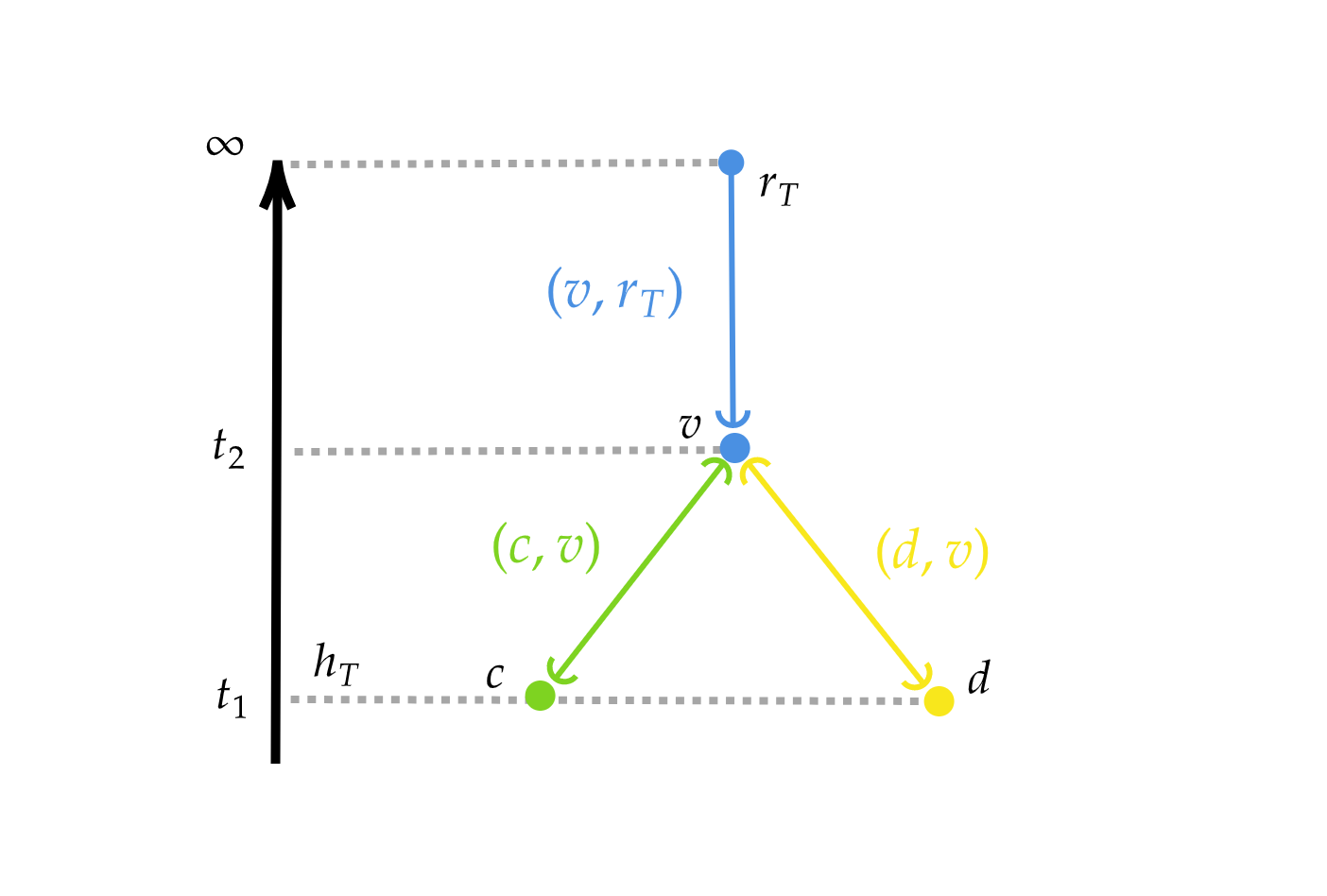}
    	\caption{The merge tree $\mathcal{M}(\T)$ with $\T$ as in \Cref{fig:regular_MT}.}
    	\label{fig:MT}
	\end{subfigure}

\caption{On the first line we see an example of a filtration along with its abstract merge tree. in In the bottom line there are a regular abstract merge tree and the associated merge tree.
 The colors are used throughout the plots to highlight the relationships between the different objects.}
\label{fig:preliminary}
\end{figure}

\begin{defi}[\cite{curry2021decorated}]
A filtration of topological spaces is a (covariant) functor 
$X_{\Bigcdot}:\mathbb{R}\rightarrow \Top$ from the poset $(\mathbb{R},\leq )$ to $\Top$, the category of topological spaces with continuous functions, such that:  
$X_t \rightarrow X_{t'}$, for $t<t'$,
are injective maps. 
\end{defi}

\subparagraph*{Example}
Given a real valued function $f:X\rightarrow \mathbb{R}$ the \emph{sublevel set} filtration is given by $X_{t}=f^{-1}((-\infty,t])$ and $X_{t<t'}=i:f^{-1}((-\infty,t])\hookrightarrow f^{-1}((-\infty,t'])$.

\subparagraph*{Example}
Given a finite set $C\subset \mathbb{R}^n$ its the \emph{C\'ech} filtration is given by $X_{t}=\bigcup_{c\in C}B_{t}(c) $. With $B_{t}(c)=\{x\in\mathbb{R}^n \mid \parallel c-x \parallel < t\}$. As before: $X_{t<t'}=i:\bigcup_{c\in C}B_{t}(c)\hookrightarrow \bigcup_{c\in C}B_{t'}(c)$.

\bigskip

Given a filtration $\X$ we can compose it with the functor $\pi_0$ sending each topological space into the set of its path connected components.
We recall that, according to standard topological notation, $\pi_0(X)$ is the set of the path connected components of $X$ and, given a continuous functions $q:X\rightarrow Y$, $\pi_0(q):\pi_0(X)\rightarrow \pi_0(Y)$ is defined as:

\[
U\mapsto V  \text{ such that }q(U)\subset V.
\]

\begin{defi}[\cite{carlsson2013classify, curry2018fiber}]
A persistent set is a functor $S:\mathbb{R}\rightarrow \Sets$. In particular, given a filtration of topological spaces $\X$, the persistent set of components of $\X$ is $\pi_0 \circ \X$. A (one dimensional persistent module) is a functor $S:\mathbb{R}\rightarrow \Vect_{\mathbb{K}}$ with values in the category of vector spaces 
$\Vect_{\mathbb{K}}$.
\end{defi}

By endowing a persistent set with the discrete topology, every persistence set can be seen as the persistence set of components of a filtration. Thus a general persistent set $S$ can be written as $\pi_0(\X)$ for some filtration $\X$. 

Based on the notion of constructible persistent sets found in \cite{patel2018generalized, curry2021decorated} one then builds the following objects.

\begin{defi}[\cite{pegoraro2024finitelyfunc}]\label{def:abstract_mt}
An abstract merge tree is a persistent set $S:\mathbb{R}\rightarrow \Sets$ such that there is a finite collection of real numbers $\{t_1<t_2<\ldots<t_n\}$ which satisfy:
\begin{itemize}
\item $S(t)=\emptyset$ for all $t<t_1$;
\item $S(t)=\{\star\}$ for all $t>t_n$;
\item if $t,t'\in (t_i,t_{i+1})$, with $t<t'$, then $S(t<t')$ is bijective.
\end{itemize}
The values $\{t_1<t_2<\ldots<t_n\}$ are called critical values of the tree and there is always a minimal set of critical values \citep{pegoraro2024finitelyfunc}. We always assume to be working with such minimal set.

If $S(t)$ is always a finite set, $S$ is a finite abstract merge tree. 
\end{defi}

Consider an abstract merge tree $\T$ and let $t_1<t_2<\ldots<t_n$ be its (minimal set of) critical values and let $i_{t}^{t'}:= X_{t\leq t'}:X_t\rightarrow X_{t'}$. Take $\varepsilon>0$ small enough. We have that at least one between  $\pi_0(i_{t_j-\varepsilon}^{t_j})$ and $\pi_0(i_{t_j}^{t_j+\varepsilon})$ is not bijective. So we have the following definition.

\begin{defi}[\cite{pegoraro2024finitelyfunc}]
An abstract merge tree $\T$ is said to be regular if $\pi_0(i_{t_j}^{t_j+\varepsilon})$ is bijective for every critical value $t_j$ and for every $\varepsilon>0$ small enough. 
\end{defi}

\begin{assump}
We consider only regular abstract merge trees. In \cite{pegoraro2024finitelyfunc} it is shown that this choice is non-restrictive.
\end{assump}

We report a result summarizing the relationship between abstract merge trees and merge trees. The main consequence of such result is that merge trees considered up to order $2$ vertices are an appropriate discrete tool to represent the information contained in RAMTs.   
\Cref{fig:abs_MT} and \Cref{fig:MT} can help the reader going through the following proposition.

\begin{prop}[\cite{pegoraro2024finitelyfunc}]
\label{prop:equivalence}
The following hold:
\begin{enumerate}
\item we can associate a merge tree without order $2$ vertices $\mathcal{M}(\T)$ to any regular abstract merge tree $\T$;
\item we can associate a regular abstract merge tree $\mathcal{F}((T,h_T))$ to any merge tree $(T,h_T)$. Moreover, we have $\mathcal{M}(\mathcal{F}((T,h_T)))\cong_2 (T,h_T)$;
\item given two merge trees $(T,h_T)$ and $(T',h_{T'})$, we have $\mathcal{F}((T,h_T))\cong \mathcal{F}((T,h_T))$ if and only if $(T,h_T)\cong_2 (T',h_{T'})$.
\end{enumerate}
\end{prop}

\begin{rmk}
For the sake of brevity, we don't report the explicit construction of $\mathcal{M}$ and $\mathcal{F}$, which can be found in the proof of Proposition 4 in \cite{pegoraro2024finitelyfunc}. Intuitively, $\mathcal{M}$ gives a discretization of a regular abstract merge tree $\T$, by sending it in the \virgolette{smallest} graph (in terms of vertices) needed to describe the maps $\pi_0(X_{t\leq t'})$, together with the induced function defined on the vertices. In particular, this can be done without using order $2$ vertices, which are superfluous. Viceversa, $\F$ associates to a merge tree a functor which is roughly an algebraic equivalent of its geometric realization, where one can actually consider all the individual points on the edges. Clearly, order two vertices are lost in this procedure.
\end{rmk}

\section{Interleaving Distance Between Merge Trees}
\label{sec:interl}

We adapt the definition of the interleaving distance between merge trees, which is originally stated by 
\cite{merge_interl} with a different notation.

\begin{defi}[adapted from \cite{merge_interl} and \cite{de2016categorified}]
Given $\X$ filtration and $\varepsilon>0$,
we define $\X^\varepsilon$ as $X_t^\varepsilon := X_{t+\varepsilon}$ and $X_{t\leq t'}^\varepsilon :=X_{t+\varepsilon\leq t'+\varepsilon}$.
    Coherently with the literature on Reeb graphs, we define the called the $\varepsilon$-smoothing operator $\S_\varepsilon$ as: $\S_\varepsilon(\T):=\pi_0(\X^\varepsilon)$. Note that $\S_\varepsilon$ also acts on natural transformations $\alpha:\T\rightarrow\G$:  by setting $\S_\varepsilon(\alpha_t):=\alpha_{t+\varepsilon}$ one obtains $\S_\varepsilon(\alpha):\S_\varepsilon(\T)\rightarrow\S_\varepsilon(\G)$. Lastly, we have the natural transformation $i_{\T}^\varepsilon: \T \rightarrow \S_\varepsilon(\T)$
    given by $(i_{\T}^\varepsilon)_t:=\pi_0(X_{t\leq t+\varepsilon}):\pi_0(X_t)\rightarrow \pi_0(X_{ t+\varepsilon})$.
\end{defi}

\begin{defi}[adapted from \cite{merge_interl}]
Take two abstract merge trees $\T$ and $\G$.
Two natural transformations $\alpha:\pi_0(\X)\rightarrow \S_\varepsilon(\pi_0(\Y))$, $\beta:\pi_0(\Y)\rightarrow \S_\varepsilon(\pi_0(\X))$ are $\varepsilon$-compatible if:

\begin{itemize}
\item $\S_\varepsilon(\beta)\circ \alpha = i_{\T}^{2\varepsilon}$;
\item $\S_\varepsilon(\alpha)\circ \beta = i_{\G}^{2\varepsilon}$.
\end{itemize}

Unfolding the definitions, this means:
\begin{itemize}
\item $\beta_{t+\varepsilon}\circ \alpha_t = \pi_0(X_{t\leq t+2\varepsilon})$
\item $\alpha_{t+\varepsilon}\circ \beta_t = \pi_0(Y_{t\leq t+2\varepsilon})$.
\end{itemize}

Then, the interleaving distance between $\T$ and $\G$ is:
\[
d_I(\T,\G) = \inf \{\varepsilon>0\mid \exists \alpha,\beta \text{ }\varepsilon\text{-compatible}\}. 
\]

We also say that $\T$ and $\G$ are $d_I(\T,\G)$-interleaved. 
\end{defi}

To lighten the notation we may also write $d_I(T,G)$, implying $d_I(\T,\G)$, with $T=\mathcal{M}(\T)$ and $G=\mathcal{M}(\G)$.

\begin{rmk}
    We have introduced the interleaving distance with this notation, because we use smoothing operators in \Cref{prop:locally}.
    The name smoothing operator comes from the interleaving distance between Reeb graphs. Introducing Reeb graphs and such distance to motivate our definition is outside the scope of this work. However, for the reader familiar with such notions, we point out the following facts. 
   Merge trees can be obtained as Reeb graphs of the epigraph of a function $f:X\rightarrow \R$ with the projection on $\R$: consider $f:X\rightarrow\R$ and define $\Gamma_f =\{(x,t)\in X\times \R \mid f(x)\leq t\}$. Then one can take the projection to $\R$ and obtain: $F:\Gamma_f \rightarrow \R$.
    Then $F^{-1}(t)=\{(x,t)\mid f(x)\leq t\}=f^{-1}((-\infty,t])$. Thus, the merge tree of $f$ is equivalently represented by the Reeb cosheaf \citep{de2016categorified} $U\mapsto \pi_0(F^{-1}(U))$. The smoothed version of this Reeb cosheaf,
    as defined in \cite{de2016categorified}, is induced by $F_\varepsilon^{-1}(t):=F^{-1}((t-\varepsilon,t+\varepsilon))=f^{-1}((-\infty,t+\varepsilon])$.  Which is implies that the smoothed Reeb cosheaf corresponds to the smoothed merge tree in the sense that we define. 
    This also implies that the interleaving distance between merge trees can be induced by the one between Reeb graphs. 
\end{rmk}

\section{Wasserstein Distance Between Persistence Diagrams}\label{sec:wass}

Given two persistence diagrams $D_1$ and $D_2$, the expression of the $p$-Wasserstein distance between them is the following:
\[
W_p(D_1,D_2) = \left( \inf_\gamma  \sum_{x \in D_1} \parallel x-\gamma(x)\parallel^p_\infty \right)^{1/p},
\]
where $\gamma$ ranges over the functions partially matching points between diagrams $D_1$ and $D_2$, and matching the remaining points of both diagrams with the line $y=x$ on the plane (for details see \cite{cohen_PD}). In other words we measure the distances between the points of the two diagrams, pairing each point of a diagram either with a point on the other diagram, or with a point on $y=x$. Each point (also called persistence pair) can be matched once and only once. The minimal cost of such matching provides the distance. The case $p=\infty$ is usually referred to as the bottleneck distance $d_B$.

\section{Constrained Edit Distance}
\label{rmk:constr}
    The constrained edit distance presented in \cite{zhang1996constrained} is obtained adding to the mapping properties of the classical edit distance $D_E$ an additional requirement, which we now report. Consider a weighted tree $(T,w_T)$ and consider $A\subset V_T$. The least common ancestor (LCA) of the vertices in $A$ is defined as:
    \[
    \text{LCA}(A)=\min \bigcap_{v\in A} \{w\in V_T \mid w\geq v\}.
    \]
    Aa mapping $M$ between a weighted tree $(T,w_T)$ and a weighted tree $(T',w_{T'})$ for the constraint edit distance \cite{zhang1996constrained} is a mapping for the classical edit distance with an additional requirement that can be formulated as:
    \begin{itemize}
        \item[(C)]  for every $a,b,c\in E_T$ and $a',b',c'\in E_{T'}$ such that $(a,a'),(b,b'),(c,c')\in M$, we have:
        \[
        \text{LCA}(\{a,b\})\geq c \text{ if and only if } \text{LCA}(\{a',b'\})\geq c'.
        \]
    \end{itemize}
    There are not obstacles to adding requirement (C) to the properties (M1)-(M4), defining a constrained version of $d_E$. Since incorporating the ghosting edits into $D_E$ does not significantly increase its computational cost, as shown by \Cref{prop:complexity}, we strongly believe that also the computational complexity of the constrained edit distance is only marginally affected by such edit operations. And so, in particular, a constrained version of $d_E$ should admit a polynomial time algorithm.

\section{Comparison with Other Distances For Merge Trees}
\label{sec:edit_comparison}

In this section we illustrate some behaviors of the metric $d_E$, we comment on its potential polynomial time approximation, and extensively compare $d_E$ with other definitions of distances between merge trees which appeared in literature. In this way we can better portrait how the variability between merge trees is captured by the proposed edit distance. We consider separately the different metrics we compare $d_E$ against, plus, at the end, we consider some topics which are common to all/most of the different comparisons.  We avoid a comparison with \cite{sridharamurthy2021comparative} as the goal of such metric is comparing subtrees of merge trees, and so it behaves very differently from all other distances.

\subsection{Editing a Merge Tree}
\label{sec:order_2}

We devote this subsection to explore with some easy examples the definitions and results given in  \Cref{sec:truncation} and \Cref{sec:MTED}.

First note that, by construction, $\Tr_K((T,h_T))$ - for $K$ big enough - is a representation of the merge tree $(T,h_T)$ coherent with the metric $d_E$ and thus can be used also to visually compare two merge trees.
We can then consider a merge tree $(T,h_T)$ and edit $(T,w_T)$ according to the rules in \Cref{sec:weights_edits}: via shrinking, deletions and ghosting of vertices and the inverse operations. We look at the results of the edits in light of the merge tree $(\Tr_K)^{-1}((T,w_T))$.

Let $(T,h_T)=\mathcal{M}(\T)$ and consider an edge $e=(v,v')\in E_T$, with $t=h_T(v)$ and $t'=h_T(v')$. Since the height of the root is fixed and equal to $K$, shrinking $e$ reducing its weight by some value $\varepsilon>0$  (with $\varepsilon<w_T(e)$) amounts to \virgolette{moving upwards} $sub_T(v')$ by $\varepsilon$, that is changing $h_T(v'')\mapsto h_T(v'')+\varepsilon$ - as in \Cref{fig:ghostings}$\rightarrow$\Cref{fig:shrinkings} or in \Cref{fig:compare_1} and \Cref{fig:compare_recap} (left). Having $\varepsilon > w_T(e)$ means deleting $e$. Similarly, increasing $w_T(e)$ by $\varepsilon$, amounts to lowering 
$sub_T(v')$ by $\varepsilon$ - as it partially happens in \Cref{fig:insertions} inserting the red internal vertex. 
Consider now the splitting of the edge $e$ into $e_1=(v,v'')$ and $e_2=(v'',v')$ with $T'$ being the novel tree structure and $w_{T'}(e_1)=\varepsilon_1$ and $w_{T'}(e_2)=\varepsilon_2$ - as for any of the yellow vertices in \Cref{fig:splittings}. We must have $\varepsilon_i>0$ and $w_T(e)= \varepsilon_1 + \varepsilon_2$. This clearly induces a well defined height function $h_{T'}(v'')$. The merge tree $(T',w_{T'})$ differs from $(T,w_T)$ by the order two vertex $v''$, while the height function on $V_{T'}-\{v''\}$ is still the same. And, accordingly, the associated RAMTs are the same $\T=\mathcal{F}((T,h_T))\cong \mathcal{F}((T',h_{T'}))$ (with $\mathcal{F}$ bein as in \Cref{prop:equivalence}). Thus we have changed the graph structure of $T$ without changing the topological information it represents.

\begin{rmk}\label{rmk:order_p}
In light of this paragraph, it would be natural to try to define a family of metrics indexed by integers $p\geq 1$ by saying that the costs of an edit path the $p$-th root of sum of the costs of the edit operations to the $p$-th power. But now we can easily see that for any $p>1$ this has no hope of being a meaningful pseudo metric for weighted trees. In fact, consider the case of a weighted tree made by two vertices and one edge with weight $1$. The cost of shrinking the $p$-metric would be $\parallel 1\parallel_p =1$. At the same time one can split it in half with $0$ cost and the cost of shrinking this other tree would be $\parallel(1/2,1/2)\parallel_p<1$. Splitting the segment again and again will make its shrinking cost go to $0$. In other words all weighted trees, if considered up to order $2$ vertices, would be at distance zero from the tree with no branches.
\end{rmk}
 
\subsection{Interleaving Distance \citep{merge_interl} and Persistence Diagrams}

\label{sec:vs_interl}

We extensively discussed similarities and differences between $d_I$ and $d_E$ in \Cref{sec:stability}. We complement such analysis with a simulation to further showcase the differences between $d_E$ and $d_I$ (see also \Cref{fig:sim_1_fn}).

For $i=0,\ldots,9$, let $g_i:[0,11]\rightarrow \mathbb{R}$ be such that $g_i \equiv 0$ on 
$[0,11]- [i+1/3,i+2/3]$ while, on 
$[i+1/3,i+2/3]$, \(g_i\) is the linear interpolation of $(i+1/3,0)$, $(i+1/2,1)$ and $(i+2/3,0)$. Then, for $i=0,\ldots,9$, define $G_i$  as $G_i\equiv 0$ on
$[0,11]- [i+2/3,i+1]$, while, on $[i+2/3,i],$  $G_i$ is the 
linear interpolation of $(i+2/3,0)$, $(i+3/4,5)$ and $(i+1,0)$.

Then $f_i$, $i=0,\ldots,9$, is obtained as follows:
\[
f_i = G_i+\sum_{j=0}^{9} g_j. 
\]
See \Cref{fig:sim_1_fn} to better visualize this data set: we have a constant set of lower peaks at height $1$ and an higher peak with height $5$ which is shifting left to right as $i$ increases.
In this way we are just changing the left-right distribution
of the smaller peaks wrt the highest one. 

We obtain the associated merge trees and then compute the pairwise distances between the merge trees with $d_E$. The results are represented in 
\Cref{fig:sim_1_dist}:
the shortest edit path between the $i$-th merge tree and $i+1$-th is given by the deletion of one leaf in each tree to make the disposition of leaves coincide between the two trees. The more the peaks' disposition is different between the two trees, the more one needs to delete leaves in both trees to find the path between them. Note that
the first function (the one in which the highest peak is the second peak) and the last function (the one in which the highest peak is the second-last peak) can be obtained one from the other via a $y$-axis symmetry and translation. Similarly, the second function is equal, up to homeomorphisms of the domain, to the third-last one, etc. (see also \cite{pegoraro2024functional}). Thus the merge trees are the same. To sum up the situation depicted in the first row of \Cref{fig:sim_1_dist}, first we get (left-to-right) farther away from the first merge tree, and then we return closer to it. This intuition is confirmed by looking at the multidimentional scaling (MDS) embedding in $\mathbb{R}^2$ of the pairwise distance matrix (see \Cref{fig:sim_1_MDS} - note that the shades of gray reflect, from white to black, the ordering of the merge trees). The discrepancies between the couple of points which should be identified are caused by numerical errors.

First, it is very easy to observe that all such functions can't be distinguished by PDs, since they all share the PD in \Cref{fig:sim_1_PD}. 
Second, the interleaving distance between any two merge trees representing two functions $f_i$ and $f_j$ is $1/2$ if $i\neq j$ and $0$ otherwise. Thus the metric space obtained with $d_I$ from the data set $\{f_i\}_{i=0}^9$ is isometric to the discrete metric space on $5$ elements, where each point is on the radius $1$ sphere of any other point. 

We point out that there are applications in which it would be important to separate $f_0$ and $f_4$ more than $f_0$ and $f_1$, because they differ by \virgolette{an higher amount of edits}: for instance in \cite{cavinato2022imaging} merge trees are used to represent tumors, with leaves being the lesions, and it is well known in literature that the number of lesions is a non-negligible factor in assessing the severity of the illness \citep{ost2014prognostic}, and thus a metric more sensible to the cardinality of the trees is more suitable than $d_I$.

\begin{figure}

    \begin{subfigure}[c]{\textwidth}
    	\centering
    	\includegraphics[width = 0.45\textwidth]{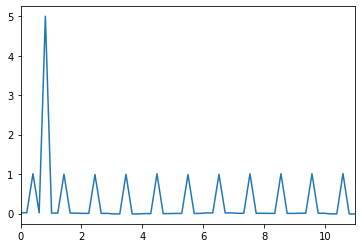}
	    \includegraphics[width = 0.45\textwidth]{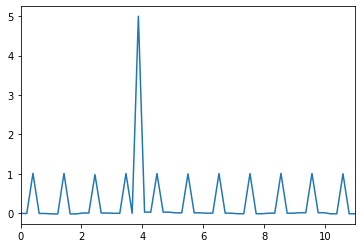}
    	\caption{The functions $f_0$, $f_3$ belonging to the simulated data set decribed in \Cref{sec:vs_interl}.}
    \label{fig:sim_1_fn}
	\end{subfigure}
		
    \begin{subfigure}[c]{\textwidth}
    	\centering
    	\includegraphics[width = 0.45\textwidth]{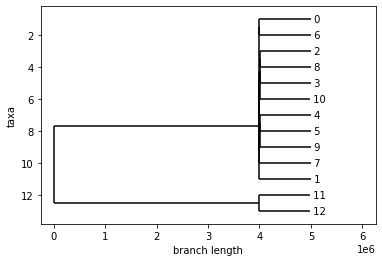}
	    \includegraphics[width = 0.45\textwidth]{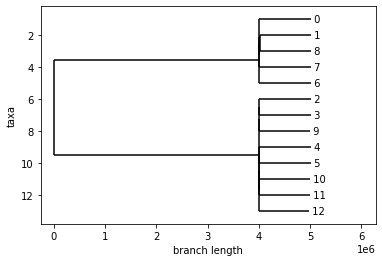}
    	\caption{The merge trees $(T_{f_0},h_{f_0})$ and $(T_{f_3},h_{f_3})$ associated to the functions in \Cref{fig:sim_1_fn}.}
   	\label{fig:sim_1_trees}
	\end{subfigure}
		
	\begin{subfigure}[c]{0.32\textwidth}
		\centering
		\includegraphics[width = \textwidth]{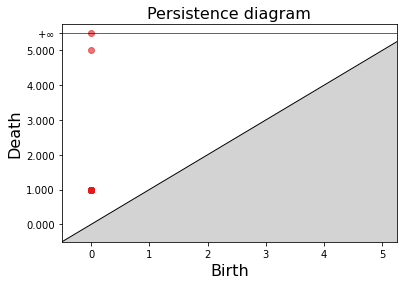}
		\caption{The persistence diagram representing the functions in \Cref{fig:sim_1_fn}. The point $(0,1)$ has multiplicity equal to the number of local minima minus $1$.}
		\label{fig:sim_1_PD}
	\end{subfigure}
	\begin{subfigure}[c]{0.32\textwidth}
    	\centering
    	\includegraphics[width = \textwidth]{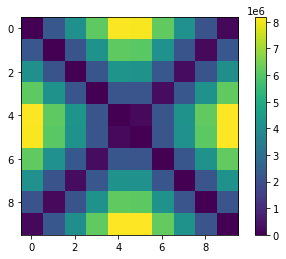}
    	\caption{Matrix of pairwise distances of the merge trees obtained from $\{f_i\}_{i=0}^{10}$.}
   	\label{fig:sim_1_dist}
    \end{subfigure}
    \begin{subfigure}[c]{0.32\textwidth}
    	\centering
    	\includegraphics[width = \textwidth]{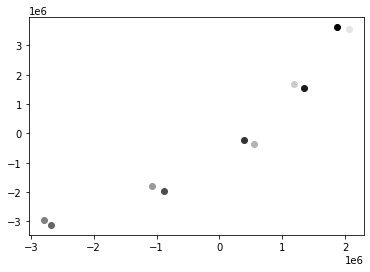}
    	\caption{Multidimensional Scaling Embedding in $\mathbb{R}^2$ of the matrix of pairwise distances shown in \Cref{fig:sim_1_dist}. The shades of gray describe, from white to black, the ordering of the trees.}
   	\label{fig:sim_1_MDS}
    \end{subfigure}

\caption{
Plots related to the simulated scenario presented in \Cref{sec:vs_interl}.}
\label{fig:ex_I}
\end{figure}

\subsection{Edit Distance Between Merge Trees \citep{merge_farlocca} and Wasserstein Distance \citep{merge_wass}}

The edit distance in \cite{merge_farlocca} is similar to classical edit distances, with the edit operations being restricted to insertion and deletion of vertices and with a \emph{relabeling} operation which is equivalent to our shrinking operation.
There is however the caveat that vertices are in fact understood as persistence pairs  $(m,s)$, with $m$ being the leaf representing the local minimum giving birth to the component and $s$ the internal vertex representing the saddle point where the components merge with an earlier born component and thus dies according to the elder rule. 
There is thus a one-to-one correspondence between persistence pairs in the merge tree and in the associated persistence diagram. Editing a vertex $m$ implies editing also its saddle point $s$: deleting $(m,s)$ means deleting all vertices $m'$ such that their persistent pair $(m',s')$ satisfies $s'<s$. If then $s$ becomes of order $2$, it is removed.
 In particular the authors highlight the impossibility to make any deletion - with the word \virgolette{deletion} to be understood according to our notation - on internal vertices, without deleting a portion of the subtree of the vertex. So they cannot delete and then insert edges to swap father-children relationships. 
To mitigate the effects of such issue they remove in a bottom-up fashion, as a preprocessing step, all saddles $s\in V_T$ such that $w_T((s,s'))<\varepsilon$  
for a certain positive threshold $\varepsilon$. All persistent pairs of the form $(m,s)$ are turned into $(m,s')$. Such issue is further discussed in \Cref{sec:preprocessing}.

Two merge trees are then matched via mappings representing these edit operations. To speed 
up the computations, the set of mappings between the trees is constrained so that disjoint subtrees are matched with disjoint subtrees.
The cost of the edit operation on an edit pair $(m,s)$ is equal to the edit operation being applied on the corresponding points in the associated PDs with the $1$-Wasserstein distance: deleting a persistent pair has the cost of matching the corresponding point to the diagonal and relabeling a persistent pair with another in the second tree has the cost of matching the two points of the two diagrams - see \cite{merge_farlocca} Section 4.3.1.

A closely related metric between merge trees is the Wasserstein distance defined in \cite{merge_wass}, which extends the metric by \cite{merge_farlocca} producing also further analysis on the resulting metric space of merge trees by addressing the problem of barycenters and geodesics. In this work the authors rely on a particular \emph{branch decomposition} of a merge tree, as defined in  \cite{merge_farlocca_2}, from which they induce the branch decomposition tree (BDT - \cite{merge_wass}, Section 2.3) used to encode the hierarchical relationships between persistence pairs. A branch decomposition is roughly a partition of the graph $T$ of a merge tree $(T,h_T)$ via ordered sequences of adjacent vertices \citep{merge_farlocca_2}. The chosen branch decomposition is the one induced by the elder rule and persistence pairs. Edit operations on such BDTs entail improved matchings and deletions between persistence pairs. To obtain the (squared) $2$-Wasserstein distance the vertices of two BDTs are matched and the resulting costs are squared and then added. 
However, the authors then explain that with this first definition geodesics cannot by found via 
linear interpolation of persistence pairs for the hierarchical structure of the merge tree can be broken. To mitigate that, they employ a normalization which shrinks all the branches on $[0,1]$, irrespective of their original persistence - \cite{merge_wass}, Section 4.2 - leading to simple geodesics obtained with linear interpolation between persistence pairs. To mitigate for this invasive procedure they introduce yet another preprocessing step artificially modifying small persistence features to reduce the normalization effects.

Some of the limitations of this approaches are listed in  
\cite{merge_wass}, Section 7.3. Also \cite{merge_farlocca_2}, Section 3.3 adds on that with further details and examples. 
 Namely, the restricted space of possible matchings between trees - which is key to obtain the computational performances of the metrics - forces unstable behaviours: issues with saddle swaps (see \cite{merge_wass}, Section 4.4 and Fig. 10)
 and instability of persistence pairs, so that elder ruled-based matchings may force very high persistence features to be matched with other very high persistence features even in situations where this implies making many unreasonable changes in the tree structures as in \Cref{fig:problem}
 (see also \cite{merge_farlocca_2}, Figure 1, Figure 2 b), Section 3.3, and \cite{merge_wass}, Section 7.3). Moreover, \cite{merge_wass} does not address the interactions between the normalization and the two preprocessing steps.

\subsection{Branch Decomposition-Independent Edit Distances for Merge Trees \citep{merge_farlocca_2}}

The work \cite{merge_farlocca_2} starts from the shortcomings of \cite{merge_farlocca} and \cite{merge_wass} trying to overcome them. Namely it defines branch decompositions, - for instance, the persitence pairs of \cite{merge_farlocca}, induce one such branch decomposition - and in order to avoid issues related with the instability of the persistence pairs, \citep{merge_farlocca_2} introduce also the possibility to optimize the chosen branch decomposition. The only big issue with such approach is that it does not define a metric on the space of merge trees, for the triangle inequality is not satisfied - see \cite{merge_farlocca_2} Theorem 2 and Figure 3.

\subsection{A Deformation-based Edit Distance for Merge Trees \citep{wetzels2022deformation}}

We consider the metric defined in \cite{wetzels2022deformation}.
The approach of \cite{wetzels2022deformation} shares many similarities with the starting point of ours, to the point that we believe that the metric $d_E$ between \emph{weighted} trees and the metric in \cite{wetzels2022deformation} should coincide in some cases. In fact, \cite{wetzels2022deformation} requires that, any time we obtain an order $2$ vertex via some deletions, that vertex is removed from the tree, with a ghosting-like procedure. However 1) they don't explicitly model mathematically the ghosting of an order $2$ vertex 2) the metric in \cite{wetzels2022deformation} is applied on merge trees without any truncation strategy and by simply turning merge trees into weighted trees by taking differences between the heights of the vertices and removing the edge at infinity (calling these objects \emph{abstract merge trees}, see Definition 1 in \cite{wetzels2022deformation}).

These fact have three consequences: 
\begin{enumerate}
    \item  it results in an unstable metric on merge trees, due to the issues presented in \Cref{fig:unpleasant};
    \item  it results in a metric which is completely insensitive to vertical translation, that is, it is not able to distinguish between the merge tree of $f$ and the merge tree of $f+h$, for any $h\in \R$; 
    \item even considering just weighted trees, we believe that explicitly modeling order $2$ vertices removal and addition is fundamental to study the geometric and topological properties of the space of merge trees, and designing data analysis tools. See \Cref{fig:refinement}.
\end{enumerate}

\begin{figure}
\centering
\includegraphics[width = 0.4\textwidth]{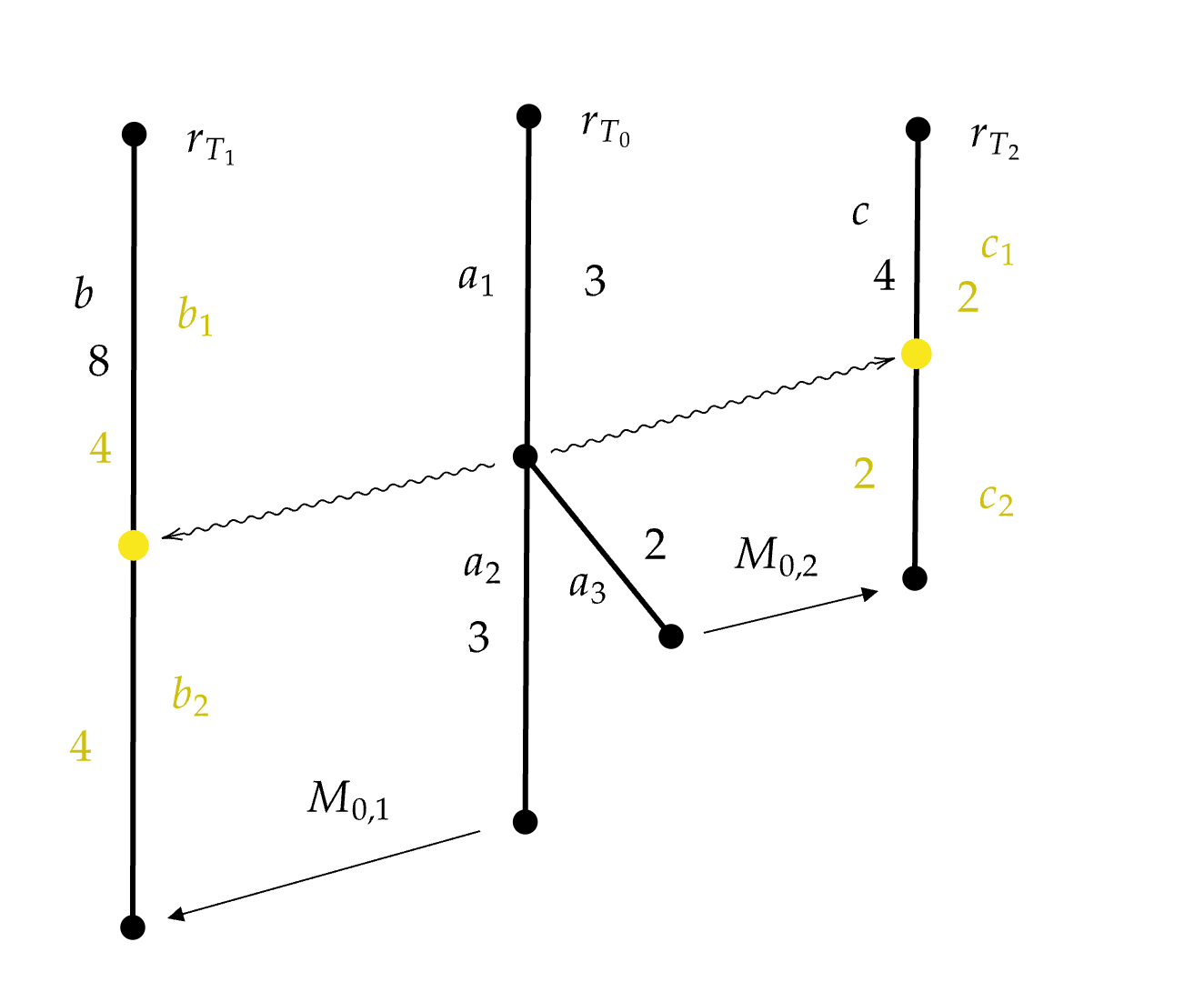}
\caption{In this figure we show the potential of explicitly modeling the addition and removal of order $2$ vertices. Consider the weighted trees $T_0$, $T_1$, and $T_2$ represented in the figure. Let $M_{0,i}$, with $i=1,2$ be two mappings $M_{0,i}$ defined, respectively, by the matchings given by the straight arrows (via the usual edge-vertex identification). $M_{0,1}$ is then completed by the deletion of $a_3$ and the ghosting of its father. And similarly, $M_{0,2}$ is then completed by the deletion of $a_2$ and the ghosting of its father. Suppose now that we split $T_1$ and $T_2$ inserting the yellow vertices and obtaining the yellow edges and the weighted trees $T'_1$, and $T'_2$. At this point we would be able to match $b_1$ with $a_1$ and $b_2$ with $a_2$ without raising the cost of $M_{0,1}$ and, similarly to match $c_1$ with $a_1$ and $c_2$ with $a_3$ without raising the cost of $M_{0,2}$. Thus we \virgolette{refined} our weighted trees and obtained mappings as expensive as the starting ones, but which involve only deletions and shrinkings, which are the operations giving the classical edit distance. Such distance is much more regular compared to $d_E$ (and to the one defined in \cite{wetzels2022deformation}). For instance, with the identifications $a_1\sim b_1 \sim c_1$, $a_2\sim b_2 $, and $a_3\sim  c_2$, we can represent $T_0$ with the vector $v_0 = (3,3,2)$, $T'_1$ with the vector $v_1 = (4,4,0)$ and $T'_2$ with the vector $v_0 = (2,0,2)$ and the mappings we described between those trees are straight lines in $\R^3$ with $\parallel\cdot \parallel_1$. This cannot be done with the edits in \cite{wetzels2022deformation}.}
\label{fig:refinement}
\end{figure}

\subsection{Heights vs Weights}
\label{sec:h_T_vs_w_T}

In this section we try to better understand the different behavior of $d_E$ when compared to the persistence-based metrics presented in the previous sections. The basic idea is that $w_T$ encodes the reciprocal positions of the merging points, instead of having the persistence pairs being free to move independently - at least locally - inside a constrained space.

Using $d_B(PD(T),PD(G))\leq d_I(T,G)$ \cite{merge_interl} and \Cref{prop:d_I<d_E}, one obtains:
\[
W_1(PD(T),PD(G))\leq (\rank(T)+\rank(G))d_E(T,G),
\]
with $PD(T)$ being the persistence diagram associated to the merge tree $T$. As we will see shortly, this bound cannot be improved.
It is thus evident that working with weights $w_T(e)$ instead of persistence pairs, as PDs do, creates differences in how the variability between trees is captured by these two metrics, despite the similarity in their stability properties.

\begin{figure}
\centering
	\begin{subfigure}[t]{0.49\textwidth}
    	\includegraphics[width = \textwidth]{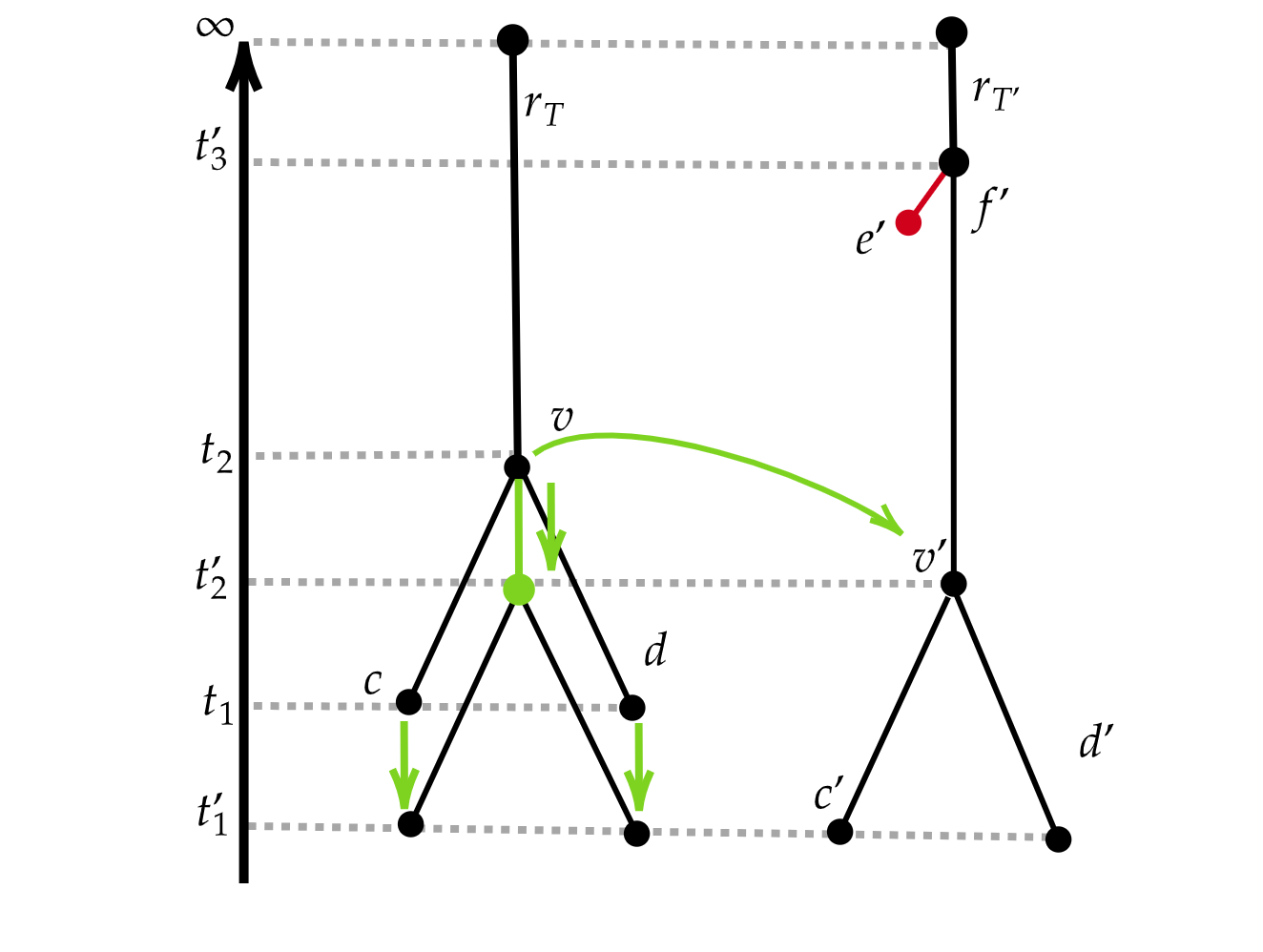}
    	\caption{Two merge trees $T$ (left) and $T'$ (right). We first delete $(e',f')$ and then edit $T$ with the shrinking defined by $(v,v')$.}
    	\label{fig:compare_1}
    \end{subfigure}
    \begin{subfigure}[t]{0.49\textwidth}
    	\centering
    	\includegraphics[width = \textwidth]{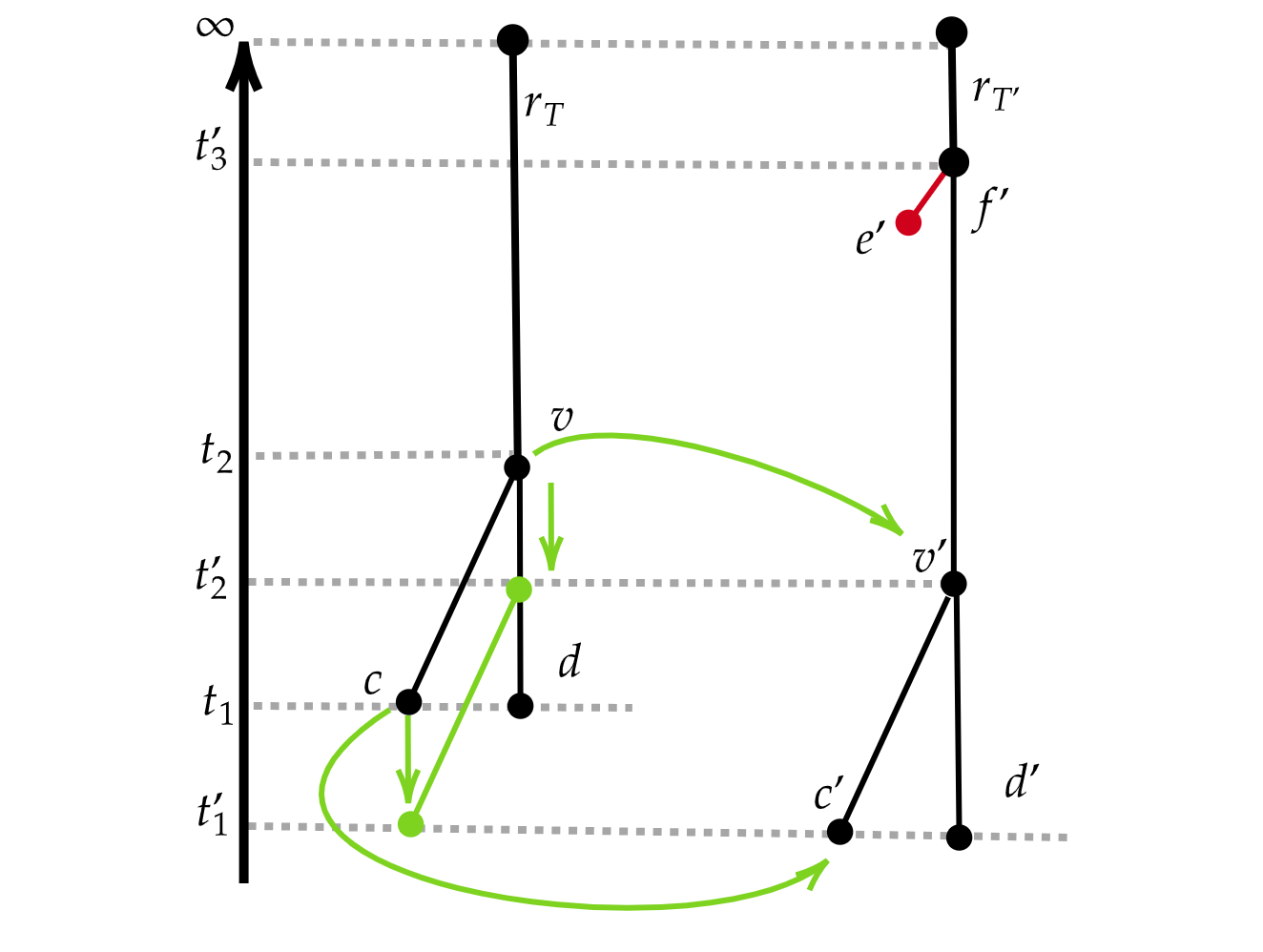}
    	\caption{The same merge trees as in \Cref{fig:compare_1} - with some edges drawn at different angles for visualization purposes. We edit $T$ matching the persistence pair $(c,v)$ with $(c',v')$ according to \cite{merge_farlocca, merge_wass}.}
    	\label{fig:compare_2}
    \end{subfigure}
    
    \centering
   	\begin{subfigure}[t]{0.49\textwidth}
    	\centering
    	\includegraphics[width = \textwidth]{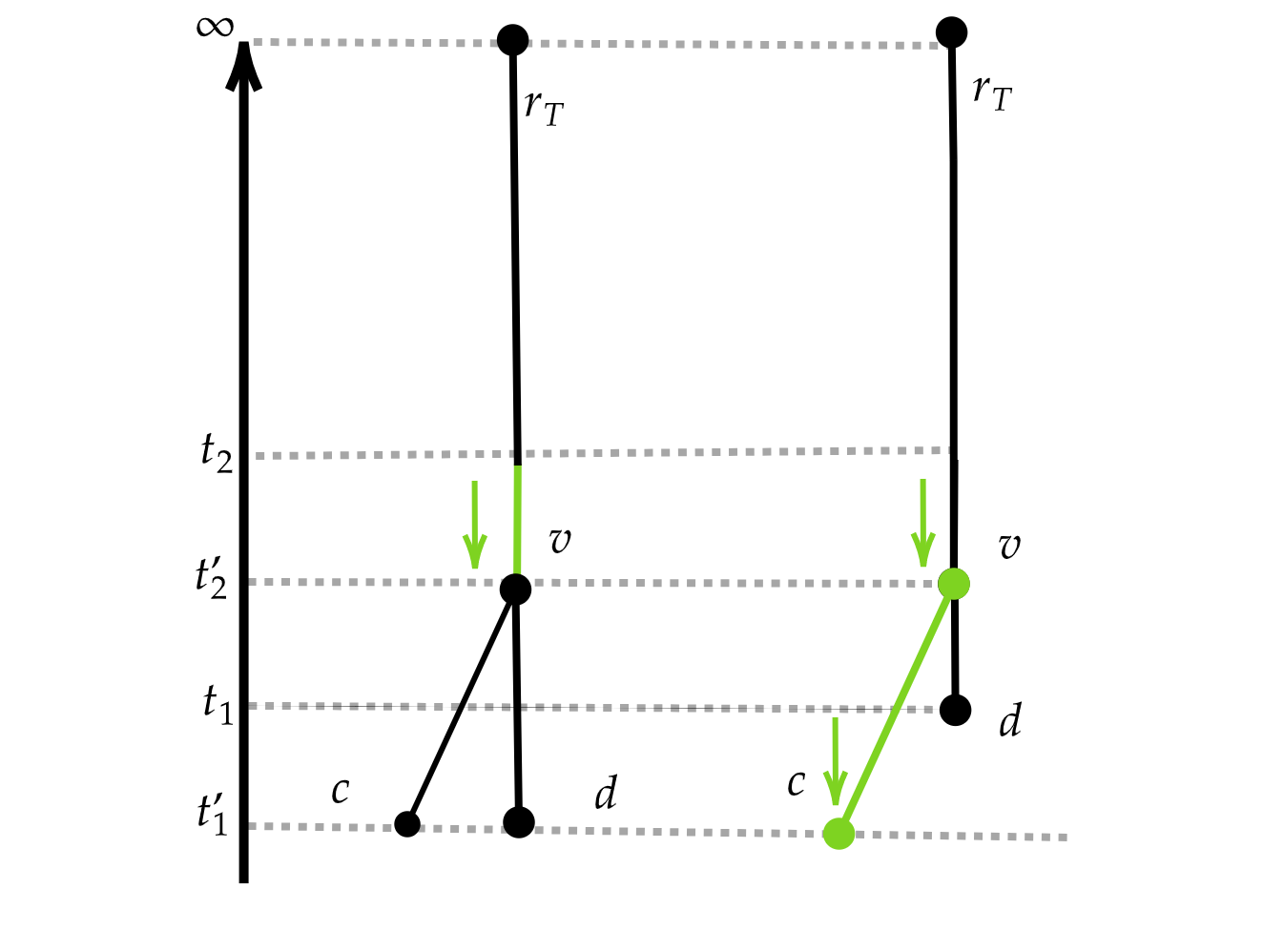}
    	\caption{The results of the edits applied on $T$ in \Cref{fig:compare_1} (left) and in \Cref{fig:compare_1} (right)}.
    	\label{fig:compare_recap}
    \end{subfigure}
\caption{A comparison between the weight based edits on which is based $d_E$, and the persistence based edits in \cite{merge_farlocca, merge_wass}.}
\label{fig:compare}
\end{figure}

Take, for instance,
the merge trees in \Cref{fig:unpleasant_1}: the persistence pairs are $(c,v)$ and $(d,r_T)$. The pairs of the rightmost tree instead are $(c',v'),(e',f')$ and $(d',r_{T'})$. Deleting $(e',f')$ according to \cite{merge_farlocca, merge_farlocca_2, merge_wass} 
amounts to deleting $e'$ according to our framework. Instead, shrinking $(v,r_T)$ on $(v',r_{T'})$, after the deletion of $(e',f')$, for us means lengthening $(v,r_T)$ by $t_2-t'_2$ and so lowering $v$ and all the vertices below $v$
by the same amount, as in \Cref{fig:compare_1}.
On the other hand, matching the persistence pair $(c,v)$ with $(c',v')$ for \cite{merge_farlocca, merge_wass} 
is equivalent to shifting the edge $(c,v)$ downwards towards $(c',v')$ leaving all other vertices fixed, as in \Cref{fig:compare_2}.
We can compare the results of such edits in \Cref{fig:compare_recap}.

Thus, given a persistence pair $(m,s)$ inside a merge tree, corresponding to the point $(b,d)$ inside the persistence diagram, moving $s$ upwards by some $\varepsilon>0$ such that $h_T(s)+\varepsilon<h_T(s')$ costs $\varepsilon$ in terms of the other edit distances, of the interleaving distance and in terms of the $1$-Wasserstein distance between PDs and leaves all other points with the same height. In terms of the edit distance $d_E$, instead, moving $s$ upwards by $\varepsilon>0$ shortens the edge above $s$ keeping all other edges of a fixed lengths and thus moves the vertices upwards by $\varepsilon$.

We can interpret this fact as the metrics in
\cite{merge_interl}, \cite{merge_farlocca} and \cite{merge_wass} 
better capturing similarities between trees via the location of their vertices in terms of heights; while $d_E$ better captures the \virgolette{shape} of the tree i.e. relative
positions of its vertices, being less sensible on the height variability since you can move more vertices at one. And this is exactly what happens in the example $f:I\rightarrow \mathbb{R}$ and $g(x)=f(x)+h$, $h\in\R$ - as in \Cref{fig:compare}: we have $W_1(PD(T),PD(G))=d_E(T,T')(\rank(T)+\rank(G))$.

\subsection{Stability vs Preprocessing}
\label{sec:preprocessing}

As already mentioned, the metrics in 
\cite{merge_farlocca} and \cite{merge_wass}
lack stability properties for similar reasons, which means that there are certain situations in which such metrics may perceive as very far trees which are in fact very close in terms of interleaving distance. In particular they are unable to model \virgolette{saddle swaps}, that is, with our terminology, deleting and inserting internal vertices to change father-children relationships. As noted also by the authors themselves, this issue needs to be addressed in some way. As already mentioned, to do so authors resort to a computational solution implemented as a preprocessing step: they fix a threshold $\varepsilon>0$ for any couple of persistence pairs $(m,s)$ and $(m',s')$ with $s'>s$ and $s'-s<\varepsilon$, they change $(m,s)$ into $(m,s')$ merging the two saddle points - in a bottom-up recursive fashion. 

In this section, we produce a brief investigation of such procedure, which is absent in the aforementioned works.
We represent the possible outcomes of this preprocessing in \Cref{fig:pers_pairs}. All merge trees in \Cref{fig:pers_pairs} are drawn with colors representing persistence pairs, and similar colors yields persistence pairs with the same persistence throughout the whole figure. They ideally should be matched by the metrics to achieve optimal distances as the differences between edges of different colors could be arbitrary big. 

In \Cref{fig:persistence_pairs_1} (left) we suppose that the edges marked with a red cross represent distances between saddles smaller or equal than $\varepsilon>0$, thus we recursively merge each saddle point with the higher one, starting from the bottom and going upwards. In this way we obtain the merge tree $T'$ in     
\Cref{fig:persistence_pairs_1} (right) which is then used as input for the metric.

In \Cref{fig:persistence_pairs_2} we see two merge trees $T$ and $G$ such that their interleaving distance is $\varepsilon+\varepsilon'$, for we need to move points of $T$ upwards by $\varepsilon+\varepsilon'$ to match persistence pairs of the same color in $G$. Their edit distance $d_E$ would be $3(\varepsilon+\varepsilon')$ for we need to delete and then insert back three small edges in $T$ to swap father-children relationships in a suitable way.
Note that one can replicate analogous situations to make the interleaving distance between the two merge trees high at will: informally speaking it is enough to add other persistence pairs as needed and force matchings between pairs in very \virgolette{different positions}.

In \Cref{fig:persistence_pairs_3} we see a possible output of the preprocessing routine.
If the preprocessing threshold is bigger than $\varepsilon$ and $\varepsilon'$ then the pre-processed trees $T'$ and $G'$ are equal. This in some sense is the desired output of the authors 
of  
\cite{merge_farlocca} and \cite{merge_wass} for now their metric can match the colors of the persistence pairs.
They suggest that such loss of variability -  $d(T',G')=0$ even if $T\nsim_2 G$ - could be mitigated by adding to the distance between the preprocessed trees the fixed threshold as many times as there are saddles merging in the procedure. Note that, if this artificial addition approximately matches the variability removed from the pairs attached to the orange vertical pair, it introduces artificial variability at the level of the branches attached to the brown edge, for they do not need any modification to be matched correctly.

In \Cref{fig:persistence_pairs_4}, instead, we suppose that the preprocessing threshold is    
smaller than $\varepsilon$ and $\varepsilon'$.
The metrics then are forced to match persistence pairs with different colors causing an excess of variability which is pictured with red edges in the figure. We point out that, depending on the weights of the persistence pairs involved, these edges could be of arbitrary length.

In \Cref{fig:persistence_pairs_5} we represent the last possible output of the preprocessing procedure, which is a situation in which the threshold we fix is greater that $\varepsilon'$ but smaller then $\varepsilon$. This is possibly the worst scenario: pairs are matched in an optimal way but we have introduced a lot of artificial variability, symbolized with red edges, in $G$. Again, because of additive phenomena caused by recursive saddle merging, these edges can be arbitrarily long.

Note that the preprocessing does not fix the issues presented in \Cref{fig:problem}, in \cite{merge_farlocca_2}, Figure 1, Figure 2 b), Section 3.3, and \cite{merge_wass}, Section 7.3.

\begin{figure}[!htbp]
	\begin{subfigure}[c]{0.49\textwidth}
    	\centering
    	\includegraphics[width = \textwidth]{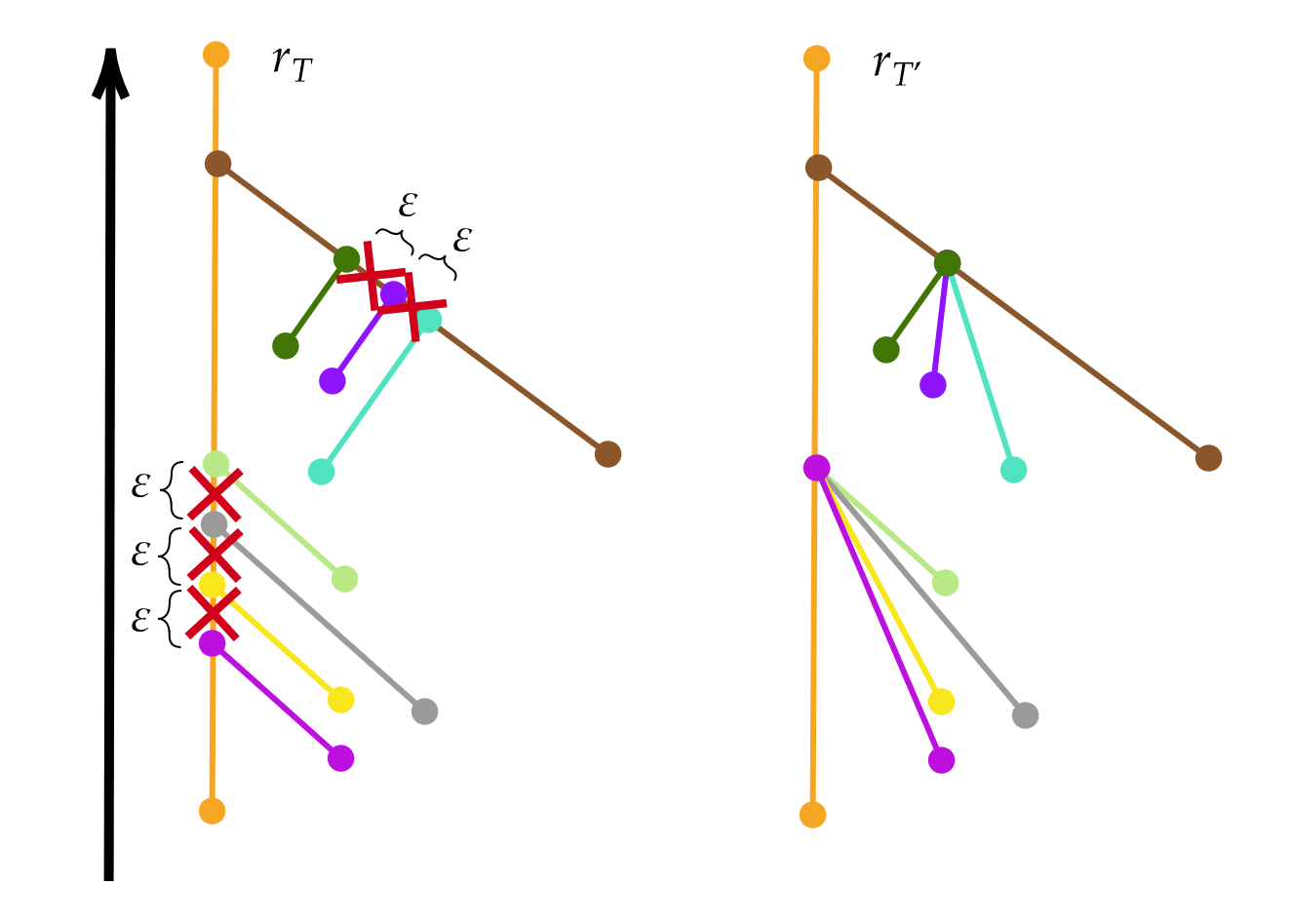}
    	\caption{A merge tree $T$ (left) with its pre-processed version (right) according to the procedure in \cite{merge_farlocca, merge_wass}.}
    	\label{fig:persistence_pairs_1}
    \end{subfigure}
   	\begin{subfigure}[c]{0.49\textwidth}
    	\centering
    	\includegraphics[width = \textwidth]{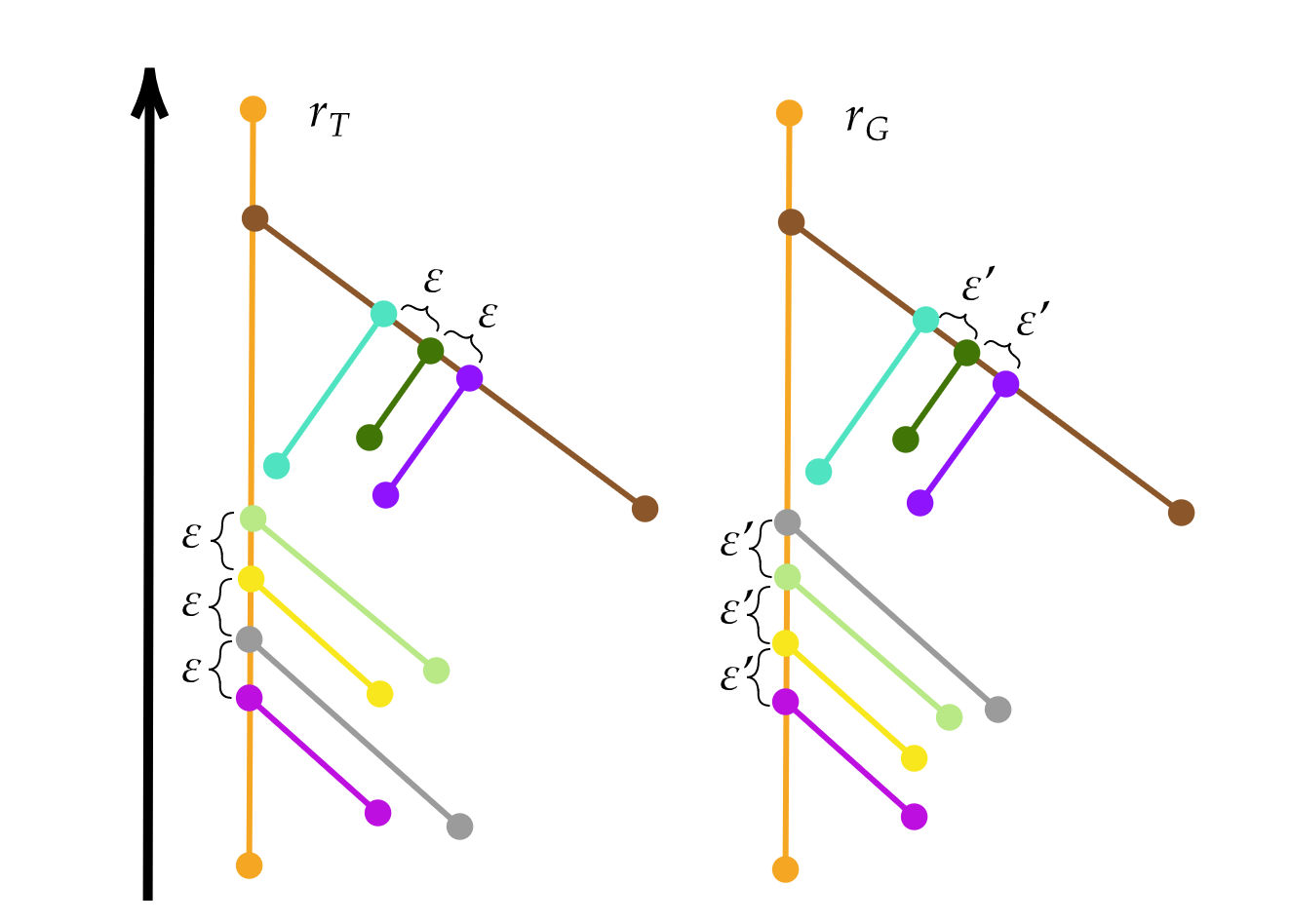}
    	\caption{Two merge trees $T$ and $G$ which are going to be preprocessed as in \Cref{sec:preprocessing}.}
    	\label{fig:persistence_pairs_2}
    \end{subfigure}

	\begin{subfigure}[c]{0.49\textwidth}
		\centering
		\includegraphics[width = \textwidth]{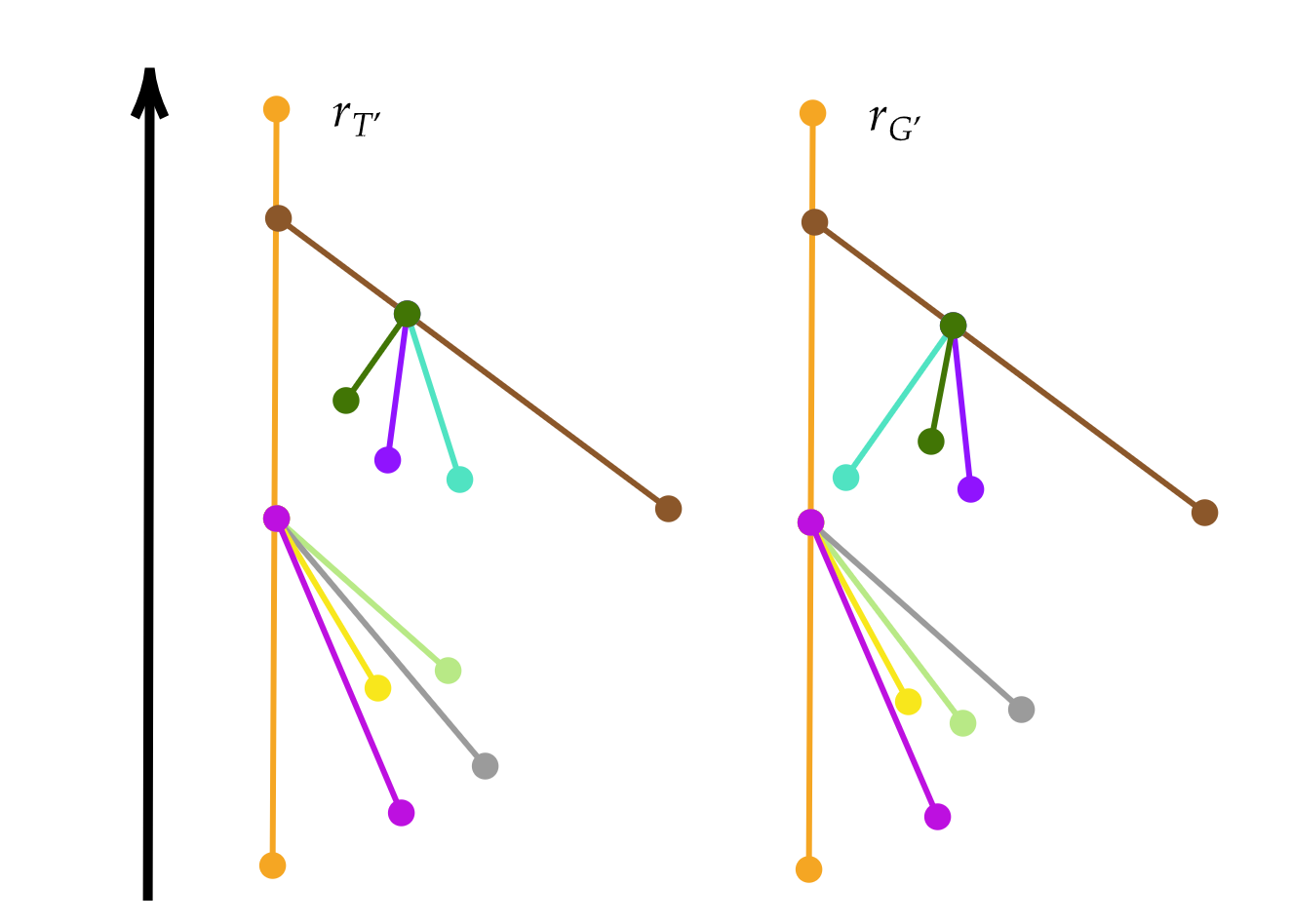}
		\caption{A possible output of the preprocessing: the preprocessing threshold is bigger than $\varepsilon$ and $\varepsilon'$.}
		\label{fig:persistence_pairs_3}
	\end{subfigure}
	\begin{subfigure}[c]{0.49\textwidth}
    	\centering
    	\includegraphics[width = \textwidth]{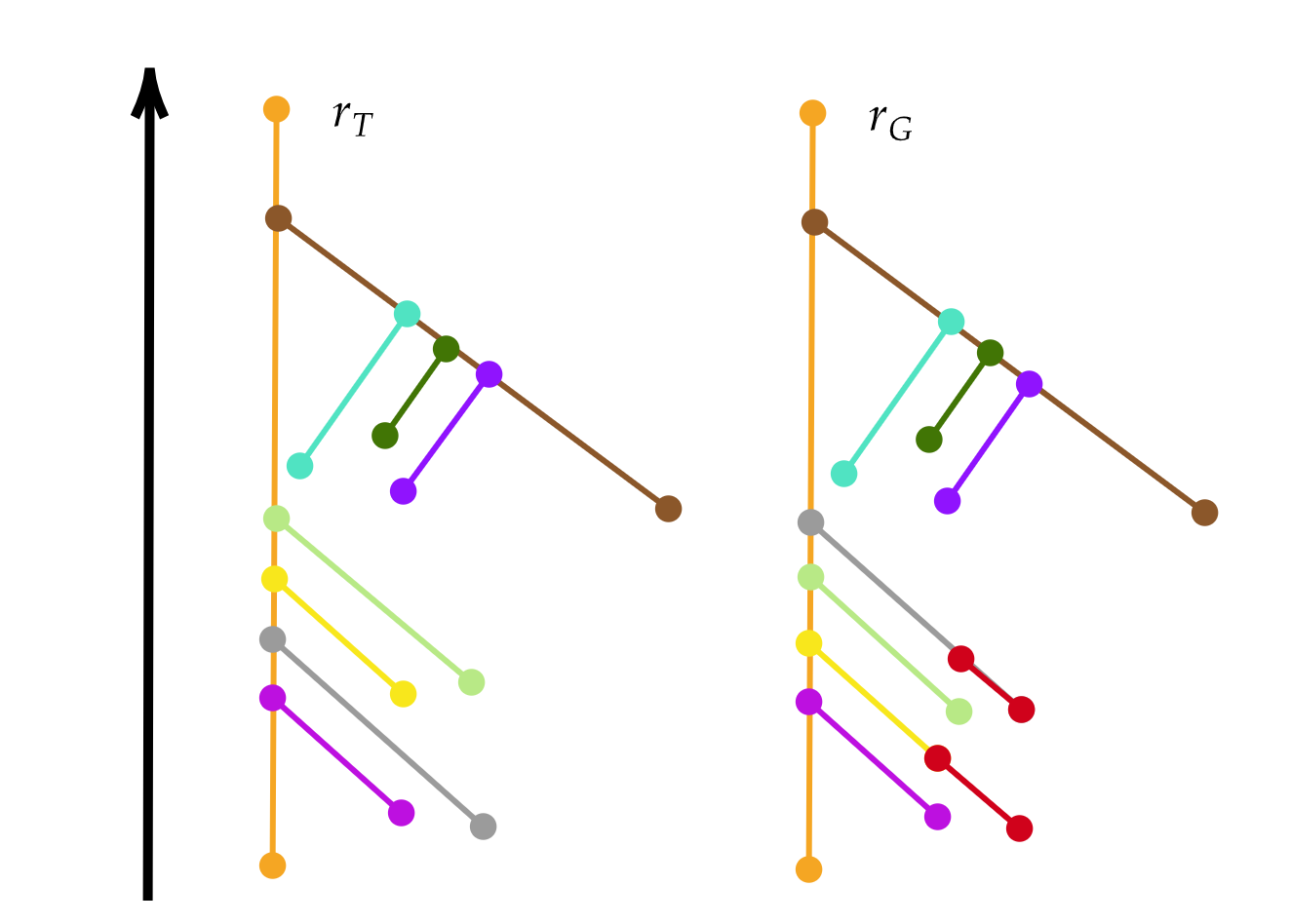}
    	\caption{Another possible output of the preprocessing: the preprocessing threshold is smaller than $\varepsilon$ and $\varepsilon'$.}
    	\label{fig:persistence_pairs_4}
    \end{subfigure}

    \begin{subfigure}[c]{0.49\textwidth}
    	\centering
    	\includegraphics[width = \textwidth]{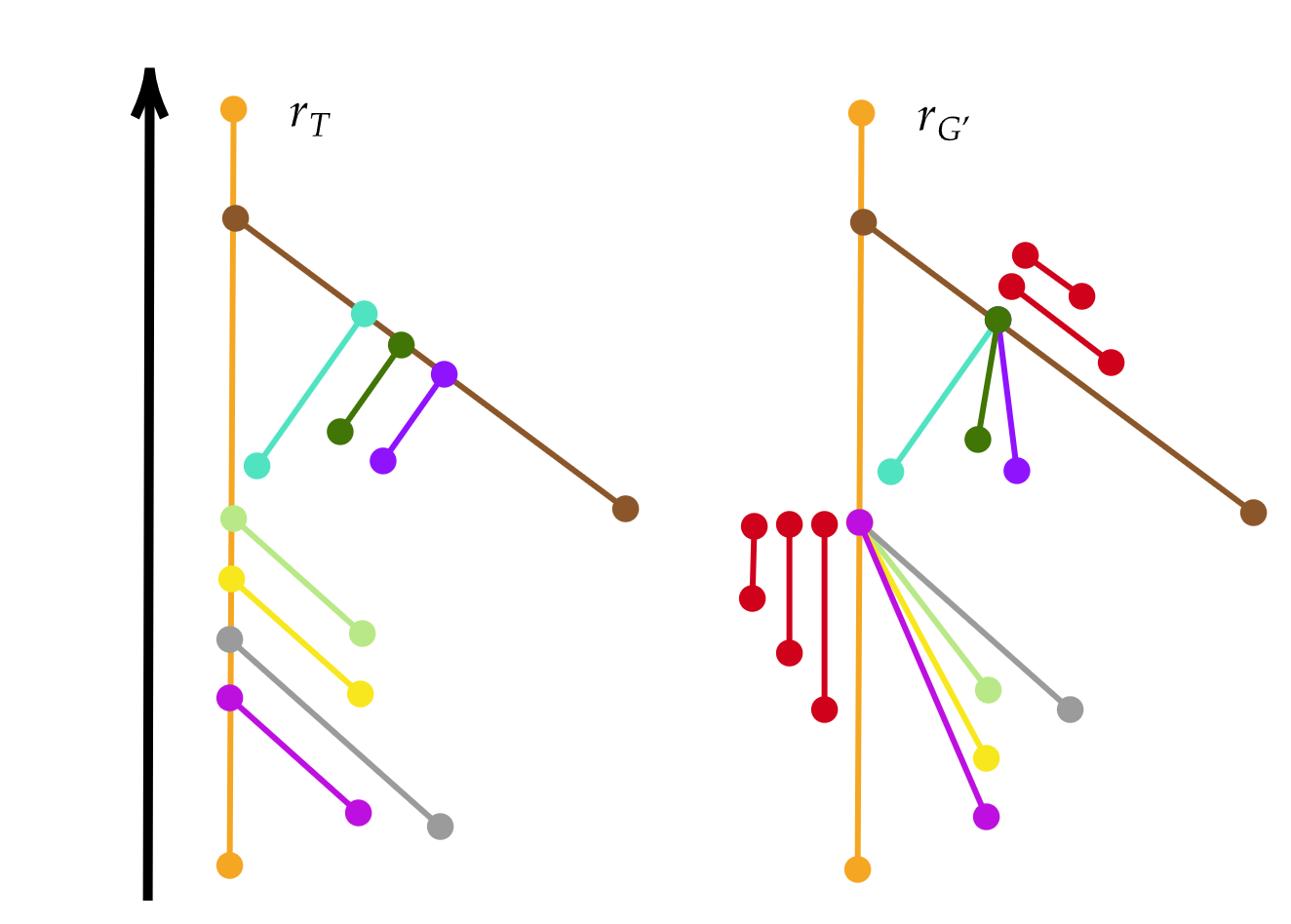}
    	\caption{The last possible output of the preprocessing: the preprocessing threshold is smaller than $\varepsilon$ and but bigger than $\varepsilon'$.}
    	\label{fig:persistence_pairs_5}
    \end{subfigure}
    \begin{subfigure}[c]{0.49\textwidth}
    	\centering
    	\includegraphics[width = \textwidth]{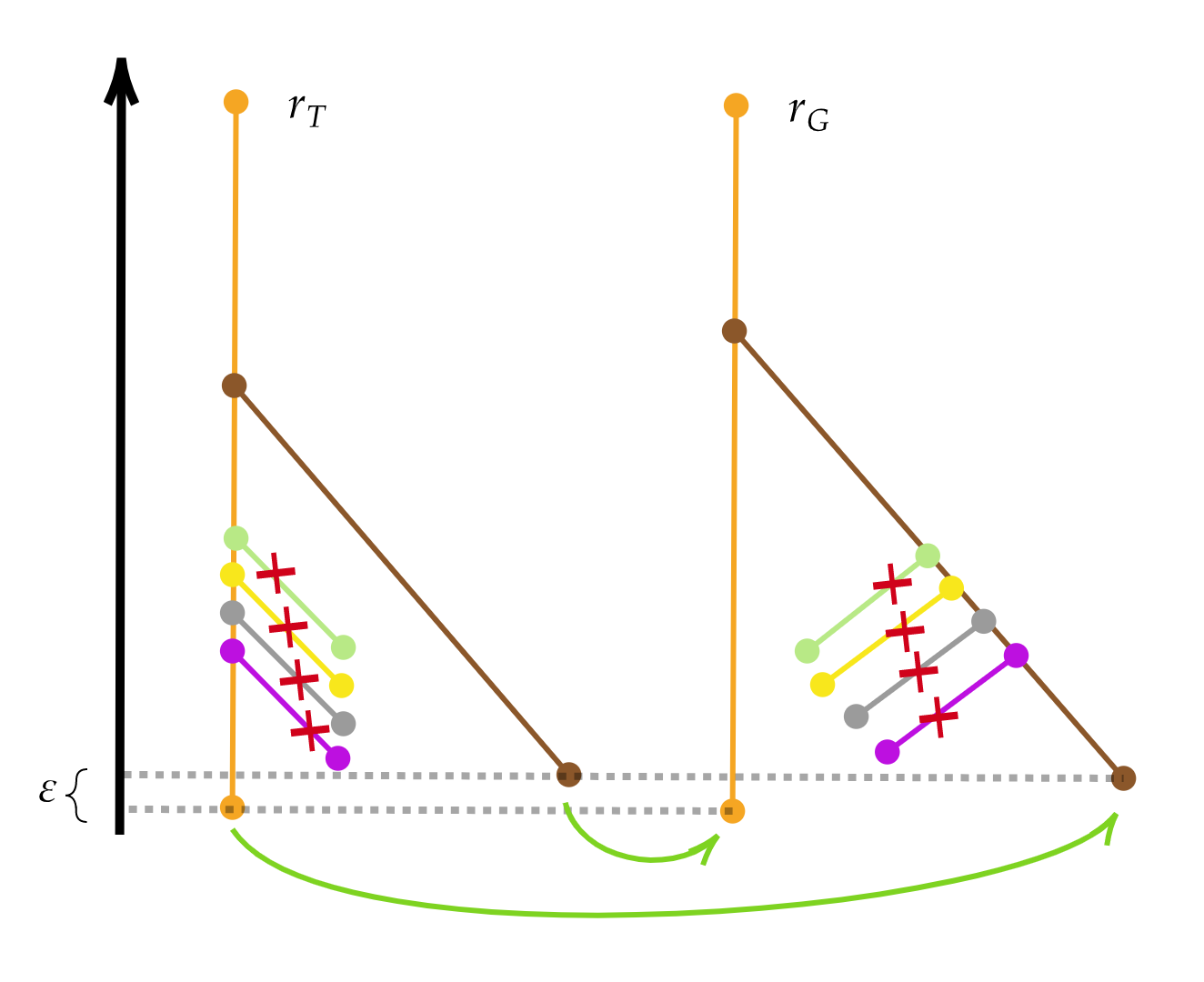}
    	\caption{As noted by \cite{merge_farlocca_2} and \cite{merge_wass}, metrics which are based on persistence pairs may force unnatural scenarios like deleting all the pairs with the red cross, in order to match the orange and the brown persistence pairs of $T$, respectively, with the orange and the brown of $G$. Instead, it would be much more natural to match the orange and the brown ones according to the green arrows and then matching the other pairs according to colors.}
    	\label{fig:problem}
    \end{subfigure}
    \hfill
\caption{Plots related to \Cref{sec:preprocessing} }
\label{fig:pers_pairs}
\end{figure}

To conclude, we point out the following fact. We deem fundamental to have metrics being able to measure different kinds of variability in a data set, especially if the information captured is interpretable, however, any attempt to adapt our edit distance $d_E$ to work with heights instead of weights, would face the problem of coherently define deletion edits. While it seems reasonable to have $w_T((v,v'))$ as a cost to delete and edge - for instance, this cost amounts to the persistence of the feature if $(v,v')$ is a persistence pair - 
any change in the height of $v'$ changes the weight of $(v,v')$ and thus the cost of its deletion, invalidating all the results about mappings in \cite{pegoraro2024finitelyfunc} and the algorithm therein. Similarly, adapting \cite{merge_farlocca} or  \cite{merge_wass} to handle saddle swaps, would mean at least removing from their mappings the property that would be (M3) in our notation, creating many issues from the theoretical and computational point of view. In fact, internal edges are not really available in the representation used by such metrics, where all points are instead understood as part of a persistence pair.

\section{Proofs}\label{sec:proofs_geom}

\bigskip\noindent
\underline{\textit{Proof of}  
\Cref{prop:truncation}.}

\smallskip\noindent

Consider a minimizing mapping $M \in M_2(G,G')$. 

Apply the deletions and ghostings described by $M$ both on $G$ and on $G'$ obtaining, respectively, the merge trees $G_M$ and $G'_M$. Note that, since $M \in M_2(G,G')$, 
if $v'$ is coupled with another edge by $M$, then $v$ must be deleted, otherwise $v'$ would be of order $2$ after the deletions, and thus it should be ghosted. Similarly, if $v'$ is deleted then $v$ must be deleted as well as: otherwise 
$v'$ would be of order $2$ after all the other deletions, and thus it cannot be deleted by $M$.
The same of course applies to $w'$ and $w$.

Let $w_{T'}(v')= n'$, $w_{T'}(v)= n$, $w_{G'}(w')= m'$  and $w_{G'}(w)=m$. By construction $n'=m'$ and thus the cost of the couple $(v',w')$, if it can be added a mapping, is zero. Since $M \in M_2(G,G')$, we know $(v',w')\notin M$. Our goal is to build a mapping containing $(v',w')$ with cost equal to the cost of $M$.

\begin{itemize}
    \item First suppose that $v$ is not deleted and $w$ is not deleted.
    
Then there exist $a = v_1 \leq \ldots \leq v_k \leq v'$ and $b = w_1 \leq \ldots \leq w_r \leq w'$ vertices in $T$ and $T'$ respectively such that:
\begin{itemize}
    \item $(a,b)\in M$;
    \item $v_i$ and $w_i$ are ghosted for all $i$ (as $M \in M_2(G,G')$ they cannot be deleted).
\end{itemize}
Let $n_i = w_T(v_i)$ and $m_i = w_{T'}(w_i)$ for all $i$.

We have seen that if $a=v'$, then $v$ must be deleted, which is absurd. Similarly we cannot have $b=w'$. 
Thus: $a= v_1 \leq \ldots \leq v_k = v <v'$
and $b= w_1 \leq \ldots \leq w_r = w<w'$. 
We have:
\[
\mid n' + \sum_{i=1}^{k} n_i - (m' + \sum_{i=1}^{r} m_i\mid = \mid n' - m' \mid + \mid \sum_{i=1}^{k} n_i-\sum_{i=1}^r m_i\mid.  
\]
This implies that we can add the couple $(v',w')$ to $M$ without increasing the cost of $M$.

\item Suppose now $v$ (and so $v'$) is deleted but $w$ is not. For what we have said before, it means that $w'$ is ghosted and $w$ is either coupled or deleted. 
This implies that the root of $G'_M$ is of order $1$. But then also the root of $G_M$ must be of order $1$. In other words, there exist $a = v_1 \leq \ldots \leq v_k < v$ and $b = w_1 \leq \ldots \leq w_r \leq w$ vertices in $T$ and $T'$ respectively such that:
\begin{itemize}
    \item $(a,b)\in M$;
    \item $v_i$ and $w_i$ are ghosted for all $i$ (as $M \in M_2(G,G')$ they cannot be deleted).
\end{itemize}

But this also implies that, after all the other deletions $v$ is of order $2$, and the same for $v'$. But this cannot happen.
\item Suppose lastly, $v$ and $w$ are all deleted. Which implies that also $v'$ and $w'$ are deleted. In this case we can add $(v,w)$ and $(v',w')$ to $M$ decreasing its cost.
\end{itemize}
Thus we can always add $(v',w')$ to $M$ and since the cost of such couple is zero, we have $d_E(T,T')=d_E(G,G')=d_E(sub_{G}(v'),sub_{G'}(w'))$.

\hfill$ \blacksquare $

\bigskip\noindent
\underline{\textit{Proof of}  
\Cref{teo:d_K}.}

\smallskip\noindent

Consider $(T,h_T)$ and $(T',h_{T'})$ such that 
$\max h_T< K'$ and 
$\max h_{T'}< K'$. 
Consider $(G,w_G)=\Tr_{K'}((T,h_{T}))$ 
and $(G',w_{G'})=\Tr_{K'}((T',h_{T'}))$. 
For any $K$ such that
$\max h_T< K<K'$ and 
$\max h_{T'}< K<K'$ there is a splitting of $(v,r_G)\in E_G$ with a vertex $w$ and a splitting of $(v',r_{G'})$ with a vertex $w'$ such that the obtained weighted trees $(R,w_R)$ and $(R',w_{R'})$ satisfy: $w_{R}((w,r_G))=w_{R'}((w',r_{G'}))$. Thus by \Cref{prop:truncation} we obtain $d_E(G,G')=d_E(R,R')=d_E(sub_R(w),sub_{R'}(w'))$. 
Moreover, $sub_R(w) = \Tr_K((T,h_{T}))$ and $sub_{R'}(w') = \Tr_K((T',h_{T'}))$.

Being $K,K'$ arbitrary we have that $d_K((T,h_T),(T',h_{T'}))$  does not depend on $K$, for $K> \max\{\max h_T,\max h_{T'}\}$. We need to prove the case $K=\max\{\max h_T,\max h_{T'}\}$.

Let $(G,w_G)=\Tr_{K}((T,h_{T}))$ and $(G',w_{G'})=\Tr_{K'}((T',h_{T'}))$, for $K'>K$. And consider now the weighted tree $\star=(\{\star\},w(\emptyset)=0)$,

We have that $d_E(G,\star)=\sum_{e\in E_{G}} w_{G}(e)$ and $d_E(G',\star)=\sum_{e\in E_{G'}} w_{G'}(e)$. Thus by the reversed triangle inequality:
\[
\mid d_E(G,\star)-d_E(G',\star)\mid=K'-K \leq d_E(G,G'),\] 
So we have: 
\[
d_E(\Tr_{K'}((T,h_{T})),\Tr_{K}((T,h_{T}))) = K'-K.\]
and we can set $K=\max\{\max h_T,\max h_{T'}\}$ and take $K'$ arbitrarily close to $K$ to finish the proof. 

\hfill$ \blacksquare $

\bigskip\noindent
\underline{\textit{Proof of}  
\Cref{prop:d_I<d_E}.}

\smallskip\noindent

Let $\T$ and $\G$ be two RAMTs which are represented by the merge trees $(T,h_T)$ and $(G,h_G)$, respectively. Then we obtain the weighted trees $(T,w_T)$ and $(G,w_G)$ via the some $\Tr_K$. The proof is quite long and requires some technical definitions. We split it into different parts to make it more easily readable.

\subparagraph*{Introduction - Display Posets}

For ease of notation, we introduce the following objects.

\begin{defi}[\cite{curry2021decorated}]Given a persistent set $S:\mathbb{R}\rightarrow \Sets$ we define its display poset as:
\[
D_S:=\bigcup_{t\in\mathbb{R}}S(t)\times\{t\}. 
\]
The set $D_S$ can be given a partial order with $(a,t)\leq (b,t')$ if $S(t\leq t')(a)=b$. 
\end{defi}

Let $D_{\T}$ and $D_{\G}$ be the display posets induced by $\T$ and $\G$. 
We call $h_{\T}$ the height function of $D_{\T}$ and $h_{\G}$ the height function of $D_{\G}$.
This construction is natural, in the sense that natural transformation between persistent sets become order preserving maps between the display posets. Given $\alpha: S\rightarrow S' $ natural transformation between persistent sets, then we have:
$D_\alpha: D_S \rightarrow D_{S'}$ defined by:
\[
D_\alpha((a,t))=(\alpha_t(a),t).
\]
Given $(a,t)\leq (b,t')$ we have that:
\[
(\alpha_t(a),t)\leq (\alpha_{t'}(b),t'),
\]
because $\alpha$ is a natural transformation.
For more details, see \cite{curry2021decorated}.

\subparagraph*{Introduction - Couplings}
Here we leverage on the notion of couplings between merge trees defined in 
\cite{pegoraro2024interleaving}.
Before recalling such definition that we highlight a subtle difference in merge trees as defined in  \cite{pegoraro2024interleaving} compared to the definition we give here. In \cite{pegoraro2024interleaving} the edge going to infinity which we require in our merge trees - $(v,r_T)$ with $h_T(r_T)=\infty$ - is not needed and thus such edge is removed. In other words, a merge tree $(T,h_T)$ is the sense of \cite{pegoraro2024interleaving} is such that $ord_T(r_T)>1$ and $h_T(r_T)\in \mathbb{R}$.  
We state the results in \cite{pegoraro2024interleaving} with the notation of the present manuscript.

Consider $C\subset V_T \times V_{G}$ and the projection $\pi_T:V_T \times V_{G} \rightarrow V_T$, we define the multivalued map $\Lambda_C^T:V_T\rightarrow V_T$ as follows: 
\[ 
\Lambda_C^T(v) =
\begin{cases} 
      \max_{v'<v} \pi_T(C) & \text{if }\#\{ v'\in V_T \mid v'<v \text{ and } v'\in \pi_T(C)\}>0 \\
      \emptyset & \text{otherwise.} 
   \end{cases}
\]
The term \virgolette{multivalued} means that $\Lambda_C^T$ is a function $\Lambda_C^T:V_T\rightarrow \mathcal{P}(V_T)$ with $\mathcal{P}(V_T)$ being the power set of $V_T$.

A coupling between the merge trees $(T,h_T)$ and $(G,h_G)$ - with the definition given in the present work - is then a set $C\subset V_T-\{r_T\} \times V_{G}-\{r_G\}$ such that:
\begin{itemize}
\item[(C1)] $\#\max C=1$ or, equivalently, $\#\max \pi_T(C)=\#\max \pi_G(C)=1$ 
\item[(C2)] the projection $\pi_T:C \rightarrow V_T$ is injective (the same for $G$)
\item[(C3)] given $(a,b)$ and $(c,d)$ in $C$, then $a<c$ if and only if $b<d$
\item[(C4)] $a\in \pi_T(C)$ implies $\#\Lambda_C^T(a)\neq 1$ (the same for $G$). 
\end{itemize}
Similarly to mappings each couple in $C$ is associated to a cost, and $\parallel C\parallel_\infty$ is defined as the highest of such costs. In \cite{pegoraro2024interleaving} it is proven that $d_I(T,G)\leq \parallel C\parallel_\infty$ for all couplings $C$.  

Define the following functions as in \cite{pegoraro2024interleaving}:

\begin{itemize}
\item define $\varphi_C^T:V_T\rightarrow V_T$ so that 
 $\varphi_C^T(x)=\min \{v\in V_T \mid v> x \text{ and }\#\Lambda(v)\neq 0\}$.
Note that since the set 
$\{v\in V_T\mid v> x\}$ is totally ordered, $\varphi_C^T(x)$ is well defined;
\item similarly, define $\delta_C^T:V_T\rightarrow V_T$, defined as $\delta_T(x) = \min\{v \in V_T \mid v\geq x \text{ and } v\in \pi_T(C)\}$;
\item set $\gamma^C_T:V_T-D_C^T\rightarrow V_G$ to be:
\[ 
\gamma^C_T(x) =
\begin{cases} 
      \arg\min\{g(w) \mid (v,w) \in C, v < x \} & \text{if }\#\{g(w) \mid (v,w) \in C, v < x \}>0 \\
      \emptyset & \text{otherwise.} 
   \end{cases}
\]
Note that if $\#\{g(w) \mid (v,w) \in C, v < x \}>0$, by \emph{(G)}, $\gamma^C_T(x)$  is uniquely defined. Moreover, $\gamma^C_T(\varphi^C_T(x))$ is well defined for any $v\in V_T$; 
\item lastly, set $\eta_C^T(x):=\gamma_C^T(\varphi_C^T(x))$.
\end{itemize}

To lighten the notation, when clear from the context, we may omit subscripts and superscripts.
The costs of the elements a coupling are defined in \cite{pegoraro2024interleaving} as follows.

\begin{itemize}
\item if $(x,y)\in C$, $cost_C(x)=|h_T(x)-h_G(y)|$;
\item if $x\notin \pi_T(C)$, we have two different scenarios:
	\begin{itemize}
	\item if $\#\Lambda(x)=0$, then  $cost_C(x)= \max\{\left( h_T(\varphi(x))-h_T(x)\right) /2,h_G(\eta(x))-f(x)\}$; 
	\item if $\#\Lambda(x)>1$, we have $cost_C(x)= \mid h_T(x)- h_G(w)\mid $ with $(\delta(x),w)\in C$;
	\item zero otherwise.
	\end{itemize}
\end{itemize}

Lastly, via Remark 3, Theorem 1, and Theorem 3 \cite{pegoraro2024interleaving} shows that,  $C$ induces 
\[
\alpha_C:D_{\T}\rightarrow D_{\G}\text{ and }\beta_C:D_{\G}\rightarrow D_{\T}\] 
such that 
\[
cost_C(v)=\mid h_{\G}(\alpha_C(v))-h_{\T}(v)\mid
\]
 for any $v\in V_T$ 
and 
\[
cost_C(w)=\mid h_{\T}(\beta_C(w))-h_{\G}(w)\mid
\] 
for any $w\in V_G$.

\subparagraph*{Main body of the proof}

We now want to establish relations between mappings between $\Tr_K((T,h_T))$ and $\Tr_K((G,h_G))$ - with $K>\max h_T,\max h_G$ - and couplings between $(T,h_T)$ and $(G,h_G)$.
When working with $\Tr_K((T,h_T))$ and $\Tr_K((G,h_G))$ we can apply the following corollary of \Cref{prop:truncation}.

\begin{cor}\label{cor:max_mapp}
Given $(T,w_T)$ and  $(T',w_{T'})$ weighted trees, if $r_T$ and $r_{T'}$ are of order $1$, there exist a minimizing mappings $M$ with $\# \max M = 1$.
\end{cor}

Consider a minimizing mapping $M\in M_2(T,T')$ with $\# \max M = 1$. Thanks to \Cref{lemma:M_2}, we have that $a\in \pi_T(M)$ implies $\#\Lambda_M^T(a)\neq 1$ (the same for $G$). In fact (C4) is equivalent to having no vertices of order two after deletions and ghostings. This means that $C_M:=\{(a,b)\in M \mid a\in V_T \text{ and }b \in V_G\}$ is a coupling. 

Conversely, given a minimizing coupling $C$, the set 
\begin{align*}
M_C:=C&\bigcup\{(a,"D")\mid a\notin\pi_T(C)\text{ and } \#\Lambda_C^T(a)\neq 1\}\\
&\bigcup\{("D",b)\mid
 b\notin\pi_G(C)\text{ and } \#\Lambda_C^G(b)\neq 1 \}\\
&\bigcup\{(a,"G")\mid a\notin\pi_T(C)\text{ and } \#\Lambda_C^T(a)= 1 \}\\
&\bigcup\{("G",b)\mid b\notin\pi_G(C)\text{ and } \#\Lambda_C^G(b)= 1\}
\end{align*}
is a mapping. Clearly $M_{C_M}=M$ and $C_{M_C}=C$.

Now we prove the following lemma.

\begin{lem}\label{lem:diseq}
Let $f:[a,K]\rightarrow [b,K]$ be a monotone function such that $f(K)=K$, $K\in\mathbb{R}$. For every $\{x_1<\ldots<x_{n+1}=K\}\subset [a,K]$:
\begin{equation}\label{eq:diseq}
\max_{i=1,\ldots,n+1} d(f(x_i),x_i) \leq \sum_{i=1}^n \mid d(x_i,x_{i+1}) - d(f(x_i),f(x_{i+1}))\mid 
\end{equation}
 
\begin{proof}
Let $v_i = f(x_i)-x_i$. And let $m\in \{1,\ldots,n+1\}$ be such that $\mid f(x_m)-x_m\mid = \max \mid f(x_i)-x_i\mid$.

Then we set:
\begin{align*}
v_1 =& v_m + \varepsilon_1 \\
v_2 =& v_m + \varepsilon_1 + \varepsilon_2\\
\ldots\\
v_i =& v_m + \sum_{j=1}^i\varepsilon_j\\
\ldots\\
v_n =& v_m + \sum_{j=1}^n\varepsilon_j.
\end{align*}

Thus we can write \Cref{eq:diseq}, which we need to prove, as:
\[
\mid v_m\mid \leq \sum_i\mid x_i-f(x_i)-(x_{i+1}-f(x_{i+1}))\mid = \sum_i\mid v_i-v_{i+1}\mid = \sum_{i=1}^n\mid \varepsilon_{i+1}\mid.
\]
Clearly we have some constraints on $\varepsilon_i$, in fact: 
\[
x_i + v_m + \sum_{j=1}^i\varepsilon_j \leq
x_{i+1} + v_m + \sum_{j=1}^{i+1}\varepsilon_j 
\]
which means:
\[
x_i -x_{i+1}\leq \varepsilon_{i+1}. 
\]

In particular $x_n+v_n=x_n=1$ means that $v_m + \sum_{j=1}^n\varepsilon_j=0$ i.e. $-v_m=\sum_{j=1}^n\varepsilon_j
$. Thus:
\[
\mid v_m\mid =\mid\sum_{j=1}^n\varepsilon_j\mid
\leq \sum_{i=1}^n\mid \varepsilon_{i+1}\mid.
\] 
\end{proof}
\end{lem}

Consider now a leaf $l\in L_T$ and take $\zeta_l=\{p\in V_T\mid v\geq l\}$. Consider the interval $[h_T(l),K]$. The map $v\in \zeta_l \mapsto h_T(v)$ gives a $1:1$ correspondence between $\zeta_l$ and a finite collection of points in $[h_T(l),K]$. 
We define $f(h_T(v))=h_{\G}(\alpha_C(v))$ for $v\in \zeta_l-\{r_T\}$. And $f(K)=K$. Note that $\alpha_C$ is constructed such that $\alpha_C(v)=w$ if $(v,w)\in C$. 
Thus $f(h_T(v))\leq K$ for $v\in \zeta_l$.

We extend $f$ on  $[h_T(l),K]$ via linear interpolation.
Since $\alpha_C$ is monotone wrt the partial order on $D_{\T}$, then $f$ is monotone on $[h_T(l),K]$.

Now consider $T$ and apply on it all the deletions and ghostings which involve points which are not in $\zeta_l$.  
We call $\{v_1=l<\ldots<v_{n+1}=r_T\}$ the coupled points in $\zeta_l$. 
Then: 
\begin{equation}\label{eq:aux_1}
\mid h_T(v_i)-h_T(v_{i+1})-f(h_T(v_i))-f(h_T(v_{i+1}))\mid \leq cost_{M_C}([v_i,v_{i+1}]\mapsto [\alpha_C(v_i),\alpha_C(v_{i+1})])
\end{equation}
With $cost_{M_C}([v_i,v_{i+1}]\mapsto [\alpha_C(v_i),\alpha_C(v_{i+1})])$ being the cost of the edits associated to points $v_i\leq p\leq v_{i+1}$ and  $\alpha_C(v_i)\leq q\leq \alpha_C(v_{i+1})$. In fact we have equality in \Cref{eq:aux_1} if there are no deletions for any $p\in V_T$ such that
$v_i\leq p\leq v_{i+1}$ and any $q\in V_{T'}$ such that $\alpha_C(v_i)\leq q\leq \alpha_C(v_{i+1})$. Otherwise the total cost of the deletions and shrinking exceeds 
$\mid h_T(v_i)-h_T(v_{i+1})-f(h_T(v_i))-f(h_T(v_{i+1}))\mid$ since:
\[
\mid n_1 + n_2 - (n_3+n_4) \mid \leq n_1 + n_3 + \mid n_2-n_4\mid  
\]
 for any $n_1,\ldots,n_4\in\mathbb{R}_{\geq 0}$.
 Lastly note that
 $\mid h_T(v_i)-f(h_T(v_i))\mid = cost_C(v)$.

By \Cref{lem:diseq}, we have:
\begin{equation}\label{eq:eq_finale}
\begin{aligned}
\max cost_C(v_i)=&\max h_T(v_i)-f(h_T(v_i))\\
 \leq &
\sum_{i=1}^{n}\mid h_T(v_i)-f(h_T(v_i))+f(h_T(v_{i+1})- h_T(v_{i+1})\mid \\
 \leq & cost_{M_C}([l,r_T]\mapsto [w,r_G])
\end{aligned} 
\end{equation}
with $(v,w)\in M_C$ and $v=\min \zeta_l \cap \pi_T(C)$.

\subparagraph*{Conclusion}
Let $\parallel C\parallel_\infty = cost_C(x)$. WLOG $x \in V_T$. Then $x\in \zeta_l$ for some $l\in L_T$. By applying \Cref{eq:eq_finale} we obtain the result.

\hfill$ \blacksquare $

\bigskip\noindent
\underline{\textit{Proof of}  
\Cref{prop:locally}.}

\smallskip\noindent

    In the following proof, we employ display posets, which are introduced in the proof of
    \Cref{prop:d_I<d_E}.  

    Consider $\mathcal{M}(\T)$ and $\pi_0(\X^n)$ such that $T=\mathcal{M}(\T)$ and $T_n=\mathcal{M}(\pi_0(\X^n))$. Note that we may often deliberately confuse
    merge trees and RAMTs to lighten the notation.

    Take any $\varepsilon>0$ and then take $N$ such that $d_I(T,T_n)\leq \varepsilon$ for every $n>N$.
    Let $G_n:=\S_{2\varepsilon}(T_n)$.
    Clearly, we always have $d_I(T_n,G_n)=2\varepsilon$. 
    Suppose that $\sup_{n\in N} \dim(G_n) = +\infty$.
    If this is the case, for every $K>0$, there is a merge tree $T_n$, with more than $K$ leaves $v$ such that $i^{2\varepsilon}_{\pi_0(\X^n)}(v)$ (more precisely,
    $(i^{2\varepsilon}_{\pi_0(\X^n)})_{h_{T_n}(v)}(v)$)
    is a leaf for $G_n$ and thus $h_{T_n}(father(v))-h_{T_n}(v)> 2\varepsilon$. 
    We now show that if $K$ is bigger than the number of leaves of $T$, we obtain a contradiction as, in order to have $d_I(T,T_n)\leq \varepsilon$, we cannot contract edges which are longer than $2\varepsilon$. 
    
    More formally,
    let $\alpha,\beta$ be $\varepsilon$-compatible maps between $T$ and $T_n$. We know they exists as $d_I(T,T_n)\leq \varepsilon$. These natural transformations, as mentioned in 
    the proof of
    \Cref{prop:d_I<d_E}, induce order preserving maps between display posets. In particular:
    \begin{align*}
    &D_\alpha:D_{\T}\rightarrow D_{\S_{\varepsilon}(\pi_0(\X^n))},\\
    &D_\beta:D_{\pi_0(\X^n)}\rightarrow D_{\S_{\varepsilon}(\T)}.
    \end{align*}
    Note that $G_n =\S_{2\varepsilon}(\pi_0(\X^n)) = \S_{\varepsilon}(\S_{\varepsilon}(\pi_0(\X^n)))$ and, in particular, $\S_{\varepsilon}(\pi_0(\X^n))$ cannot have less leaves than $G_n$.
    
    Since $T$ has less leaves that $G_n$ and so of $\S_{\varepsilon}(\pi_0(\X^n))$, we claim that there exist $v,w$ leaves of $T_n$, with $t=h_{T_n}(w)$ and $t'=h_{T_n}(v)$, such that:
    \begin{enumerate}
        \item $(i^{\varepsilon}_{\pi_0(\X^n)})_t(w)$ and $(i^{\varepsilon}_{\pi_0(\X^n)})_{t'}(v)$ are leaves of $\S_{\varepsilon}(\pi_0(\X^n))$, and $(i^{2\varepsilon}_{\pi_0(\X^n)})_t(w)$ and $(i^{2\varepsilon}_{\pi_0(\X^n)})_{t'}(v)$ are leaves of $G_n$;
        \item $h_{T_n}(father(v))-t'> 2\varepsilon$ and $h_{T_n}(father(w))-t> 2\varepsilon$;
        \item $D_\beta((w,t))\leq D_\beta((v,t'))$.
    \end{enumerate}

    We already know that 1. and 2. need to be satisfied for some leaves. 
    Not being able to find two such leaves which also satisfy 3. would imply $T$ having at least $K$ leaves.
        
    Since $D_\beta((w,t))\leq D_\beta((v,t'))$, we have
    $D_\alpha(D_\beta((w,t))\leq D_\alpha(D_\beta((v,t'))$. By definition, we know that 
    $D_\alpha(D_\beta((v,t'))=(X_{t'\leq t'+2\varepsilon}(v),t'+2\varepsilon)$. In particular, this means that  
    $(v,t'),(w,t')\leq (X_{t'\leq t'+2\varepsilon}(v),t'+2\varepsilon)$. This is absurd because we have found a point on the abstract merge tree (associated to) $T_n$, which is a common ancestor of both $v$ and $w$ and whose height differs from the ones of $v$ and $w$ by exactly $2\varepsilon$. But we know that the edges $(v,father(v))$ and $(w,father(w))$ are both longer than $2\varepsilon$.

    Thus $\sup_{n>N} \dim(G_n)=K$ for some $K\in \R$ such that $K \leq \dim(T)$. WLOG we can suppose $K=\dim(T)$. Since $d$ is finitely stable, for some $C>0$ we have:
    \begin{align*}
    d(T,G_n)\leq C\cdot 2K d_I(T,G_n) &\leq C\cdot 2K (d_I(T,T_n)+d_I(T_n,G_n))\\
    &=C\cdot 2K (d_I(T,T_n)+\varepsilon).        
    \end{align*}
    Thus by choosing $\varepsilon$ such that:
    \[
    \max\{2\varepsilon, 2KC\varepsilon\}\leq \epsilon,
    \]
 we obtain the thesis.

\hfill$ \blacksquare $

}

\bibliography{references}

\end{document}